\input amstex
\documentstyle{amsppt}
\magnification=\magstep1
\vsize =21 true cm 
\hsize =16 true cm 
\loadmsbm
\topmatter
\centerline{\bf Trace formulas for a class of 
                compact complex surfaces}
\author{\smc Lei Yang}\endauthor
\endtopmatter
\document

\centerline{\bf Contents}
$$\aligned
\text{1.} \hskip 0.1 cm 
          &\text{Introduction}\\
          &\text{1.1. Notation and terminology}\\
          &\text{1.2. The main results}\\
\text{2.} \hskip 0.1 cm
          &\text{The spectral theory of $L^2(\Gamma \backslash
                {\frak S}_{2}, k)$}\\
          &\text{2.1. The construction of hyperbolic elements}\\
          &\text{2.2. The basic integral operator of weight $0$}\\
          &\text{2.3. The basic integral operator of weight $k$}\\
          &\text{2.4. The spectral theory of $L^2(\Gamma \backslash 
                      {\frak S}_{2}, k)$}\\
\text{3.} \hskip 0.1 cm
          &\text{Trace formulas and zeta functions of weight $k$}\\
          &\text{3.1. Basic notions}\\
          &\text{3.2. The computation of the identity}\\ 
          &\text{3.3. The computation of hyperbolic elements}\\
          &\text{3.4. The trace formulas of weight $k$}\\
          &\text{3.5. The zeta functions of weight $k$}\\
          &\text{3.6. The functional equations of zeta functions}\\
          &\text{3.7. The analogue of the Riemann hypothesis}
\endaligned$$

\vskip 0.5 cm
\centerline{\bf 1. Introduction}
\vskip 0.5 cm

  In his epoch-making paper \cite{Se}, Selberg developed a 
trace formula for discrete subgroups of $SL(2, {\Bbb R})$ and 
proved the analogue of the Riemann hypothesis on compact 
Riemann surfaces of genus $g \geq 2$. There are at least 
two ways to generalize the Selberg trace formula. The first 
way involves passing from the upper half plane ${\Bbb H}$ 
to $SL(2, {\Bbb R})$ or its universal covering space. Duflo 
and Labesse \cite{DL} studied the trace formula on 
$GL(2)$. Later on, Arthur \cite{Ar} developed the 
trace formula for semisimple Lie groups of rank one. On the 
other hand, Gangolli and Warner \cite{Ga}, \cite{GW} 
gave the zeta functions of Selberg's type for symmetric 
spaces of rank one. The other approach concerns with the 
index theorem (see Millson \cite{M}, Singer \cite{S}, 
Barbasch and Moscovici \cite{BM} and Stern \cite{St}). 
For example, Millson \cite{M} gave the functional equation 
of the zeta function on a compact oriented $4n-1$ dimensional 
Riemannian manifold of constant negative curvature and found 
a connection between the zeta function and the $\eta$-invariant 
introduced by Atiyah, Patodi and Singer \cite{APS}. The 
second way is based on some interesting results of Maass 
and Roelcke. Maass \cite{Ma} studied some very basic partial 
differential operators connected with nonholomorphic 
automorphic forms on $\Gamma \backslash {\Bbb H}$.
Roelcke \cite{Ro} gave the spectral theory of these
operators (see also the work of Elstrodt \cite{El}). 
Hejhal \cite{H} obtained the trace formula for the 
automorphic forms of weight $m$ associated to these 
operators. Recently, D'Hoker and Phong \cite{DP}, 
Sarnak \cite{Sa}, Voros \cite{V}, Cartier and Voros 
\cite{CV} and Kurokawa \cite{Kr} calculated the 
regularized determinants of these operators, which are 
crucial in the string theory, and also come up in 
Ray-Singer's analytic torsion \cite{RS}.  

  In the present paper, motivated by the second way as above 
(see the book of Hejhal \cite{H} for more details) and the 
work of Singer \cite{S} and Millson \cite{M}, we give the 
trace formulas of weight $k$ for cocompact, torsion-free discrete 
subgroups of $SU(2, 1)$ and prove the analogue of the Riemann 
hypothesis on compact complex surfaces $M$ with $c_1^2(M)=3 
c_2(M)$, where $c_i(M)$ is the $i$-th Chern class of $M$, 
$c_2(M)$ is a multiple of $3$ and $c_2(M)>0$. 

\vskip 0.4 cm  
1.1. {\it Notation and terminology}
\vskip 0.2 cm
  
  Let $SU(2, 1)=\{ g \in SL(3, {\Bbb C}): g^{*} J g=J \}$ 
where $J=\left(\matrix
    &   & -1\\
    & 1 &   \\
 -1 &   &	
   \endmatrix\right)$. The corresponding complex domain 
${\frak S}_{2}=\{(z_1, z_2) \in {\Bbb C}^2: z_1+\overline{z_1}
-z_2 \overline{z_2}>0 \}$ is the simplest example of a bounded 
symmetric domain which is not a tube domain. Let us introduce 
the following differential operator
$$\aligned
  \Delta_{k}
&=(z_1+\overline{z_1}-z_2 \overline{z_2})\\
&\times \left[ (z_1+\overline{z_1}) \frac{\partial^2}{\partial z_1
 \partial \overline{z_1}}+\frac{\partial^2}{\partial z_2 \partial
 \overline{z_2}}+z_2 \frac{\partial^2}{\partial \overline{z_1}
 \partial z_2}+\overline{z_2} \frac{\partial^2}{\partial z_1
 \partial \overline{z_2}}-k(\frac{\partial}{\partial z_1}-
 \frac{\partial}{\partial \overline{z_1}}) \right], 
\endaligned\tag 1.1$$
where $2k$ is an integer (see \cite{Y} for more details). 
In particular, $\Delta=\Delta_0$ is the Laplace-Beltrami 
operator of $SU(2, 1)$ on ${\frak S}_{2}$. 

  Let $\Gamma \subset SU(2, 1)$ be a discrete subgroup. 
Assume that $\Gamma$ is properly discontinuous, 
$\Gamma \backslash {\frak S}_{2}$ is compact, and $\Gamma$ acts 
freely on ${\frak S}_{2}$. Set $M=\Gamma \backslash {\frak S}_{2}$, 
then $M$ is a compact complex algebraic surface. Such $\Gamma$ 
is a uniform, torsion-free discrete subgroup of $SU(2, 1)$. 
Hence, every $\gamma \in \Gamma$ which is not the identity 
is hyperbolic. The Jordan canonical form of hyperbolic 
elements takes the form:
$\gamma=\gamma(\mu, \theta)=\text{diag}(\mu e^{i \theta},                
 e^{-2 i \theta}, \mu^{-1} e^{i \theta})$, 
where $\mu>1$ and $\theta \in [0, 2 \pi)$. 

  An automorphic form of weight $2k$ on ${\frak S}_2$ is 
given by 
$$f(\gamma(Z))=J(\gamma, Z)^{2k} f(Z),\tag 1.2$$ 
where 
$$J(\gamma, Z)=\frac{j(\gamma, Z)}{|j(\gamma, Z)|} \quad \text{with}
  \quad j(\gamma, Z)=c_1 z_1+c_2 z_2+c_3\tag 1.3$$
for $\gamma=\left(\matrix
     * &   * & *\\
     * &   * & *\\
   c_1 & c_2 & c_3  
   \endmatrix\right) \in \Gamma$. 
 
  A holomorphic function $g(Z)$ is called a modular form of
weight $k$ if it satisfies that
$$g(\gamma(Z))=j(\gamma, Z)^{k} g(Z), \quad \text{for} \quad
  \gamma \in \Gamma.\tag 1.4$$
In fact, if $g(Z)$ is a modular form of weight $2k$, then
$f(Z)=\rho(Z)^{k} g(Z)$ with $\rho(Z)=z_1+\overline{z_1}-z_2
\overline{z_2}$ is an automorphic form of weight $2k$.
Moreover, 
$$\Delta_{k} f(Z)=k(k-2) f(Z).\tag 1.5$$
Conversely, if $f(Z)$ is an automorphic form of weight $2k$, 
then $g(Z)=\rho(Z)^{-k} f(Z)$ satisfies that
$$g(\gamma(Z))=j(\gamma, Z)^{2k} g(Z).$$

  Let ${\Cal F}$ be the fundamental region (i.e. the fundamental
domain) of $\Gamma$. Set
$$L^2(\Gamma \backslash {\frak S}_{2}, k)=\{f \in L^2({\Cal F}): 
  f(\gamma(Z))=J(\gamma, Z)^{2k} f(Z) \quad \text{for} \quad 
  \gamma \in \Gamma\}.\tag 1.6$$
We study the spectral theory of $L^2(\Gamma \backslash 
{\frak S}_{2}, k)$, and prove the completeness of eigenfunction 
expansions (see Theorem 2.11, Proposition 2.12, Proposition 2.13 
and Theorem 2.14).  

  Let $\Phi \in C_{c}^{2}({\Bbb R})$ and consider the function 
$\Phi(u(Z, Z^{\prime}))$, where $u(Z, Z^{\prime})=\frac{|\rho(Z, 
Z^{\prime})|^2}{\rho(Z) \rho(Z^{\prime})}-1$ with
$\rho(Z, Z^{\prime})=\overline{z_1}+z_1^{\prime}-\overline{z_2} 
 z_2^{\prime}$, 
$\rho(Z)=\rho(Z, Z)$ for $Z=(z_1, z_2) \in {\frak S}_{2}$,
$Z^{\prime}=(z_1^{\prime}, z_2^{\prime}) \in {\frak S}_{2}$. 
In fact, $u(Z, Z^{\prime})$ is a point-pair invariant, i.e., 
$u(g(Z), g(Z^{\prime}))=u(Z, Z^{\prime})$ for every
$g \in SU(2, 1)$ (see \cite{Y} for more details.)
Set $H_{k}(Z, Z^{\prime})=\frac{\rho(Z, Z^{\prime})^{2k}}
{|\rho(Z, Z^{\prime})|^{2k}}$. Let 
$K_{k}(Z, Z^{\prime})=H_k(Z, Z^{\prime}) \Phi(u(Z, Z^{\prime}))$.
The integral operator of weight $k$ is given by
$$L_{k} f(Z)=\int_{{\frak S}_{2}} K_{k}(Z, Z^{\prime}) 
         f(Z^{\prime}) dm(Z^{\prime}).\tag 1.7$$
We have (see Theorem 2.5) 
$$L_{k} f(Z)=\Lambda_{k}(\lambda) f(Z),\tag 1.8$$ 
where $f$ is an eigenfunction of $\Delta_k$ on ${\frak S}_{2}$: 
$\Delta_k f=\lambda f$ and $\Lambda_{k}(\lambda)$ depends only 
on $\lambda$, $k$ and $\Phi$. A representative set of 
eigenfunctions is given by $\Delta_k \rho^s=\lambda \rho^s$ 
with $\lambda=s(s-2)$. Let $s=1+ir$, we will use the $r$ 
instead of the $\lambda$ as parameter. Set 
$h_{k}(r)=\Lambda_{k}(-1-r^2)$. When $k=0$, the connection 
between $\Phi(u)$ and $h(r):=h_0(r)$ is given by the relations.
$$\aligned
 P(v)=\int_{v}^{\infty} \Phi(u) \arccos \left(\frac{1-u+2v}{1+u}
      \right) du, \quad
 &\Phi(u)=\frac{1}{\pi} \int_{u}^{\infty} \frac{1}{\sqrt{v-u}} 
       d(\sqrt{v+1} P^{\prime}(v)).\\ 
 P(v)=g(\eta) \quad \text{with} \quad 
 &v=\frac{1}{4} (e^{\eta}+e^{-\eta}-2).\\ 
 h(r)=2 \pi \int_{-\infty}^{\infty} g(\eta) e^{ir \eta} d \eta, 
 \quad &g(\eta)=\frac{1}{4 \pi^2} \int_{-\infty}^{\infty} 
 h(r) e^{-ir \eta} dr.
\endaligned\tag 1.9$$	   
When $k=1$, the connection between $\Phi(u)$ and $h_1(r)$ 
is given by   
$$\aligned
  P_1(v)=\int_{v}^{\infty} \frac{2 \Phi(u)}{u+1}
  \sqrt{(v+1)(u-v)} du, \quad &
  \Phi(u)=\frac{u+1}{\pi} \int_{u}^{\infty}
  \frac{1}{\sqrt{v-u}} d \left[ \left(\frac{P_1(v)}
  {\sqrt{v+1}}\right)^{\prime}\right].\\
  P_1(v)=g_{1}(\eta), \quad \text{with} \quad 
  &v=\frac{1}{4}(e^{\eta}+e^{-\eta}-2),\\
  h_1(r)=2 \pi \int_{-\infty}^{\infty} g_1(\eta) e^{ir \eta} d \eta, 
  \quad 
 &g_1(\eta)=\frac{1}{4 \pi^2} \int_{-\infty}^{\infty} h_1(r)
  e^{-ir \eta} dr.  
\endaligned\tag 1.10$$  
When $k=\frac{1}{2}$, the connection between $\Phi(u)$ and 
$h_{\frac{1}{2}}(r)$ is given by
$$\aligned
  Q_{0}(v)=2 \int_{v}^{\infty} \Phi(u) \sqrt{\frac{u-v}{u+1}} du,
  \quad 
 &\Phi(u)=\frac{1}{\pi} \sqrt{u+1} \int_{u}^{\infty}
  \frac{d Q_{0}^{\prime}(v)}{\sqrt{v-u}}.\\
  Q_{0}(v)=g_{\frac{1}{2}}(\eta), \quad \text{with} \quad
  &v=\frac{1}{4}(e^{\eta}+e^{-\eta}-2),\\
  h_{\frac{1}{2}}(r)=2 \pi \int_{-\infty}^{\infty} 
  g_{\frac{1}{2}}(\eta) e^{ir \eta} d \eta, \quad
 &g_{\frac{1}{2}}(\eta)=\frac{1}{4 \pi^2} \int_{-\infty}^{\infty}
  h_{\frac{1}{2}}(r) e^{-ir \eta} dr. 
\endaligned\tag 1.11$$

  The zeta function of weight $k$ associated with the group 
$<\gamma_1> \times <\gamma_2>$ generated by $\gamma_1$ and 
$\gamma_2$, where $<\gamma_1>$ is an infinite cyclic group, 
$<\gamma_2>$ is a finite cyclic group, is given by 
$$\aligned
 &Z_{<\gamma_1> \times <\gamma_2>}(k, s)\\
=&\prod_{j=0}^{\infty} \prod_{l=0}^{\infty} \prod_{n=0}^{\infty} 
  \left[1-N(\gamma_1)^{-(s+j+\frac{l+n}{2})}\right]
  ^{\frac{3}{\text{ord}(\gamma_2)} \sum_{q=1}^{\text{ord}(\gamma_2)} 
  \exp([2k+3(l-n)]iq \theta(\gamma_2))},
\endaligned\tag 1.12$$
for $\text{Re}(s)>2$ and $2k \in {\Bbb Z}$. For $\Gamma$ a 
uniform, torsion-free discrete subgroup of $SU(2, 1)$, we 
define
$$Z_{k}(s)=Z_{\Gamma}(k, s)=\prod_{\{(\gamma_1, \gamma_2)\}} 
  Z_{<\gamma_1> \times <\gamma_2>}(k, s) \tag 1.13$$
where the pair $\{(\gamma_1, \gamma_2) \}$ runs through a 
set of representatives of primitive conjugacy classes of 
$\Gamma$ (see Proposition 2.1 for the details).

\vskip 0.4 cm
1.2. {\it The main results}
\vskip 0.2 cm
  
  Note that $c_2(M)$ is equal to the Euler number 
$e(M)=\sum_{i=0}^{4} (-1)^{i} \dim H^{i}(M, {\Bbb R})$ 
since $M$ is compact. Here $M$ is considered as a 
Riemannian manifold of dimension four. Now, we can 
 state the main theorem in the present paper.

{\smc Theorem 1.1 (Main Theorem 1)}. {\it Assume that  
$h_{k}(r)$ is an analytic function on $|\text{Im}(r)| 
\leq 1+\delta$, $h_{k}(-r)=h_{k}(r)$ and $|h_{k}(r)| 
\leq M [1+|\text{Re}(r)|]^{-4-\delta}$, where $\delta$ 
and $M$ are positive constants. Moreover, suppose that
$$g_{k}(\eta)=\frac{1}{4 \pi^2} \int_{-\infty}^{\infty}
  h_{k}(r) e^{-ir \eta} dr, \quad \eta \in {\Bbb R}.$$
Then the trace formulas of weight $k$ for $SU(2, 1)$ 
take the following form:
$$\aligned
  \text{Tr}(L_{k})=
 &\sum_{n=0}^{\infty} h_{k}(r_n)
 =\frac{2}{3} c_2(M) \int_{-\infty}^{\infty} h_{k}(r) 
  r (r^2+k^2) \frac{\cosh 2 \pi r+\cos 2k \pi}{\sinh 2 \pi r} dr+\\
 &+2 \pi \sum_{\{(\gamma_1, \gamma_2)\}} \frac{3}{\text{ord}(\gamma_2)}  
  \sum_{q=1}^{\text{ord}(\gamma_2)} \sum_{m=1}^{\infty} 
  \frac{\ln N(\gamma_1)}{[N(\gamma_1)^{\frac{m}{2}}-
  N(\gamma_1)^{-\frac{m}{2}}]}\\
 &\times \frac{e^{2kiq \theta(\gamma_2)}}{|N(\gamma_1)^{
  \frac{m}{4}} e^{\frac{3}{2} i q \theta(\gamma_2)}-
  N(\gamma_1)^{-\frac{m}{4}} e^{-\frac{3}{2} i 
  q \theta(\gamma_2)}|^2} g_{k}(m \ln N(\gamma_1)),
\endaligned\tag 1.14$$
where $c_2(M) \geq 3$, $k=0, \frac{1}{2}, 1$ and 
$h_0=h$, $g_0=g$.}

  Using the concept of Poincar\'{e} map, the trace formulas 
can be expressed as (see Theorem 3.2):
$$\aligned
 &\text{Tr}(L_{k})=\frac{2}{3} c_2(M) \int_{-\infty}^{\infty} 
  h_{k}(r) r (r^2+k^2) \frac{\cosh 2 \pi r+\cos 2 k \pi}
  {\sinh 2 \pi r} dr+\\
 &+2 \pi \sum_{\{(\gamma_1, \gamma_2)\}} \frac{3}{\text{ord}(\gamma_2)} 
  \sum_{q=1}^{\text{ord}(\gamma_2)} \sum_{m=1}^{\infty} 
  \frac{\ln N(\gamma_1)}{|\det(I-P_{\gamma_1^{m} 
  \gamma_2^{q}})|^{\frac{1}{2}}} e^{2kiq \theta(\gamma_2)} 
  g_{k}(m \ln N(\gamma_1)).
\endaligned\tag 1.15$$

  It is well-known that the $c$-function of Harish-Chandra 
(see \cite{Ga})
$$c(r)^{-1}=\frac{\Gamma(\frac{1}{2}(p+q))}{\Gamma(p+q)}
            \frac{\Gamma(ir+\frac{p}{2})}{\Gamma(ir)}
			\frac{\Gamma(\frac{ir}{2}+\frac{p}{4}+\frac{q}{2})}
			{\Gamma(\frac{ir}{2}+\frac{p}{4})}.$$
In the case of $SU(2, 1)$, $p=2$, $q=1$, one has
$$c(2r)^{-1} c(-2r)^{-1}=\frac{\pi}{4} \frac{r^3}{\tanh \pi r}.
  \tag 1.16$$
Up to a constant, it coincides with the value of
$r(r^2+k^2) \frac{\cosh 2 \pi r+\cos 2k \pi}{\sinh 2 \pi r}$
at $k=0$. By Theorem 1.1, we can derive the functional equations 
of the zeta functions $Z_{k}(s)$. In fact, we conjecture that
the above trace formulas are valid for all $k$ with $2k \in 
{\Bbb Z}$. In the case of $k=0$, $k=\frac{1}{2}$ and $k=1$, we 
can write down the trace formulas explicitly. For every 
$k$ with $2k \in {\Bbb Z}$, there is an analytic continuation 
and a functional equation for the zeta function $Z_{k}(s)$. 
In particular, we can write down them explicitly for $Z_0(s)$,
$Z_{\frac{1}{2}}(s)$ and $Z_1(s)$ according to the weight $k=0$,
$k=\frac{1}{2}$ or $k=1$, respectively.
  
{\smc Theorem 1.2 (Main Theorem 2)}. {\it The functional 
equations of the zeta functions of weight $k$ are given by
$$\aligned
 &Z_k(s)=Z_k(2-s) \exp \left[\frac{8}{3} \pi c_2(M) 
  \int_{0}^{s-1} v(v^2-k^2) \cot \pi v dv \right], \quad k=0, 1;\\
 &Z_k(s)=Z_k(2-s) \exp \left[\frac{8}{3} \pi c_2(M)
  \int_{0}^{s-1} v(k^2-v^2) \tan \pi v dv \right], \quad 
  k=\frac{1}{2},  
\endaligned\tag 1.17$$
where $M=\Gamma \backslash {\frak S}_{2}$ is a compact complex
algebraic surface with $c_1^2(M)=3 c_2(M)$.}

  The zeros of zeta functions $Z_k(s)$ ($k=0, \frac{1}{2}, 1$) have 
the following properties. In particular, the last property 
tells us that the analogue of the Riemann hypothesis is true.
  
{\smc Theorem 1.3 (Main Theorem 3)}. {\it Let 
$$Z_{k}(s)=\prod_{\{(\gamma_1, \gamma_2)\}} \prod_{j=0}^{\infty} 
  \prod_{l=0}^{\infty} \prod_{n=0}^{\infty} 
  \left[1-N(\gamma_1)^{-s-j-\frac{l+n}{2}}\right]
  ^{\frac{3}{\text{ord}(\gamma_2)} \sum_{q=1}^{\text{ord}(\gamma_2)} 
  \exp([2k+3(l-n)]iq \theta(\gamma_2))}$$
where $\text{Re}(s)>2$ and $k=0, \frac{1}{2}, 1$. Then	   
\roster
\item $Z_{k}(s)$ are actually entire functions.
\item $Z_{k}(s)$ $(k=0, 1)$ have trivial zeros $s=-m$, $m \geq 0$, 
      with multiplicity $\frac{8}{3} c_2(M) (m+1)^{3}$;
	  $Z_k(s)$ $(k=\frac{1}{2})$ has trivial zeros $s=-m+\frac{1}{2}$,
	  $m \geq 0$, with multiplicity $\frac{4}{3} c_2(M) m(m+1)(2m+1)$. 
\item The nontrivial zeros of $Z_{k}(s)$ are located at $1 \pm i r_n$.
\endroster}

  It is interesting that we can deduce that $c_2(M)$ is a multiple
of $3$ by (2) of Theorem 1.3. Thus $c_2(M) \geq 3$, which is the
first main theorem in \cite{HP}.

  This paper consists of three chapters. In chapter two, we study 
the spectral theory of $L^2(\Gamma \backslash {\frak S}_{2}, k)$, 
and prove the completeness of eigenfunction expansions. We give
the integral transformations and their inverse which are related 
to the group $SU(2, 1)$ and the symmetric space ${\frak S}_{2}$.
These transformations play an important role in the trace formulas.
In chapter three, we obtain the trace formulas of weight $k$ on 
$SU(2, 1)$ for cocompact, torsion-free discrete subgroups. 
In the rest of this chapter, we give the definition of the 
zeta functions of weight $k$ and their expansions. Using 
the trace formulas, we get the functional equations of the 
zeta functions of weight $k$. At last, we prove the analogue 
of the Riemann hypothesis for a class of compact complex surfaces.  

\vskip 0.5 cm
\centerline{\bf 2. The spectral theory of $L^2(\Gamma \backslash
                   {\frak S}_{2}, k)$}
\vskip 0.5 cm

2.1. {\it The construction of hyperbolic elements}
\vskip 0.2 cm
  
  At first, let us recall some basic results about the discrete 
subgroups of Lie groups (see \cite{R}). Let $G$ be a connected
Lie group and $\Gamma$ a discrete subgroup such that $G/\Gamma$
is compact. Then $\Gamma$ is finitely presentable (see \cite{R},
p.95, Theorem 6.15). In fact, let $M$ be a compact connected 
manifold. Then the fundamental group $\pi_1(M, x)$ of $M$ at
a point $x \in M$ is finitely presentable (see \cite{R}, p.95,
Theorem 6.16). Hence, the uniform torsion-free subgroup $\Gamma
\subset SU(2, 1)$ is finitely presentable. If $\Gamma \subset G$ 
is a uniform, torsion-free discrete subgroup then every 
$\gamma \in \Gamma$ which is not the identity is hyperbolic.
Now, let us recall the classifications of linear transformations
in $SU(2, 1)$ (see \cite{G}). There are four classes: the identity 
element, the hyperbolic elements, the elliptic elements and the 
parabolic elements. In the ball model, the hyperbolic elements 
take the form:   
$$T=\left(\matrix
   e^{i \theta} \cosh u & 0 & e^{i \theta} \sinh u\\
   0                    & e^{-2 i \theta} & 0     \\
   e^{i \theta} \sinh u & 0 & e^{i \theta} \cosh u
  \endmatrix\right).$$
Now, we give the expression of hyperbolic elements in the
unbounded realization ${\frak S}_{2}$. By the Cayley transformation 
$C: {\Bbb B}^2=\{(w_1, w_2) \in {\Bbb C}^2: |w_1|^2+|w_2|^2<1 \} 
 \to {\frak S}_{2}=\{(z_1, z_2) \in {\Bbb C}^2: z_1+\overline{z_1}
 -z_2 \overline{z_2}>0 \}$, where
$C=\left(\matrix 
    -\overline{\omega} &   & -\omega\\
                       & 1 &        \\
    -1                 &   & 1 
   \endmatrix\right)$, 
$C(w_1, w_2)=\left(\frac{-\overline{\omega} w_1-\omega}{-w_1+1},
 \frac{w_2}{-w_1+1}\right)$ and 
$C^{-1}=\left(\matrix
         1 &   & \omega\\
           & 1 &       \\ 
         1 &   & -\overline{\omega}
        \endmatrix\right)$,
$C^{-1}(z_1, z_2)=\left(\frac{z_1+\omega}{z_1-\overline{\omega}},
 \frac{z_2}{z_1-\overline{\omega}}\right)$, we have
$$g=C T C^{-1}=\left(\matrix
  e^{u} e^{i \theta} &  & (\omega-\overline{\omega}) \sinh u 
                          \cdot e^{i \theta}\\
                     & e^{-2 i \theta} &  \\
  0                  &                 & e^{-u} e^{i \theta} 
           \endmatrix\right) \in SU(2, 1),$$  
where $SU(2, 1)=\{g \in SL(3, {\Bbb C}): g^{*} J g=J \}$ with
$J=\left(\matrix
        &   & -1\\     
        & 1 &   \\
     -1 &   &
   \endmatrix\right)$.
The characteristic polynomial of $g$ is 
$\det(\lambda I-g)=(\lambda-e^{u} e^{i \theta})
(\lambda-e^{-2 i \theta})(\lambda-e^{-u} e^{i \theta})$.
Hence, the Jordan canonical form of $g$ takes the form
$\gamma=\text{diag}(e^{u} e^{i \theta}, e^{-2i \theta},
e^{-u} e^{i \theta}) \in SU(2, 1)$. Let 
$h=\left(\matrix
 \alpha &       & -\frac{\sqrt{3}}{2} i \overline{\alpha}^{-1}\\
        & \beta &                                             \\ 
        &       & \overline{\alpha}^{-1}
   \endmatrix\right)$
where $\alpha, \beta \in {\Bbb C}$, $|\beta|=1$ and
$\alpha \beta \overline{\alpha}^{-1}=1$. Then $h \in SU(2, 1)$
and $h^{-1} g h=\gamma$, i.e., $g$ is conjugate in $SU(2, 1)$ 
to $\gamma$. Without loss of generality, we can assume that 
$u>0$. Set $\mu=e^{u}$, then $\mu>1$ and  
$\gamma=\text{diag}(\mu e^{i \theta}, e^{-2i \theta},                           
        \mu^{-1} e^{i \theta})$.  

{\smc Proposition 2.1}. {\it Suppose that $\gamma=\text{diag}(\mu
e^{i \theta}, e^{-2 i \theta}, \mu^{-1} e^{i \theta}) \in \Gamma$
is hyperbolic. Then the centralizer $Z(\gamma)$ is a direct product
of two cyclic subgroups of $\Gamma$, $Z(\gamma)=<\gamma_1> \times
<\gamma_2>$, where $\gamma_1=\text{diag}(r_0, 1, r_0^{-1})$ $(r_0>1)$
is hyperbolic and
$\gamma_2=\text{diag}(e^{i \varphi_0}, e^{-2 i \varphi_0},
e^{i \varphi_0})$ with order $\text{ord}(\gamma_2)$ and
$\varphi_0=\frac{2 \pi}{\text{ord}(\gamma_2)}$. Moreover, 
$\gamma=\gamma_1^m \gamma_2^n$ with $m \geq 1$ and $1 \leq n
\leq \text{ord}(\gamma_2)$. The two generators $\gamma_1$,
$\gamma_2$ are uniquely determined by $\gamma$.}		
		
{\it Proof}. If $g \in Z(\gamma)$, then $g$ has the form		
$g=\text{diag}(r e^{i \varphi}, e^{-2 i \varphi}, r^{-1} 
e^{i \varphi})$ with $r>0$. Note that $\Gamma$ is 
discontinuous, so $Z(\gamma)=<\gamma_1> \times <\gamma_2>$,		
i.e., $g=\gamma_1^m \gamma_2^n$, $m \geq 1$ and $1 \leq n
\leq \text{ord}(\gamma_2)$.
The uniqueness of $\gamma_1$ and $\gamma_2$ is proved by the
same argument since $m \geq 1$ and $1 \leq n \leq 
\text{ord}(\gamma_2)$. 
\flushpar		
$\qquad \qquad \qquad \qquad \qquad \qquad \qquad \qquad \qquad
 \qquad \qquad \qquad \qquad \qquad \qquad \qquad \qquad \qquad
 \quad \square$		

  $N(\gamma_1):=r_0^2$ in Proposition 2.1 is the minimal norm 
of a hyperbolic element from $Z(\gamma)$. Following Selberg 
\cite{Se}, we call $\gamma_1$ a primitive element for $\gamma$ 
in $\Gamma$. The pair $(\gamma_1, \gamma_2)$ in Proposition 2.1 
is called a primitive pair for $\gamma$ in $\Gamma$.  
		
\vskip 0.4 cm
2.2. {\it The basic integral operator of weight $0$}
\vskip 0.2 cm
  
  Put
$$L f(Z)=\int_{{\frak S}_2} k(Z, Z^{\prime}) f(Z^{\prime})
         dm(Z^{\prime}).\tag 2.1$$
Let $f(Z)$ be an eigenfunction of the Laplace-Beltrami operator
$\Delta$ on ${\frak S}_{2}$: $\Delta f(Z)=\lambda f(Z)$.
According to \cite{Se} (or see Theorem 2.5 in the special case
$k=0$), for an integral operator we may write
$$\int_{{\frak S}_{2}} k(Z, Z^{\prime}) f(Z^{\prime}) dm(Z^{\prime})
 =\Lambda(\lambda) f(Z),\tag 2.2$$
where $\Lambda(\lambda)$ depends only on $\lambda$ and $k$.
A representative set of eigenfunctions is given by $\rho(Z)^{s}$
since $\Delta \rho(Z)^s=\lambda \rho(Z)^s$ with $\lambda=s(s-2)$.
Set $s=1+ir$ with $r \in {\Bbb C}$, then $s(s-2)=-1-r^2$.
Let $h(r)=\Lambda(-1-r^2)$. The connection between 
$\Phi(u)$ and $h(r)$ is given by the following theorem. 
  
{\smc Theorem 2.2}. {\it Let $\Phi(u) \in C_{c}^{2}({\Bbb R})$ 
be a positive function of real variable. Set
$$P(v)=\int_{v}^{\infty} \Phi(u) \arccos \left(\frac{1-u+2v}{1+u}
       \right) du.\tag 2.3$$
Moreover, define $g(\eta)$ by $P(v)=g(\eta)$ with $v=\frac{1}{4}
(e^{\eta}+e^{-\eta}-2)$. Then one has the following integral transform
$$h(r)=2 \pi \int_{-\infty}^{\infty} g(\eta) e^{ir \eta} d \eta, \quad
  g(\eta)=\frac{1}{4 \pi^2} \int_{-\infty}^{\infty} h(r) e^{-ir \eta}
  dr,\tag 2.4$$ 
and
$$\Phi(u)=\frac{1}{\pi} \int_{u}^{\infty} \frac{1}{\sqrt{v-u}}
          d(\sqrt{v+1} P^{\prime}(v)).\tag 2.5$$}
  
{\it Proof}. Let us start with the formula:
$$\int_{{\frak S}_{2}} k(Z_0, Z) \rho(Z)^{s} dm(Z)=\Lambda(\lambda) 
  \rho(Z_0)^{s}, \quad \lambda=s(s-2).$$
Let $Z_0=(\frac{1}{2}, 0)$ and 
$Z=(z_1, z_2)=(\frac{1}{2}(\rho+|z|^2)+it, z)$, where
$\rho>0$, $t \in {\Bbb R}$ and $z \in {\Bbb C}$. Then $\rho(Z_0)=1$, 
$|\rho(Z_0, Z)|^2=\frac{1}{4}(1+\rho+|z|^2)^2+t^2$. Thus
$$\int_{0}^{\infty} \int_{-\infty}^{\infty} \int_{{\Bbb C}}
  \Phi \left(\frac{\frac{1}{4}(1+\rho+|z|^2)^2+t^2}{\rho}-1\right)
  dz d \overline{z} dt \rho^{s-3} d \rho=\Lambda(\lambda).$$
Let $z=\sqrt{r} e^{i \theta}$, $r>0$, $\theta \in [0, 2 \pi)$, then
$dz \wedge d \overline{z}=-i dr \wedge d \theta$ and
$$4 \pi \int_{0}^{\infty} \int_{0}^{\infty} \int_{0}^{\infty}
  \Phi \left(\frac{\frac{1}{4}(1+\rho+r)^2+t^2}{\rho}-1\right)
  dr dt \rho^{s-3} d \rho=\Lambda(\lambda).$$
Set $\xi=\frac{\frac{1}{4}(1+\rho+r)^2}{\rho}-1$, 
$\eta=\frac{t^2}{\rho}$. Then 
$\xi \in [\frac{(\rho-1)^2}{4 \rho}, \infty)$, $\eta \in [0, \infty)$
and $dr=\frac{\sqrt{\rho}}{\sqrt{\xi+1}} d \xi$,
$dt=\frac{\sqrt{\rho}}{2 \sqrt{\eta}} d \eta$.
We have
$$2 \pi \int_{0}^{\infty} \int_{0}^{\infty} 
  \int_{\frac{(\rho-1)^2}{4 \rho}}^{\infty} \Phi(\xi+\eta)
  \frac{d \xi}{\sqrt{\xi+1}} \frac{d \eta}{\sqrt{\eta}}
  \rho^{s-2} d \rho=\Lambda(\lambda).$$
Put $w=\xi+\eta$, then
$$2 \pi \int_{0}^{\infty} \int_{\frac{(\rho-1)^2}{4 \rho}}^{\infty}
  \int_{\xi}^{\infty} \Phi(w) \frac{d w}{\sqrt{w- \xi}}
  \frac{d \xi}{\sqrt{\xi+1}} \rho^{s-2} d \rho=\Lambda(\lambda).$$
Here
$$\aligned
 &\int_{\frac{(\rho-1)^2}{4 \rho}}^{\infty} \int_{\xi}^{\infty}
  \Phi(w) \frac{d w}{\sqrt{w- \xi}} \frac{d \xi}{\sqrt{\xi+1}}\\
=&\int_{\frac{(\rho-1)^2}{4 \rho}}^{\infty} \Phi(w) dw
  \int_{\frac{(\rho-1)^2}{4 \rho}}^{w} \frac{d \xi}
  {\sqrt{(w- \xi)(\xi+1)}}\\
=&\int_{\frac{(\rho-1)^2}{4 \rho}}^{\infty} \Phi(w) 
  \arccos \left(\frac{\frac{(\rho-1)^2}{2 \rho}-w+1}{w+1}\right) dw   
=P \left(\frac{(\rho-1)^2}{4 \rho}\right).
\endaligned$$
We conclude that
$$\Lambda[s(s-2)]=2 \pi \int_{0}^{\infty} \rho^{s-2} P \left(
  \frac{(\rho-1)^2}{4 \rho} \right) d \rho.$$
Set $s=1+ir$ and $g(\eta)=P(\frac{1}{4}(e^{\eta}+e^{-\eta}-2))$, then
$\Lambda(\lambda)=\Lambda(-1-r^2)=h(r)$. On the other hand,
let $\rho=e^{\eta}$, then
$$\aligned
 &2 \pi \int_{0}^{\infty} \rho^{s-2} P \left(\frac{(\rho-1)^2}
  {4 \rho}\right) d \rho
 =2 \pi \int_{-\infty}^{\infty} e^{(s-2) \eta} P \left(\frac{1}{4}
  (e^{\eta}+e^{-\eta}-2)\right) e^{\eta} d \eta\\ 
=&2 \pi \int_{-\infty}^{\infty} e^{(s-1) \eta} g(\eta) d \eta 
 =2 \pi \int_{-\infty}^{\infty} e^{ir \eta} g(\eta) d \eta. 
\endaligned$$
Consequently,
$h(r)=2 \pi \int_{-\infty}^{\infty} e^{ir \eta} g(\eta) d \eta$.
So, $g(\eta)=\frac{1}{4 \pi^2} \int_{-\infty}^{\infty} h(r) 
e^{-ir \eta} dr$.
By
$$P(v)=\int_{v}^{\infty} \Phi(u) \arccos \left(\frac{1-u+2v}{1+u}\right)
       du,$$
we have
$$P^{\prime}(v)=-\frac{1}{\sqrt{v+1}} \int_{v}^{\infty} 
  \frac{\Phi(u)}{\sqrt{u-v}} du.$$
Set
$Q(v)=-\sqrt{v+1} P^{\prime}(v)$,
then
$Q(v)=\int_{v}^{\infty} \frac{\Phi(u)}{\sqrt{u-v}} du$.
By \cite{Ku}, p.56, Theorem 5.3.1 or \cite{H}, p.15, 
Proposition 4.1, we have
$\Phi(u)=-\frac{1}{\pi} \int_{u}^{\infty} \frac{dQ(v)}{\sqrt{v-u}}$.
Thus,
$$\Phi(u)=\frac{1}{\pi} \int_{u}^{\infty} \frac{1}{\sqrt{v-u}}
          d(\sqrt{v+1} P^{\prime}(v)).$$
$\qquad \qquad \qquad \qquad \qquad \qquad \qquad \qquad \qquad
 \qquad \qquad \qquad \qquad \qquad \qquad \qquad \qquad \qquad
 \quad \square$ 

{\smc Theorem 2.3}. {\it Let $k(Z, Z^{\prime})$ be a point pair
invariant as in Theorem 2.2, and put 
$$n(a, b)=\left(\matrix
          1 & a & \frac{|a|^2}{2}+ib\\
		    & 1 & \overline{a}      \\
			&   & 1
          \endmatrix\right), \quad a \in {\Bbb C}, b \in {\Bbb R}.$$
Then
$$\int_{\Bbb R} \int_{\Bbb C} k(Z, n(a, b)Z^{\prime}) 
  da d \overline{a} db
 =2 \pi \rho \rho^{\prime} g(\log \rho^{\prime}-\log \rho),\tag 2.6$$
where $\rho=\rho(Z)$ and $\rho^{\prime}=\rho(Z^{\prime})$ and
$g$ is as in Theorem 2.2.}

{\it Proof}. Let
$$n=n(a, b)=\left(\matrix
            1 & a & \frac{|a|^2}{2}+ib\\
              & 1 & \overline{a}\\
              &   & 1
            \endmatrix\right), \quad a \in {\Bbb C}, b \in {\Bbb R},$$
then
$$n(a, b)(Z^{\prime})=(z_1^{\prime}+a z_2^{\prime}+\frac{|a|^2}{2}
  +ib, z_2^{\prime}+\overline{a}).$$
So
$$\rho(Z, n(a, b)Z^{\prime})=\rho(Z, Z^{\prime})+a z_2^{\prime}
 -\overline{a} \overline{z_2}+\frac{|a|^2}{2}+ib.$$
Set $Z=(\frac{\rho+|z|^2}{2}+it, z)$ and 
$Z^{\prime}=(\frac{\rho^{\prime}+|z^{\prime}|^2}{2}+it^{\prime},
 z^{\prime})$, then 
$$\rho(Z, Z^{\prime})=\frac{1}{2}(\rho+\rho^{\prime}+|z-z^{\prime}|^2)
 +i(t^{\prime}-t+\text{Im} z \overline{z^{\prime}}).$$
Put $a=u+iv$, $u, v \in {\Bbb R}$, $z=x+iy$ and 
$z^{\prime}=x^{\prime}+i y^{\prime}$, then
$da \wedge d \overline{a}=-2i du \wedge dv$ and 
$$\aligned
  \rho(Z, n(a, b)Z^{\prime})
=&\frac{1}{2}[\rho+\rho^{\prime}+(u+x^{\prime}-x)^2+
  (v+y-y^{\prime})^2]\\
+&i[v(x+x^{\prime})+u(y+y^{\prime})+t^{\prime}-t+\text{Im} z
  \overline{z^{\prime}}+b].
\endaligned$$
Hence,
$$\aligned
 &\int_{\Bbb R} \int_{\Bbb C} k(Z, n(a, b)Z^{\prime}) 
  da d \overline{a} db\\
=&2 \int_{\Bbb R} \int_{\Bbb C} \Phi \left(\frac{|\frac{1}{2}
  (\rho+\rho^{\prime}+u^2+v^2)+ib|^2}{\rho \rho^{\prime}}-1\right) 
  du dv db\\
=&2 \int_{\Bbb R} \int_{0}^{2 \pi} \int_{0}^{\infty}
  \Phi \left(\frac{|\frac{1}{2}(\rho+\rho^{\prime}+r^2)+ib|^2}
  {\rho \rho^{\prime}}-1 \right) r dr d \theta db\\
=&2 \pi \int_{\Bbb R} \int_{0}^{\infty} 
  \Phi \left(\frac{|\frac{1}{2}(\rho+\rho^{\prime}+r)+ib|^2}
  {\rho \rho^{\prime}}-1 \right) dr db\\
=&4 \pi \int_{0}^{\infty} \int_{\rho+\rho^{\prime}}^{\infty}
  \Phi \left(\frac{\frac{1}{4} r^2+b^2}{\rho \rho^{\prime}}-1\right)
  dr db.
\endaligned$$

  Let $t=\frac{b^2}{\rho \rho^{\prime}}$ and 
$u=\frac{r^2}{4 \rho \rho^{\prime}}$, then
$$db=\frac{\sqrt{\rho \rho^{\prime}}}{2 \sqrt{t}} dt, \quad
  dr=\frac{\sqrt{\rho \rho^{\prime}}}{\sqrt{u}} du.$$
We have
$$\aligned
 &\int_{\Bbb R} \int_{\Bbb C} k(Z, n(a, b)Z^{\prime}) 
  da d \overline{a} db\\
=&2 \pi \rho \rho^{\prime} \int_{0}^{\infty} 
  \int_{\frac{(\rho+\rho^{\prime})^2}{4 \rho \rho^{\prime}}}^{\infty}
  \Phi(u+t-1) \frac{du}{\sqrt{u}} \frac{dt}{\sqrt{t}}\\
=&2 \pi \rho \rho^{\prime} \int_{0}^{\infty} 
  \int_{\frac{(\rho-\rho^{\prime})^2}{4 \rho \rho^{\prime}}}^{\infty}
  \Phi(v+t) \frac{dv}{\sqrt{v+1}} \frac{dt}{\sqrt{t}}.  
\endaligned$$

  Set $\rho^{\prime}/\rho=e^{\eta}$, then
$\frac{(\rho-\rho^{\prime})^2}{4 \rho \rho^{\prime}}
=\sinh^2 \frac{\eta}{2}$. Thus,
$$\aligned
 &\int_{\Bbb R} \int_{\Bbb C} k(Z, n(a, b)Z^{\prime}) 
  da d \overline{a} db\\
=&2 \pi \rho \rho^{\prime} \int_{\sinh^2 \frac{\eta}{2}}^{\infty}
  \int_{v}^{\infty} \Phi(u) \frac{du}{\sqrt{u-v}} 
  \frac{dv}{\sqrt{v+1}}\\  
=&2 \pi \rho \rho^{\prime} \int_{\sinh^2 \frac{\eta}{2}}^{\infty}
  \Phi(u) du \int_{\sinh^2 \frac{\eta}{2}}^{u} 
  \frac{dv}{\sqrt{(u-v)(v+1)}}\\
=&2 \pi \rho \rho^{\prime} \int_{\sinh^2 \frac{\eta}{2}}^{\infty} \Phi(u) 
  \arccos \left(\frac{1-u+2 \sinh^2 \frac{\eta}{2}}{1+u}\right) du\\
=&2 \pi \rho \rho^{\prime} P(\sinh^2 \frac{\eta}{2})
 =2 \pi \rho \rho^{\prime} g(\eta)
 =2\pi \rho \rho^{\prime} g(\log \rho^{\prime}-\log \rho).
\endaligned$$
$\qquad \qquad \qquad \qquad \qquad \qquad \qquad \qquad \qquad 
 \qquad \qquad \qquad \qquad \qquad \qquad \qquad \qquad \qquad
 \quad \square$

\vskip 0.4 cm
2.3. {\it The basic integral operator of weight $k$}
\vskip 0.2 cm
  
  Now, we give the integral transformation of weight $k$. 
Let ${\Cal F}$ be the fundamental region (i.e. the fundamental 
domain) of $\Gamma$. Since $\Gamma$ is a cocompact discrete 
subgroup of $SU(2, 1)$, ${\Cal F}$ is compact. Set 
$$H_{k}(Z, W)=\frac{\rho(Z, W)^{2k}}{|\rho(Z, W)|^{2k}},$$ 
where 
$\rho(Z, W)=\overline{z_1}+w_1-\overline{z_2} w_2$.
Let 
$$K_{k}(Z, W)=H_k(Z, W) \Phi \left(\frac{|\rho(Z, W)|^2}
  {\rho(Z) \rho(W)}-1\right).$$

{\it Definition} {\smc 2.4}. An automorphic form of weight $2k$ on 
${\frak S}_2$ is given by
$$f(\gamma(Z))=J(\gamma, Z)^{2k} f(Z),\tag 2.7$$
where $J(\gamma, Z)=\frac{j(\gamma, Z)}{|j(\gamma, Z)|}$ with
$j(\gamma, Z)=c_1 z_1+c_2 z_2+c_3$
for $\gamma=\left(\matrix
     * &   * & *\\
     * &   * & *\\
   c_1 & c_2 & c_3  
   \endmatrix\right) \in \Gamma$.
Let 
$$L^2(\Gamma \backslash {\frak S}_{2}, k)=\{f \in L^2({\Cal F}): 
  f(\gamma(Z))=J(\gamma, Z)^{2k} f(Z) \quad \text{for} \quad 
  \gamma \in \Gamma\}.$$ 

  Following the method in \cite{H} (see \cite{H}, pp.364-369,
Proposition 2.14), we give the next theorem.

{\smc Theorem 2.5}. {\it Let $f \in C^2({\frak S}_{2})$ be any
eigenfunction of $\Delta_k$: $\Delta_k f=\lambda f$. Then the
integral operator of weight $k$
$$L_{k} f(Z):=\int_{{\frak S}_{2}} K_{k}(Z, W) f(W) dm(W)=
  \Lambda_{k}(\lambda) f(Z),\tag 2.8$$
where $\Lambda_{k}(\lambda)$ depends only on $(\Phi, k, \lambda)$.
The value of $\Lambda_{k}(\lambda)$ is independent of $f$.}   
   
{\it Proof}. We define the operator $P_{\gamma}$ by   
$$(P_{\gamma} u)(Z)=u(\gamma(Z)) J(\gamma, Z)^{-2k}, \quad
  \gamma \in SU(2, 1).\tag 2.9$$
It is known that (see \cite{Y}, Proposition 4.2)   
$\Delta_{k} P_{\gamma}=P_{\gamma} \Delta_{k}$.   
Let $K$ be the maximal compact subgroup of $SU(2, 1)$. Every
$\gamma \in K$ can be given by
$$\gamma=C \left(\matrix
           a & b & 0\\
           c & d & 0\\
           0 & 0 & e^{i \theta}
      \endmatrix\right) C^{-1},$$   
where $\left(\matrix
       a & b\\
       c & d
	   \endmatrix\right) \in U(2)$, $\theta \in [0, 2 \pi)$
and $(ad-bc) e^{i \theta}=1$. 
Write explicitly,	   
$$\gamma=\left(\matrix
    -\overline{\omega} a-\omega e^{i \theta}
   &-\overline{\omega} b & -a+e^{i \theta}\\	  
    c & d & \omega c\\
    -a+e^{i \theta} & -b & 
	-\omega a-\overline{\omega} e^{i \theta}
    \endmatrix\right).$$
We have $\gamma(-\omega, 0)=(-\omega, 0)$.

  We now claim that without loss of generality $Z=(-\omega, 0)$.
Suppose, in fact, that $Z=(-\omega, 0)$ is OK. For general $Z$,
we write $Z=\gamma(-\omega, 0)$ and check that (using \cite{Y},
Proposition 4.7)
$$\aligned
 &\int_{{\frak S}_{2}} K_{k}(Z, W) f(W) dm(W)\\
=&\int_{{\frak S}_{2}} K_{k}(\gamma(-\omega, 0), W) f(W) dm(W)\\
=&\int_{{\frak S}_{2}} K_{k}(\gamma(-\omega, 0), \gamma(\xi))
  f(\gamma(\xi)) dm(\xi)\\
=&\int_{{\frak S}_{2}} J(\gamma, (-\omega, 0))^{2k}
  K_{k}((-\omega, 0), \xi) J(\gamma, \xi)^{-2k} f(\gamma(\xi))
  dm(\xi)\\
=&J(\gamma, (-\omega, 0))^{2k} \int_{{\frak S}_{2}} 
  K_{k}((-\omega, 0), \xi) P_{\gamma} f(\xi) dm(\xi)\\
=&J(\gamma, (-\omega, 0))^{2k} \Lambda_{k}(\lambda)
  (P_{\gamma}f)(-\omega, 0)\\
=&\Lambda_{k}(\lambda) f[\gamma(-\omega, 0)]
 =\Lambda_{k}(\lambda) f(Z),
\endaligned$$   
since $\Delta_{k} P_{\gamma} f=P_{\gamma} \Delta_{k} f=
\lambda (P_{\gamma} f)$.   
   
  We may restrict ourselves to the case $Z=(-\omega, 0)$.
We claim that without loss of generality
$$P_{\gamma} f=P_{\gamma(\theta)} f=e^{-2ki \theta} f, 
  \quad \gamma \in K.$$       
Suppose, in fact, that the result is known to be true for
such symmetric $f$. We then take the original eigenfunction
$f$ and form
$$A=\int_{{\frak S}_{2}} K_{k}((-\omega, 0), W) f(W) dm(W).$$   
It is clear that
$$\aligned
A&=\int_{{\frak S}_{2}} K_{k}(\gamma(-\omega, 0), W) f(W) dm(W)\\
 &=\int_{{\frak S}_{2}} K_{k}(\gamma(-\omega, 0), \gamma(\xi))
   f(\gamma(\xi)) dm(\xi)\\
 &=\int_{{\frak S}_{2}} J(\gamma, (-\omega, 0))^{2k}
   K_{k}((-\omega, 0), \xi) J(\gamma, \xi)^{-2k}
   f(\gamma(\xi)) dm(\xi).
\endaligned$$
Here, 
$$J(\gamma, (-\omega, 0))=\frac{(-a+e^{i \theta})(-\omega)+
  (-\omega a-\overline{\omega} e^{i \theta})}{|(-a+
  e^{i \theta})(-\omega)+(-\omega a-\overline{\omega}
  e^{i \theta})|}=e^{i \theta}.$$
Thus,
$$\aligned
A&=\int_{{\frak S}_{2}} e^{2ki \theta} K_{k}((-\omega, 0), \xi)
   P_{\gamma(\theta)} f(\xi) dm(\xi)\\
 &=\frac{1}{2 \pi} \int_{0}^{2 \pi} \int_{{\frak S}_{2}}
   K_{k}((-\omega, 0), \xi) P_{\gamma(\theta)} f(\xi)
   e^{2ki \theta} dm(\xi) d \theta\\
 &=\int_{{\frak S}_{2}} K_{k}((-\omega, 0), \xi) F(\xi) dm(\xi),
\endaligned$$
where
$$F(\xi)=\frac{1}{2 \pi} \int_{0}^{2 \pi} P_{\gamma(\theta)} f(\xi)
         e^{2ki \theta} d \theta.$$
Note that
$$\aligned
  (P_{\gamma_1} P_{\gamma_2} u)(Z)
=&P_{\gamma_1}(u(\gamma_2(Z)) J(\gamma_2, Z)^{-2k})\\
=&u(\gamma_2(\gamma_1(Z))) J(\gamma_2, \gamma_1(Z))^{-2k}
  J(\gamma_1, Z)^{-2k}\\ 
=&u(\gamma_2 \gamma_1(Z)) J(\gamma_2 \gamma_1, Z)^{-2k}\\
=&(P_{\gamma_2 \gamma_1} u)(Z).
\endaligned$$
Thus,
$$P_{\gamma(\varphi)} F(\xi)=\frac{1}{2 \pi} \int_{0}^{2 \pi}
  P_{\gamma(\theta) \gamma(\varphi)} f(\xi) e^{2ki \theta}
  d \theta.$$
Here,   
$$\gamma(\theta) \gamma(\varphi)
 =C \left(\matrix
     * &  \\
       & e^{i \theta}
    \endmatrix\right) C^{-1} \cdot
  C \left(\matrix
     * &  \\
       & e^{i \varphi} 
    \endmatrix\right) C^{-1}
 =C \left(\matrix
	 * &  \\
       & e^{i(\theta+\varphi)}  	
	\endmatrix\right) C^{-1}=\gamma(\theta+\varphi).$$
Therefore,   
$$\aligned
  P_{\gamma(\varphi)} F(\xi)
&=\frac{1}{2 \pi} \int_{0}^{2 \pi} P_{\gamma(\theta+\varphi)}
  f(\xi) e^{2ki \theta} d \theta\\
&=\frac{1}{2 \pi} \int_{\varphi}^{2 \pi+\varphi}
  P_{\gamma(\tau)} f(\xi) e^{2ki(\tau-\varphi)} d \tau\\
&=\frac{1}{2 \pi} \int_{0}^{2 \pi} P_{\gamma(\tau)} f(\xi)
  e^{2ki (\tau-\varphi)} d \tau\\
&=e^{-2ki \varphi} F(\xi).
\endaligned$$   
It is easy to see that $F \in C^2({\frak S}_{2})$ and
$F(-\omega, 0)=f(-\omega, 0)$. By $\Delta_k P_{\gamma}=
P_{\gamma} \Delta_{k}$, we see that $\Delta_k F=\lambda F$.   

  We can apply the symmetric case to deduce that
$$A=\int_{{\frak S}_{2}} K_{k}((-\omega, 0), \xi) F(\xi) dm(\xi)
   =\Lambda_{k}(\lambda) F(-\omega, 0)=\Lambda_{k}(\lambda)
    f(-\omega, 0).$$

  We now restrict not only to $Z=(-\omega, 0)$, but also to
the symmetric case $P_{\gamma(\theta)} f=e^{-2ki \theta} f$.
We introduce a change of variables:
$$W=C^{-1}(Z)=\left(\frac{z_1+\omega}{z_1-\overline{\omega}},
  \frac{z_2}{z_1-\overline{\omega}}\right);
  Z=C(W)=\left(\frac{-\overline{\omega} w_1-\omega}
  {-w_1+1}, \frac{w_2}{-w_1+1}\right),$$
where $Z \in {\frak S}_{2}$ and $W \in {\Bbb B}^{2}$.
By $u(Z)=v(W)$, one defines the operator $\widetilde{\Delta_k}$:
$\Delta_k u=\widetilde{\Delta_k} v$.

  The following equations describe the basic symmetry of
weight $k$:
$$\aligned
  u[\gamma(\theta)(Z)] J(\gamma(\theta), Z)^{-2k}
 &=e^{-2ki \theta} u(Z),\\
  v[C^{-1} \gamma(\theta)(Z)] J(\gamma(\theta), Z)^{-2k}
 &=e^{-2ki \theta} v(W),\\
  v[C^{-1} \gamma(\theta) C(W)] J(\gamma(\theta), Z)^{-2k}
 &=e^{-2ki \theta} v(W). 
\endaligned$$
Here
$C^{-1} \gamma(\theta) C=\left(\matrix
 a & b &   \\ 
 c & d &   \\
   &   & e^{i \theta}
 \endmatrix\right)$.
We have
$$v \left[\left(\matrix
  a & b &     \\
  c & d &     \\
    &   & e^{i \theta}
  \endmatrix\right)(w_1, w_2)\right]
  \frac{[(-a+e^{i \theta}) z_1-b z_2-\omega a- 
  \overline{\omega} e^{i \theta}]^{-2k}}{|\text{same}|^{-2k}}
 =e^{-2ki \theta} v(w_1, w_2),$$
Using the expressions $z_1=\frac{-\overline{\omega} w_1-\omega}
{-w_1+1}$ and $z_2=\frac{w_2}{-w_1+1}$, we have
$$v(e^{-i \theta}(a w_1+b w_2), e^{-i \theta}(c w_1+d w_2))
  \frac{|1-e^{-i \theta}(a w_1+b w_2)|^{2k}}{[1-e^{-i \theta}
  (a w_1+b w_2)]^{2k}}=v(w_1, w_2) \frac{|1-w_1|^{2k}}{(1-w_1)^{2k}}.$$
This implies that $v(w_1, w_2) (1-w_1)^{-k} (1-\overline{w_1})^{k}$
must be radially symmetric. We shall therefore write
$$u(Z)=v(W)=\psi(W) (1-w_1)^{k} (1-\overline{w_1})^{-k}
 =\psi(W) (1-w_1)^{\alpha} (1-\overline{w_1})^{\beta}.$$
Corresponding to this, we define
$\Delta_k u=\widetilde{\Delta_k} v=
 \widetilde{\widetilde{\Delta_k}} \psi$.
The associated change of variables is very well-behaved
analytically.
 
  In our previous paper \cite{Y}, the Laplace-Beltrami operator
on ${\Bbb B}^{2}$ is given by
$$\aligned
  \Delta_{{\Bbb B}^{2}}
=&(1-|w_1|^2-|w_2|^2)\\
\times &\left[(1-|w_1|^2) \frac{\partial^2}
  {\partial w_1 \partial \overline{w_1}}+(1-|w_2|^2)
  \frac{\partial^2}{\partial w_2 \partial \overline{w_2}}
  -w_1 \overline{w_2} \frac{\partial^2}{\partial w_1
  \partial \overline{w_2}}-\overline{w_1} w_2 \frac{\partial^2}
  {\partial \overline{w_1} \partial w_2}\right].
\endaligned$$  
We conclude that near the origin
$$\aligned
  \widetilde{\widetilde{\Delta_k}} \psi
=&A_{1, 1} \psi_{w_1, \overline{w_1}}+A_{1, 2} \psi_{w_1, 
  \overline{w_2}}+A_{2, 1} \psi_{w_2, \overline{w_1}}+
  A_{2, 2} \psi_{w_2, \overline{w_2}}\\
 &+B_1 \psi_{w_1}+C_1 \psi_{\overline{w_1}}
  +B_2 \psi_{w_2}+C_2 \psi_{\overline{w_2}}+D \psi
\endaligned$$
with real analytic coefficients and $A_{1, 1}(0, 0; 0, 0) \neq 0$,  
$A_{2, 2}(0, 0; 0, 0) \neq 0$, where
$A_{i, j}=A_{i, j}(w_1, \overline{w_1}; w_2, \overline{w_2})$,
$B_j=B_j(w_1, \overline{w_1}; w_2, \overline{w_2})$,
$C_j=C_j(w_1, \overline{w_1}; w_2, \overline{w_2})$,
$1 \leq i, j \leq 2$ and $D=D(w_1, \overline{w_1};
w_2, \overline{w_2})$.

  Now, we must study the differential equation
$$\widetilde{\widetilde{\Delta_k}} \psi=\lambda \psi
  (1-w_1)^{\alpha} (1-\overline{w_1})^{\beta},$$  
along the various rays $\arg(w_1)=\theta$. 
Now, we introduce polar coordinates
$w_1=r e^{i \theta_1} \cos \varphi$ and 
$w_2=r e^{i \theta_2} \sin \varphi$ with 
$\varphi \in (0, \frac{\pi}{2})$.
Then
$$\aligned
  \frac{\partial}{\partial w_1} 
=&\frac{1}{2} e^{-i \theta_1} \cos \varphi 
  \frac{\partial}{\partial r}-\frac{i}{2} 
  \frac{e^{-i \theta_1}}{r \cos \varphi} 
  \frac{\partial}{\partial \theta_1}-
  \frac{e^{-i \theta_1}}{2r} \sin \varphi
  \frac{\partial}{\partial \varphi},\\
  \frac{\partial}{\partial w_2}
=&\frac{1}{2} e^{-i \theta_2} \sin \varphi
  \frac{\partial}{\partial r}-\frac{i}{2}
  \frac{e^{-i \theta_2}}{r \sin \varphi}
  \frac{\partial}{\partial \theta_2}+
  \frac{e^{-i \theta_2}}{2r} \cos \varphi 
  \frac{\partial}{\partial \varphi}.
\endaligned\tag 2.10$$
Hence,
$$\aligned
  \Delta
=&(1-r^2) [\frac{1}{4}(1-r^2) \frac{\partial^2}{\partial r^2}
  +\frac{3-r^2}{4 r} \frac{\partial}{\partial r}+\frac{1}{4 r^2}
  \frac{\partial^2}{\partial \varphi^2}+\frac{1}{2 r^2} 
  \cot 2 \varphi \frac{\partial}{\partial \varphi}\\
 &+\frac{1}{4 r^2}(\frac{1}{\cos^2 \varphi} \frac{\partial^2}
  {\partial \theta_1^2}+\frac{1}{\sin^2 \varphi} \frac{\partial^2}
  {\partial \theta_2^2})-\frac{1}{4}(\frac{\partial}{\partial 
  \theta_1}+\frac{\partial}{\partial \theta_2})^{2}].
\endaligned\tag 2.11$$
If $\psi$ is radially symmetric, we will get a $\theta_1$-family 
of homogeneous linear ordinary differential equation with
regular singular point $r=0$. The indicial equation is
$\nu(\nu+2)=0$ for each $\theta_1$. The general theory of
differential equations tells us that $\psi$ is unique up to 
a proportionality factor whenever $\psi$ is known to be
smooth at the origin. In our case, we conclude that
$f(z_1, z_2)=f(-\omega, 0) g(z_1, z_2)$, where the
eigenfunction $g(z_1, z_2)$ is symmetric and normalized
by $g(-\omega, 0)=1$. Th function $g(z_1, z_2)$ depends 
only on $(k, \lambda)$. Hence, we can define
$$\Lambda_{k}(\lambda)=\int_{{\frak S}_{2}} K_{k}((-\omega,
  0), W) g(W) dm(W).$$
This completes the proof of Theorem 2.5.
\flushpar
$\qquad \qquad \qquad \qquad \qquad \qquad \qquad \qquad \qquad
 \qquad \qquad \qquad \qquad \qquad \qquad \qquad \qquad \qquad
 \quad \square$
    
{\smc Theorem 2.6}. {\it Let $\Phi(u) \in C_{c}^{2}({\Bbb R})$ 
be a positive function of real variable. Set
$$P_1(v)=\int_{v}^{\infty} \frac{2 \Phi(u)}{u+1}
  \sqrt{(v+1)(u-v)} du.\tag 2.12$$  
Moreover, define $g_1(\eta)$ by $P_1(v)=g_{1}(\eta)$, with 
$v=\frac{1}{4}(e^{\eta}+e^{-\eta}-2)$. 
Then one has the following integral transform
$$h_1(r)=2 \pi \int_{-\infty}^{\infty} g_1(\eta) e^{ir \eta} d \eta, 
  \quad 
  g_1(\eta)=\frac{1}{4 \pi^2} \int_{-\infty}^{\infty} h_1(r)
  e^{-ir \eta} dr,$$  
and
$$\Phi(u)=\frac{u+1}{\pi} \int_{u}^{\infty}
  \frac{1}{\sqrt{v-u}} d \left[ \left(\frac{P_1(v)}{\sqrt{v+1}}
  \right)^{\prime}\right].$$
Set
$$Q_0(v)=2 \int_{v}^{\infty} \Phi(u) \sqrt{\frac{u-v}{u+1}} du.
  \tag 2.13$$  
Moreover, define $g_{\frac{1}{2}}(\eta)$ by 
$Q_0(v)=g_{\frac{1}{2}}(\eta)$, with 
$v=\frac{1}{4}(e^{\eta}+e^{-\eta}-2)$. 
Then one has the following integral transform
$$h_{\frac{1}{2}}(r)=2 \pi \int_{-\infty}^{\infty} 
  g_{\frac{1}{2}}(\eta) e^{ir \eta} d \eta, \quad 
  g_{\frac{1}{2}}(\eta)=\frac{1}{4 \pi^2} 
  \int_{-\infty}^{\infty} h_{\frac{1}{2}}(r) e^{-ir \eta} dr,$$  
and
$$\Phi(u)=\frac{1}{\pi} \sqrt{u+1} \int_{u}^{\infty}
  \frac{d Q_0^{\prime}(v)}{\sqrt{v-u}}.$$}

{\it Proof}. By Theorem 2.5,
$$\int_{{\frak S}_{2}} K_{k}(Z, W) \rho(W)^{s} dm(W)
 =\Lambda_{k}(\lambda) \rho(Z)^{s}, \quad \lambda=s(s-2).\tag 2.14$$
Set $Z_0=(\frac{1}{2}, 0)$ and $Z=(z_1, z_2)=(\frac{1}{2}(\rho+
|z|^2)+it, z)$ with $t \in {\Bbb R}$, $\rho>0$ and $z \in {\Bbb C}$.
Then
$$\aligned
  \Lambda_{k}(\lambda)
&=\int_{{\frak S}_{2}} H_{k}((\frac{1}{2}, 0), (z_1, z_2))
  \Phi \left[\frac{|\frac{1}{2}+z_1|^2}{\rho(Z)}-1\right]
  \rho(Z)^s dm(Z)\\
&=\int_{{\frak S}_{2}} \frac{[\frac{1}{2}(1+\rho+|z|^2)+it]^{2k}}
  {|\frac{1}{2}(1+\rho+|z|^2)+it|^{2k}} \Phi \left[\frac{\frac{1}{4}
  (1+\rho+|z|^2)^2+t^2}{\rho}-1\right] \rho^s \frac{dt d \rho dz
  d \overline{z}}{\rho^3}\\
&=\pi \int_{0}^{\infty} \int_{-\infty}^{\infty}
  \int_{0}^{\infty} \frac{[(1+\rho+r)+it]^{2k}}{|(1+\rho+r)+it|^{2k}}
  \Phi \left[\frac{(1+\rho+r)^2+t^2}{4 \rho}-1\right] \rho^{s-3}
  dr dt d \rho\\
&=\pi \int_{0}^{\infty} \int_{-\infty}^{\infty}
  \int_{0}^{\infty} \frac{\sum_{m=0}^{2k} \binom{2k}{m} 
  (1+\rho+r)^{m} i^{2k-m} t^{2k-m}}{[(1+\rho+r)^2+t^2]^{k}}\\
&\times \Phi \left[\frac{(1+\rho+r)^2+t^2}{4 \rho}-1\right]
 \rho^{s-3} dr dt d \rho.    
\endaligned$$
The contribution from the odd powers of $t$ is zero. 
When $2k$ is even, 
$$\aligned
  \Lambda_{k}(\lambda)
=&\pi (-1)^{k} \int_{0}^{\infty} \int_{-\infty}^{\infty}
  \int_{0}^{\infty} \frac{\sum_{m=0}^{k} \binom{2k}{2m} 
  (1+\rho+r)^{2m} i^{-2m} t^{2k-2m}}{[(1+\rho+r)^2+t^2]^{k}}\\
 &\times \Phi \left[\frac{(1+\rho+r)^2+t^2}{4 \rho}-1\right] 
  \rho^{s-3} dr dt d \rho.  
\endaligned$$
Set $u^2=\frac{t^2}{4 \rho}$ and $v=\frac{(r+\rho+1)^2}{4 \rho}-1$,
then $u \in (-\infty, \infty)$, 
$v \in [\frac{(\rho-1)^2}{4 \rho}, \infty)$. We have
$dt=2 \sqrt{\rho} du$ and $dr=\sqrt{\frac{\rho}{v+1}} dv$. Hence,
$$\aligned
  \Lambda_{k}(\lambda)
&=2 \pi (-1)^{k} \int_{0}^{\infty} \int_{-\infty}^{\infty} 
  \int_{\frac{(\rho-1)^2}{4 \rho}}^{\infty} \sum_{m=0}^{k} 
  \binom{2k}{2m} (-1)^{m} (v+1)^{m- \frac{1}{2}} u^{2(k-m)}\\
&\times (u^2+v+1)^{-k} \Phi(u^2+v) \rho^{s-2} dv du d \rho\\  
&=4 \pi (-1)^{k} \int_{0}^{\infty} \int_{0}^{\infty} 
  \int_{\frac{(\rho-1)^2}{4 \rho}}^{\infty}  
  \sum_{m=0}^{k} \binom{2k}{2m} (-1)^{m} (v+1)^{m- \frac{1}{2}}
  u^{2(k-m)}\\ 
&\times (u^2+v+1)^{-k} \Phi(u^2+v) \rho^{s-2} dv du d \rho.
\endaligned$$
Let $w=u^2$, then
$$\aligned
  \Lambda_{k}(\lambda)
&=2 \pi (-1)^{k} \int_{0}^{\infty} \int_{0}^{\infty} 
  \int_{\frac{(\rho-1)^2}{4 \rho}}^{\infty}
  \sum_{m=0}^{k} \binom{2k}{2m} (-1)^{m} (v+1)^{m- \frac{1}{2}}
  w^{k-m- \frac{1}{2}}\\
&\times (w+v+1)^{-k} \Phi(w+v) \rho^{s-2} dv dw d \rho.  
\endaligned$$
Set $u=w+v$, then 
$$\aligned
  \Lambda_{k}(\lambda)
&=2 \pi (-1)^{k} \int_{0}^{\infty} \rho^{s-2} d \rho 
  \int_{\frac{(\rho-1)^2}{4 \rho}}^{\infty} \int_{v}^{\infty}
  \sum_{m=0}^{k} \binom{2k}{2m} (-1)^{m}\\
&\times (v+1)^{m- \frac{1}{2}} (u-v)^{k-m- \frac{1}{2}}
  (u+1)^{-k} \Phi(u) du dv\\
&=2 \pi (-1)^{k} \int_{0}^{\infty} \rho^{s-2} d \rho \sum_{m=0}^{k}
  \binom{2k}{2m} (-1)^{m} \int_{\frac{(\rho-1)^2}{4 \rho}}^{\infty}
  \Phi(u)\\
&\times (u+1)^{-k} du \int_{\frac{(\rho-1)^2}{4 \rho}}^{u} 
  (v+1)^{m- \frac{1}{2}} (u-v)^{k-m- \frac{1}{2}} dv.
\endaligned$$
Set
$$P_k(w)=(-1)^{k} \int_{w}^{\infty} \sum_{m=0}^{k} \binom{2k}{2m} 
  (-1)^{m} \Phi(u) (u+1)^{-k} du \int_{w}^{u} (v+1)^{m- \frac{1}{2}} 
  (u-v)^{k-m- \frac{1}{2}} dv.$$

  When $2k$ is odd, say that $2k=2p+1$. Similarly, we have
$$\aligned
  \Lambda_{k}(\lambda)
&=2 \pi (-1)^{p} \int_{0}^{\infty} \rho^{s-2} d \rho \sum_{l=0}^{p}
  \binom{2p+1}{2l+1} (-1)^{l} \int_{\frac{(\rho-1)^2}{4 \rho}}^{\infty}
  \Phi(u)\\
&\times (u+1)^{-p-\frac{1}{2}} du \int_{\frac{(\rho-1)^2}{4 \rho}}^{u} 
  (v+1)^{l} (u-v)^{p-l-\frac{1}{2}} dv.
\endaligned$$
Set
$$\aligned
 &P_k(w)=Q_{p}(w)\\
=&(-1)^{p} \int_{w}^{\infty} \sum_{l=0}^{p}
  \binom{2p+1}{2l+1} (-1)^{l} \Phi(u) (u+1)^{-p-\frac{1}{2}} du 
  \int_{w}^{u} (v+1)^{l} (u-v)^{p-l-\frac{1}{2}} dv.
\endaligned$$

  For $2k \in {\Bbb Z}$, we conclude that
$$\Lambda_{k}[s(s-2)]=2 \pi \int_{0}^{\infty} \rho^{s-2} 
  P_k \left(\frac{(\rho-1)^2}{4 \rho}\right) d \rho.$$
Set $s=1+ir$ and $g_{k}(\eta)=P_k \left(\frac{1}{4}(e^{\eta}+
e^{-\eta}-2)\right)$ with $\rho=e^{\eta}$. Consequently,
$$h_{k}(r)=\Lambda_{k}(-1-r^2)=2 \pi \int_{-\infty}^{\infty}
  e^{ir \eta} g_k(\eta) d \eta, \quad r \in {\Bbb C}.$$
Then
$g_k(\eta)=\frac{1}{4 \pi^2} \int_{-\infty}^{\infty} h_k(r)
  e^{-ir \eta} dr$.  

  In fact, when $k$ is a positive integer,
$$P_{k}(w)=\int_{w}^{\infty} \Phi(u) (u+1)^{-k} du
           \int_{w}^{u} \sum_{m=0}^{k} (-1)^{m} \binom{2k}{2m}
		   (v+1)^{k-m-\frac{1}{2}} (u-v)^{m-\frac{1}{2}} dv.$$
We have
$$P_{k}^{\prime}(w)=\sum_{m=0}^{k} (-1)^{m+1} \binom{2k}{2m}
  (w+1)^{k-m-\frac{1}{2}} \int_{w}^{\infty} \Phi(u) (u+1)^{-k}
  (u-w)^{m-\frac{1}{2}} du.$$
Thus,
$$\aligned
  P_{k}^{\prime}(v)
=&\int_{v}^{\infty} \frac{\sum_{m=0}^{k} \binom{2k}{2m} (u-v)^{m} 
  (-1)^{m} (v+1)^{k-m}}{(u+1)^{k}} \frac{-\Phi(u)}
  {\sqrt{(u-v)(v+1)}} du\\
=&\int_{v}^{\infty} \frac{\sum_{m=0}^{k} \binom{2k}{2m} (u-v)^{m}
  (-v-1)^{k-m}}{(-u-1)^{k}} \frac{-\Phi(u)}{\sqrt{(u-v)(v+1)}} du.
\endaligned$$
A straightforward calculation gives
$$\aligned
  P_{k}^{\prime}(v)
=&\int_{v}^{\infty} \frac{\sum_{m=0}^{k} \binom{2k}{2m} (u-v)^{m}
  [(u-v)+(-u-1)]^{k-m}}{(-u-1)^{k}} \frac{-\Phi(u)}
  {\sqrt{(u-v)(v+1)}} du\\
=&\int_{v}^{\infty} \frac{\sum_{m=0}^{k} \sum_{j=0}^{k-m} 
  \binom{2k}{2m} \binom{k-m}{j} (u-v)^{m} (u-v)^{k-m-j} 
  (-u-1)^{j}}{(-u-1)^{k}}\\
 &\times \frac{-\Phi(u)}{\sqrt{(u-v)(v+1)}} du\\
=&\sum_{m=0}^{k} \sum_{j=0}^{k-m} \binom{2k}{2m} \binom{k-m}{j}
  \int_{v}^{\infty} \frac{(u-v)^{k-j}}{(-u-1)^{k-j}}
  \frac{-\Phi(u)}{\sqrt{(u-v)(v+1)}} du\\
=&\sum_{j=0}^{k} c_{j}(k) \int_{v}^{\infty} 
  \frac{(u-v)^{k-j}}{(-u-1)^{k-j}} 
  \frac{-\Phi(u)}{\sqrt{(u-v)(v+1)}} du,
\endaligned$$
where $c_{j}(k)=\sum_{m=0}^{k-j} \binom{2k}{2m} \binom{k-m}{j}$.

  Now, we need the following lemma (see \cite{H}, p.396, Lemma 4.5):

{\smc Lemma}. {\it Suppose that $\theta(t) \in C({\Bbb R}^{+})$ and 
that $\theta(t)=O((t+1)^{-\beta})$ for some $\beta>k+\frac{1}{2}$.
Then
$$\int_{x}^{\infty} (t-x)^{k-\frac{1}{2}} \theta(t) dt
 =k! \binom{-\frac{1}{2}}k \int_{x}^{\infty} 
  \frac{\theta^{(-k)}(t)}{\sqrt{t-x}} dt.$$
It is understood here that $k$ is a non-negative integer, $x \geq 0$
and the antiderivative $\theta^{(-k)}(t)$ is computed with base point
$+\infty$.}

  Assume that $\Phi(t)=O((t+1)^{-1})$, by the above lemma, we have
$$P_{k}^{\prime}(v)=-\frac{1}{\sqrt{v+1}} \sum_{j=0}^{k} c_{j}(k)
  (k-j)! \binom{-\frac{1}{2}}{k-j} \int_{v}^{\infty} 
  \frac{\left[\frac{\Phi(u)}{(-u-1)^{k-j}}
  \right]^{(j-k)}}{\sqrt{u-v}} du.$$    
Therefore,
$$P_{k}^{\prime}(v)=-\frac{1}{\sqrt{v+1}} \int_{v}^{\infty}
  \frac{W_{k}(u)}{\sqrt{u-v}} du,$$
where
$$W_{k}(u)=\sum_{j=0}^{k} c_{j}(k) (k-j)! \binom{-\frac{1}{2}}{k-j} 
  \left[\frac{\Phi(u)}{(-u-1)^{k-j}}\right]^{(j-k)}.$$

  In particular, when $k=1$,
$$\aligned
  P_1(w)
=&\int_{w}^{\infty} \Phi(u) (u+1)^{-1} du \int_{w}^{u}
  \left(\sqrt{\frac{v+1}{u-v}}-\sqrt{\frac{u-v}{v+1}}\right) dv\\
=&\int_{w}^{\infty} \Phi(u) (u+1)^{-1} du \left[\frac{\pi}{2}(u+1)
  -2 \sqrt{(u-v)(v+1)}\right]|_{v=w}^{u}\\  
=&\int_{w}^{\infty} \Phi(u) \frac{2}{u+1} \sqrt{(u-w)(w+1)} du.
\endaligned$$
Thus,
$-\left(\frac{P_1(v)}{\sqrt{v+1}}\right)^{\prime}
 =\int_{v}^{\infty} \frac{\Phi(u)}{u+1} \frac{1}{\sqrt{u-v}} du$.  
We have  
$$\frac{\Phi(u)}{u+1}=\frac{1}{\pi} \int_{u}^{\infty}
  \frac{1}{\sqrt{v-u}} d \left[ \left(\frac{P_1(v)}{\sqrt{v+1}}
  \right)^{\prime}\right].$$  

  When $k=\frac{1}{2}$, i.e. $p=0$, we have
$$Q_0(w)=2 \int_{w}^{\infty} \Phi(u) \sqrt{\frac{u-w}{u+1}} du.$$
Hence, 
$$Q_0^{\prime}(w)=-\int_{w}^{\infty} \frac{\Phi(u)}{\sqrt{u+1}}
  \frac{1}{\sqrt{u-w}} du.$$
Thus,
$$\frac{\Phi(u)}{\sqrt{u+1}}=\frac{1}{\pi} \int_{u}^{\infty}
  \frac{d Q_0^{\prime}(w)}{\sqrt{w-u}}.$$
$\qquad \qquad \qquad \qquad \qquad \qquad \qquad \qquad \qquad
  \qquad \qquad \qquad \qquad \qquad \qquad \qquad \qquad \qquad
  \quad \square$

  Now, let us recall some basic facts about the growth of volume
of geodesic balls in a complex hyperbolic space (see \cite{Go},
pp.104-105, 3.3.4). Since the Euclidean distance $r$ from $0$ 
is related to the hyperbolic distance $\rho$ by $r=\tanh \left(
\frac{\rho}{2} \right)$, the volume of a ball of radius $\rho$ 
is given by
$$\text{Vol}(B(\rho))=\frac{8^n \sigma_{2n-1}}{2n} \sinh^{2n}
  \left(\frac{\rho}{2}\right) \sim \frac{8^n \sigma_{2n-1}}{2n}
  e^{n \rho},$$
where $\sigma_{2n-1}=\frac{2 \pi^n}{n!}$ is the Euclidean
volume of the unit sphere $S^{2n-1} \subset {\Bbb C}^{n}$.
Thus the volume of a geodesic ball has exponential growth
rate $n$. In our case, $n=2$ and 
$\text{Vol}(B(\rho)) \sim O(e^{2 \rho})$.
By the same method as in \cite{H}, p.5, Proposition 2.2,
we can prove that
$$\sharp \{ \gamma \in \Gamma: \gamma({\Cal F}) \cap 
  B(Z_0; \rho) \neq \emptyset \}=O(e^{2 \rho}),$$
where $Z_0 \in {\Cal F}$ and $B(Z_0; \rho)$ denotes the 
geodesic ball about $Z_0$ with radius $\rho$. 

{\it Definition} {\smc 2.7}. The automorphic kernel 
function is given by
$$G_{k}(Z, W)=\sum_{\gamma \in \Gamma} K_{k}(Z, \gamma(W))
              J(\gamma, W)^{2k},\tag 2.15$$
for $(Z, W) \in {\frak S}_{2} \times {\frak S}_{2}$.

{\smc Proposition 2.8}. {\it
\roster
\item $G_{k}(Z, W)$ reduces to a sum of uniformly bounded length.
\item $G_{k}(Z, W)=\overline{G_{k}(W, Z)}$ on ${\frak S}_{2} \times
      {\frak S}_{2}$.
\item $G_{k}(\gamma_1(Z), \gamma_2(W))=J(\gamma_1, Z)^{2k}
      G_{k}(Z, W) J(\gamma_2, W)^{-2k}$ for $\gamma_1, \gamma_2
	  \in \Gamma$.
\endroster}

{\it Proof}. 
(1) The function $\Phi$ has compact support. Hence, 
$K_{k}(Z, \gamma(W))=0$ when $\delta(Z, \gamma(W)) \geq A$ (say),
where $\delta(Z, W)=\frac{1}{2}[\sigma(Z, W)+1]$ (see \cite{Y},
Proposition 3.4). The question reduces to finding an upper bound
on the number of polygons $\gamma({\Cal F})$, $\gamma \in \Gamma$,
which intersect $\{ \xi: \delta(\xi, Z) \leq A \}$. This is trivial
since ${\Cal F}$ is compact (see also the above argument about 
the growth of volume of geodesic balls).
\flushpar
(2)
$$\aligned
  \overline{G_{k}(W, Z)}
=&\sum_{\gamma \in \Gamma} \overline{K_{k}(W, \gamma(Z))}
  \overline{J(\gamma, Z)^{2k}}\\
=&\sum_{\gamma \in \Gamma} K_{k}(\gamma(Z), W) J(\gamma, Z)^{-2k}\\
=&\sum_{\gamma \in \Gamma} J(\gamma, Z)^{2k} K_{k}(Z, \gamma^{-1}(W))
  J(\gamma, \gamma^{-1}(W))^{-2k} J(\gamma, Z)^{-2k}\\
=&\sum_{\gamma \in \Gamma} K_{k}(Z, \gamma^{-1}(W))
  J(\gamma, \gamma^{-1}(W))^{-2k}.
\endaligned$$
By $1=J(I, W)=J(\gamma, \gamma^{-1}(W)) J(\gamma^{-1}, W)$, we have
$$\overline{G_{k}(W, Z)}=\sum_{\gamma \in \Gamma}
  K_{k}(Z, \gamma^{-1}(W)) J(\gamma^{-1}, W)^{2k}
  =G_{k}(Z, W).$$
(3) Note that
$$\aligned
  G_{k}(Z, g(W))
=&\sum_{\gamma \in \Gamma} K_{k}(Z, \gamma g(W)) J(\gamma, g(W))^{2k}\\
=&\sum_{\gamma \in \Gamma} K_{k}(Z, \gamma g(W)) 
  \frac{J(\gamma g, W)^{2k}}{J(g, W)^{2k}}\\ 
=&\sum_{h \in \Gamma} K_{k}(Z, h(W)) J(h, W)^{2k} J(g, W)^{-2k}\\
=&J(g, W)^{-2k} G_{k}(Z, W).
\endaligned$$
We have
$$\aligned
  G_{k}(\gamma_1(Z), \gamma_2(W))
=&\overline{G_{k}(\gamma_2(W), \gamma_1(Z))}\\
=&\overline{G_{k}(\gamma_2(W), Z) J(\gamma_1, Z)^{-2k}}\\
=&J(\gamma_1, Z)^{2k} \overline{G_{k}(\gamma_2(W), Z)}\\
=&J(\gamma_1, Z)^{2k} G_{k}(Z, \gamma_2(W))\\
=&J(\gamma_1, Z)^{2k} G_{k}(Z, W) J(\gamma_2, W)^{-2k}.
\endaligned$$
$\qquad \qquad \qquad \qquad \qquad \qquad \qquad \qquad \qquad
 \qquad \qquad \qquad \qquad \qquad \qquad \qquad \qquad \qquad
 \quad \square$

{\smc Proposition 2.9}. {\it $L_k$ is a bounded, self-adjoint
operator taking $L^2(\Gamma \backslash {\frak S}_{2}, k)$ into 
$L^2(\Gamma \backslash {\frak S}_{2}, k)$.}

{\it Proof}. For $f \in L^2(\Gamma \backslash {\frak S}_{2}, k)$,
we have 
$$\aligned
  L_k f(Z)
=&\int_{{\frak S}_{2}} K_{k}(Z, W) f(W) dm(W)\\
=&\sum_{\gamma \in \Gamma} \int_{\gamma({\Cal F})}
  K_{k}(Z, W) f(W) dm(W)\\
=&\sum_{\gamma \in \Gamma} \int_{{\Cal F}} K_{k}(Z, \gamma(W))
  f(\gamma(W)) dm(W)\\
=&\sum_{\gamma \in \Gamma} \int_{{\Cal F}} K_{k}(Z, \gamma(W))
  J(\gamma, W)^{2k} f(W) dm(W)\\
=&\int_{{\Cal F}} G_{k}(Z, W) f(W) dm(W).  
\endaligned$$
We have
$$|L_{k}f(Z)|^2 \leq \int_{{\Cal F}} |G_{k}(Z, W)|^{2} dm(W)
  \cdot \int_{{\Cal F}} |f(W)|^{2} dm(W).$$
Hence,
$$\int_{{\Cal F}} |L_{k}f(Z)|^{2} dm(Z) \leq \int_{{\Cal F}}
  \int_{{\Cal F}} |G_{k}(Z, W)|^{2} dm(W) dm(Z) \cdot
  \int_{{\Cal F}} |f(W)|^{2} dm(W).$$

  For $\gamma \in \Gamma$,
$$\aligned
  L_{k}f(\gamma(Z))
=&\int_{{\frak S}_{2}} K_{k}(\gamma(Z), W) f(W) dm(W)\\
=&\int_{{\frak S}_{2}} K_{k}(\gamma(Z), \gamma(\xi)) f(\gamma(\xi))
  dm(\xi)\\
=&\int_{{\frak S}_{2}} J(\gamma, Z)^{2k} K_{k}(Z, \xi) 
  J(\gamma, \xi)^{-2k} \cdot J(\gamma, \xi)^{2k} f(\xi) dm(\xi)\\
=&J(\gamma, Z)^{2k} \int_{{\frak S}_{2}} K_{k}(Z, \xi) f(\xi) dm(\xi).
\endaligned$$
So, $L_{k}f(\gamma(Z))=J(\gamma, Z)^{2k} L_{k}f(Z)$. Therefore, $L_k$
is a bounded linear operator from $L^2(\Gamma \backslash {\frak S}_{2},
k)$ to $L^{2}(\Gamma \backslash {\frak S}_{2}, k)$. Moreover,
$$\aligned
 &(L_k f, g)=\int_{{\Cal F}} L_k f(Z) \overline{g(Z)} dm(Z)\\
=&\int_{{\Cal F}} \int_{{\Cal F}} G_{k}(Z, W) f(W) 
  \overline{g(Z)} dm(W) dm(Z)\\
=&\int_{{\Cal F}} \int_{{\Cal F}} G_{k}(Z, W) f(W)
  \overline{g(Z)} dm(Z) dm(W)\\
=&\int_{{\Cal F}} f(W) \int_{{\Cal F}} G_{k}(Z, W) 
  \overline{g(Z)} dm(Z) dm(W)\\
=&\int_{{\Cal F}} f(W) \int_{{\Cal F}} \overline{G_{k}(W, Z)}
  \overline{g(Z)} dm(Z) dm(W)\\
=&\int_{{\Cal F}} f(W) \overline{L_{k} g(W)} dm(W)=(f, L_k g).
\endaligned$$
Thus, $L_k$ is a self-adjoint operator.
\flushpar
$\qquad \qquad \qquad \qquad \qquad \qquad \qquad \qquad \qquad
 \qquad \qquad \qquad \qquad \qquad \qquad \qquad \qquad \qquad
 \quad \square$

{\smc Proposition 2.10}. {\it The integral operator
$L_{k}: L^2(\Gamma \backslash {\frak S}_{2}, k) \to 
 L^2(\Gamma \backslash {\frak S}_{2}, k)$ is of 
 Hilbert-Schmidt type and is compact.}

{\it Proof}. Recall that $L_k f(Z)=\int_{{\Cal F}} G_{k}(Z, W)
f(W) dm(W)$. Since $\Phi$ has compact support,
$$\int_{{\Cal F}} \int_{{\Cal F}} |G_{k}(Z, W)|^2 dm(Z) dm(W)
  < \infty.$$
By \cite{Yo}, p.277-278, Example 2, we have that $L_k$ is of 
the Hilbert-Schmidt type and is compact.
\flushpar
$\qquad \qquad \qquad \qquad \qquad \qquad \qquad \qquad \qquad
 \qquad \qquad \qquad \qquad \qquad \qquad \qquad \qquad \qquad
 \quad \square$ 
  
\vskip 0.4 cm
2.4. {\it The spectral theory of $L^2(\Gamma \backslash 
     {\frak S}_{2}, k)$}
\vskip 0.2 cm
    
  Motivated by the method as in \cite{ElGM}, pp.145-147, Lemma 2.2, we
get the following theorem:
  
{\smc Theorem 2.11}. {\it Let 
$\varphi_{s, k}(Z, Z^{\prime})=H_{k}(Z, Z^{\prime}) \phi_{s, k}(\sigma)$ 
with
$$\phi_{s, k}(\sigma)=\sigma^{-|k|} (1-\sigma)^{|k|-s} 
  F(s-|k|, s-1-|k|; 2s-1; -\frac{1}{\sigma-1})$$
for $s \in {\Bbb C}$ (see \cite{Y}, pp.35-36), where
$$\sigma=\sigma(Z, Z^{\prime})=\frac{|\rho(Z, Z^{\prime})|^{2}}
        {\rho(Z) \rho(Z^{\prime})}, \quad  
  H_{k}(Z, Z^{\prime})=\frac{\rho(Z, Z^{\prime})^{2k}}
  {|\rho(Z, Z^{\prime})|^{2k}}.$$
In particular, when $k=0$,
$$\varphi_{s, 0}(Z, Z^{\prime})=\phi_{s}(u)=u^{-s} 
  F(s, s-1; 2s-1; -u^{-1}),$$
where $u=u(Z, Z^{\prime})=\sigma(Z, Z^{\prime})-1$ (see \cite{Y}, 
p.18). Assume that $\lambda=s(s-2)$ and 
$g \in C_{c}^{2}({\frak S}_{2})$. The up to a constant,
$$4 \pi^2 g(Z^{\prime})=\int_{{\frak S}_{2}} \varphi_{s, k}
  (Z^{\prime}, Z) (-\Delta_k+\lambda) g(Z) dm(Z).\tag 2.16$$
Moreover,
$$4 \pi^2 g(Z^{\prime})=\int_{{\Cal F}} G_{s, k}(Z^{\prime}, Z) 
  (-\Delta_k+\lambda) g(Z) dm(Z)\tag 2.17$$  
for all $g \in C^2(\Gamma \backslash {\frak S}_{2}, k):=\{ f \in
C^{2}({\Cal F}): f(\gamma(Z))=J(\gamma, Z)^{2k} f(Z)$ for 
$\gamma \in \Gamma \}$, where
$$G_{s, k}(Z^{\prime}, Z)=\sum_{\gamma \in \Gamma} \varphi_{s, k}
  (Z^{\prime}, \gamma(Z)) J(\gamma, Z)^{2k}.$$}  
  
{\it Proof}. Assume that $Z^{\prime}=(-\omega, 0)$. We transform 
${\frak S}_{2}$ isometrically onto ${\Bbb B}^{2}$ such that 
$Z^{\prime}=(-\omega, 0) \mapsto 0$, $Z \mapsto W=(w_1, w_2)$. 
On ${\Bbb B}^{2}$, let us consider the following differential 
form $\omega$: $\omega=(1-|w_1|^2-|w_2|^2)^{-2} \Omega$, with
$$\aligned
  \Omega
=&[(1-|w_1|^2)(f_{w_1} g dw_1 \wedge dw_2 \wedge d \overline{w_2}+
  f g_{\overline{w_1}} d \overline{w_1} \wedge dw_2 \wedge d 
  \overline{w_2})\\
 &-w_1 \overline{w_2} (f_{w_1} g dw_1 \wedge d \overline{w_1} \wedge
  dw_2+f g_{\overline{w_2}} d \overline{w_1} \wedge dw_2 \wedge
  d \overline{w_2})\\
 &+(1-|w_2|^2) (f_{w_2} g dw_1 \wedge d \overline{w_1} \wedge dw_2
  +f g_{\overline{w_2}} dw_1 \wedge d \overline{w_1} \wedge
  d \overline{w_2})\\
 &-\overline{w_1} w_2 (f_{w_2} g dw_1 \wedge dw_2 \wedge d 
  \overline{w_2}+f g_{\overline{w_1}} dw_1 \wedge d \overline{w_1}
  \wedge d \overline{w_2})],
\endaligned\tag 2.18$$
where $f_{w_j}:=\frac{\partial f}{\partial w_j}$ and
$f_{\overline{w_j}}:=\frac{\partial f}{\partial \overline{w_j}}$.
Then
$$d \omega=2(1-|w_1|^2-|w_2|^2)^{-3} (w_1 d \overline{w_1}+
  \overline{w_1} dw_1+w_2 d \overline{w_2}+\overline{w_2} dw_2)
  \wedge \Omega+(1-|w_1|^2-|w_2|^2)^{-2} d \Omega.$$
Here,
$$\aligned
 &(w_1 d \overline{w_1}+\overline{w_1} dw_1+w_2 d \overline{w_2}
  +\overline{w_2} dw_2) \wedge \Omega\\
=&(1-|w_1|^2-|w_2|^2)(-w_1 f_{w_1} g+\overline{w_1} f 
  g_{\overline{w_1}}-w_2 f_{w_2} g+\overline{w_2} f 
  g_{\overline{w_2}}) dw_1 \wedge d \overline{w_1} \wedge
  dw_2 \wedge d \overline{w_2}.
\endaligned$$
$$\aligned
  d \Omega
=&[2(w_1 f_{w_1} g- \overline{w_1} f g_{\overline{w_1}}+w_2 f_{w_2} g
  -\overline{w_2} f g_{\overline{w_2}})\\
 &+(1-|w_1|^2-|w_2|^2)^{-1}(f \Delta_{{\Bbb B}^{2}} g-g 
  \Delta_{{\Bbb B}^{2}} f)] dw_1 \wedge d \overline{w_1} 
  \wedge dw_2 \wedge d \overline{w_2}.
\endaligned$$
Therefore,
$$d \omega
=(1-|w_1|^2-|w_2|^2)^{-3} (f \Delta_{{\Bbb B}^{2}} g-g 
 \Delta_{{\Bbb B}^{2}} f) dw_1 \wedge d \overline{w_1} 
 \wedge dw_2 \wedge d \overline{w_2}
=(f \Delta_{{\Bbb B}^{2}} g-g \Delta_{{\Bbb B}^{2}} f) 
  dV_{{\Bbb B}^{2}}.\tag 2.19$$
Transform the integral to an integral over ${\Bbb B}^{2}$,
exclude a small ball of radius $\epsilon$ and center $0$. We
denote the inverse image of ${\Bbb B}_{\epsilon}^{2}:=\{W \in 
{\Bbb B}^{2}: |W|=|w_1|^2+|w_2|^2 \geq \epsilon \}$ in
${\frak S}_{2}$ by ${\frak S}_{2, \epsilon}$. Let $S_{\epsilon}$
be the boundary of ${\Bbb B}_{\epsilon}^{2}$ with orientation
such that the normal is directed outside. Using spherical
coordinates, we obtain the restriction of the differential form
to the sphere $S_{\epsilon}$. We have
$dw_1=e^{i \theta_1} \cos \varphi dr+ir e^{i \theta_1} \cos 
      \varphi d \theta_1-r e^{i \theta_1} \sin \varphi d \varphi$,
$dw_2=e^{i \theta_2} \sin \varphi dr+ir e^{i \theta_2} \sin
      \varphi d \theta_2+r e^{i \theta_2} \cos \varphi d \varphi$.
Hence,
$$dm(W)=(1-|w_1|^2-|w_2|^2)^{-3} dw_1 d \overline{w_1} dw_2 
  d \overline{w_2}=4r^3 (1-r^2)^{-3} \sin \varphi \cos \varphi 
  dr d \theta_1 d \theta_2 d \varphi.\tag 2.20$$
In particular, when $r=\epsilon$ (constant),
$$\aligned
  dw_1 \wedge dw_2 \wedge d \overline{w_2}
=&-2 r^3 e^{i \theta_1} \sin \varphi \cos^2 \varphi
  d \theta_1 \wedge d \theta_2 \wedge d \varphi,\\
  dw_1 \wedge d \overline{w_1} \wedge dw_2
=&-2 r^3 e^{i \theta_2} \sin^2 \varphi \cos \varphi
  d \theta_1 \wedge d \theta_2 \wedge d \varphi.   
\endaligned$$
On the other hand,
$$\aligned
  \frac{\partial f}{\partial w_1}
=&\frac{1}{2} e^{-i \theta_1} \cos \varphi \frac{\partial f}
  {\partial r}-\frac{i}{2} \frac{e^{-i \theta_1}}{r \cos \varphi}
  \frac{\partial f}{\partial \theta_1}-\frac{e^{-i \theta_1}}{2r}
  \sin \varphi \frac{\partial f}{\partial \varphi},\\  
  \frac{\partial f}{\partial w_2}
=&\frac{1}{2} e^{-i \theta_2} \sin \varphi \frac{\partial f}
  {\partial r}-\frac{i}{2} \frac{e^{-i \theta_2}}{r \sin \varphi}
  \frac{\partial f}{\partial \theta_2}+\frac{e^{-i \theta_2}}{2r}
  \cos \varphi \frac{\partial f}{\partial \varphi}.
\endaligned$$
Therefore,
$$\aligned
   T(f)
:=&(1-|w_1|^2-|w_2|^2)^{-2} ([(1-|w_1|^2) f_{w_1}-\overline{w_1} w_2
   f_{w_2}] dw_1 \wedge dw_2 \wedge d \overline{w_2}\\
  &+[(1-|w_2|^2) f_{w_2}-w_1 \overline{w_2} f_{w_1}] dw_1 \wedge 
   d \overline{w_1} \wedge dw_2)\\
 =&-2 r^3 (1-r^2)^{-1} \sin \varphi \cos \varphi (e^{i \theta_1}
   \cos \varphi f_{w_1}+e^{i \theta_2} \sin \varphi f_{w_2}) 
   d \theta_1 \wedge d \theta_2 \wedge d \varphi\\
 =&-r^3 (1-r^2)^{-1} \sin \varphi \cos \varphi
   \left[\frac{\partial}{\partial r}-\frac{i}{r} 
   \left(\frac{\partial}{\partial \theta_1}+\frac{\partial}
   {\partial \theta_2}\right)\right] f d \theta_1 \wedge
   d \theta_2 \wedge d \varphi.
\endaligned\tag 2.21$$
In fact, 
$$\omega=T(f) g+f \overline{T(g)}.\tag 2.22$$

  In the coordinates $(r, \varphi, \theta_1, \theta_2)$, let
$f=f(r)$. By $\Delta_{{\Bbb B}^{2}} f=\lambda f$ with 
$\lambda=s(s-2)$, we get
$$\frac{1}{4} (1-r^2)^2 \frac{d^2 f}{d r^2}+\frac{(1-r^2)(3-r^2)}
  {4r} \frac{df}{dr}-\lambda f=0,$$
i.e.,
$$r^2 \frac{d^2 f}{d r^2}+r \frac{3-r^2}{1-r^2} \frac{df}{dr}
  +\frac{4 r^2}{(1-r^2)^2} s(2-s)f=0.$$
Here $\frac{3-r^2}{1-r^2}=3+2r^2+2r^4+\cdots$, i.e., $a_0=3$.
The indicial equation is $\nu(\nu-1)+3 \nu=0$.
Hence, $\nu=0$ or $\nu=-2$. Since $\frac{(1-r^2)(3-r^2)}{4r}$
has a pole at $r=0$, we have $\nu=-2$.
Let $f(r)=\sum_{n=0}^{\infty} a_n r^{n-2}$. Without loss of
generality, we can assume that $a_0=1$.  Therefore,
$$\lim_{r \to 0} r^3 f(r)=0, \quad 
  \lim_{r \to 0} r^3 f^{\prime}(r)=-2.\tag 2.23$$
Set $v=r^2$ and $f(r)=\phi(v)$, one has
$$\frac{d^2 \phi}{d v^2}+\left(\frac{2}{v}+\frac{-1}{v-1}\right)
  \frac{d \phi}{d v}+\frac{s(2-s)}{v(v-1)^2} \phi=0.\tag 2.24$$
It is known that the Fuchsian equation with three regular singularities
$a, b, c$ is given by
$$\aligned
 &\frac{d^2 w}{dz^2}+\left[ \frac{1-\alpha_1-\alpha_2}{z-a}+
  \frac{1-\beta_1-\beta_2}{z-b}+\frac{1-\gamma_1-\gamma_2}{z-c}
  \right] \frac{dw}{dz}+\frac{1}{(z-a)(z-b)(z-c)}\\
 &\times \left[ \frac{\alpha_1 \alpha_2 (a-b)(a-c)}
  {z-a}+\frac{\beta_1 \beta_2 (b-c)(b-a)}{z-b}+\frac{\gamma_1 \gamma_2 
  (c-a)(c-b)}{z-c} \right] w=0,
\endaligned$$
where $(\alpha_1, \alpha_2)$, $(\beta_1, \beta_2)$, $(\gamma_1, 
\gamma_2)$ are exponents belonging to $a, b, c$, respectively, 
and, they satisfy that 
$\alpha_1+\alpha_2+\beta_1+\beta_2+\gamma_1+\gamma_2=1$.
The solution of this equation is given in Riemann P-notation by
$$w(z)=P \left\{\matrix
         a & b & c\\
		 \alpha_1 & \beta_1 & \gamma_1\\
		 \alpha_2 & \beta_2 & \gamma_2
		 \endmatrix; z \right\}.$$
When $c=\infty$, it reduces to
$$\aligned
 &\frac{d^2 w}{d z^2}+\left[\frac{1-\alpha_1-\alpha_2}{z-a}+
  \frac{1-\beta_1-\beta_2}{z-b} \right] \frac{dw}{dz}+\frac{1}
  {(z-a)(z-b)}\\
 &\times \left[\frac{\alpha_1 \alpha_2 (a-b)}{z-a}+\frac{\beta_1 
  \beta_2 (b-a)}{z-b}+\gamma_1 \gamma_2 \right] w=0.
\endaligned$$
Now, in our case,
$$a=0, b=1, c=\infty, \alpha_1=0, \alpha_2=-1, \beta_1=s,  
  \beta_2=2-s, \gamma_1=0, \gamma_2=0,$$
and
$$P \left\{\matrix
         0 &  1  & \infty\\  
         0 &  s  & 0     \\
	    -1  & 2-s & 0                        
        \endmatrix; r^2 \right\}=(r^2-1)^{s} F(s, s; 2, r^2).$$
Now, we obtain the solution 
$$f(r)=(r^2-1)^{s} F(s, s; 2s-1; 1-r^2).\tag 2.25$$

  In the next argument, we will prove that by the isometric
transformation, the radial solution of the partial differential
equation $\Delta_{{\Bbb B}^{2}, k} f_{k}=s(s-2) f_{k}$ (see (2.32)
for the definition of $\Delta_{{\Bbb B}^{2}, k}$), i.e. 
$f_{k}=f_{k}(r)$ with $r^2=|w_1|^2+|w_2|^2$,
is transformed to the function $\phi_{s, k}(\sigma)$ which is
the solution of the differential equation
$\widetilde{\Delta_k} \phi_{s, k}(\sigma)=s(s-2) \phi_{s, k}(\sigma)$ 
(see (2.31) for the definition of $\widetilde{\Delta_k}$). 
At first, let us consider the case of weight $k=0$.
   
  By \cite{Y}, Lemma 3.2 and Proposition 3.4, we have		
$$\left[u(u+1) \frac{d^2}{d u^2}+(3u+2) \frac{d}{du}+s(2-s)\right] 
  \phi_{s}(u)=0,\tag 2.26$$		
where $u=u(Z, Z^{\prime})=\sigma(Z, Z^{\prime})-1$ and
$\delta=\delta(Z, Z^{\prime})=\frac{1}{2}[\sigma(Z, Z^{\prime})+1]$.		
On the other hand, 
$$\delta(Z, Z^{\prime})=\cosh d(Z, Z^{\prime})=\cosh 
  d_{{\Bbb B}^2}(W, 0)=\cosh \left(\log \frac{1+r}{1-r}\right)
  =\frac{1+r^2}{1-r^2}.\tag 2.27$$		
Thus, $u=2(\delta-1)=\frac{4 r^2}{1-r^2}=\frac{4v}{1-v}$,
and $\frac{d}{d v}=\frac{1}{4}(u+4)^2 \frac{d}{d u}$.
The equation (2.24) is transformed to the equation
$$\left[u(u+4) \frac{d^2}{d u^2}+(3u+8) \frac{d}{d u}+s(2-s)\right]
  \widetilde{\phi}_{s}(u)=0.$$		
Set $u=4t$, then
$$\left[t(t+1) \frac{d^2}{d t^2}+(3t+2) \frac{d}{dt}+s(2-s)\right]
  \widetilde{\widetilde{\phi}}_{s}(t)=0,$$
which is the same as (2.26).
 
  Now, let us consider the general case:
$$(\Delta_k-\lambda) g_{k, \lambda}(Z, Z_{0})=0, \quad
  Z \in {\frak S}_{2}-\{ Z_{0} \}.$$  
Set 
$$\widetilde{g_{k, \lambda}}(Z, (-\omega, 0))=\left[\frac{
  \rho(Z, (-\omega, 0))}{\rho(\overline{Z}, (-\overline{\omega},
  0))}\right]^{-k} g_{k, \lambda}(Z, Z_{0}).\tag 2.28$$
Then
$$(\widetilde{\Delta_k}-\lambda) \widetilde{g_{k, \lambda}}
  (Z, (-\omega, 0))=0,\tag 2.29$$  
where
$$\widetilde{\Delta_k}=\left[\frac{\rho(Z, (-\omega, 0))}
  {\rho(\overline{Z}, (-\overline{\omega}, 0))}\right]^{-k}
  \circ \Delta_k \circ \left[\frac{\rho(Z, (-\omega, 0))}
  {\rho(\overline{Z}, (-\overline{\omega}, 0))}\right]^{k}.
  \tag 2.30$$
In fact,
$$\aligned
 &\widetilde{\Delta_k}=\left(\frac{\overline{z_1}-\omega}{z_1-
  \overline{\omega}}\right)^{-k} \circ \Delta_{k} \circ \left(
  \frac{\overline{z_1}-\omega}{z_1-\overline{\omega}}\right)^{k}
  =\Delta+(z_1+\overline{z_1}-z_2 \overline{z_2})\\
 &\times \left[k \frac{z_1+\omega}{\overline{z_1}-\omega}
  \frac{\partial}{\partial z_1}-k \frac{\overline{z_1}+
  \overline{\omega}}{z_1-\overline{\omega}} \frac{\partial}
  {\partial \overline{z_1}}+k \frac{z_2}{\overline{z_1}-\omega}
  \frac{\partial}{\partial z_2}-k \frac{\overline{z_2}}{z_1-\omega}
  \frac{\partial}{\partial \overline{z_2}}+\frac{k^2}
  {(\overline{z_1}-\omega)(z_1-\overline{\omega})}\right].
\endaligned\tag 2.31$$
Moreover,
$$\aligned
 &\widetilde{\Delta_k} \phi_{s, k}(\sigma)=\left(\frac{
  \overline{z_1}-\omega}{z_1-\overline{\omega}}\right)^{-k}
  \Delta_{k} \varphi_{s, k}(Z, (-\omega, 0))\\
=&\left(\frac{\overline{z_1}-\omega}{z_1-\overline{\omega}}
  \right)^{-k} s(s-2) \varphi_{s, k}(Z, (-\omega, 0))
 =s(s-2) \phi_{s, k}(\sigma).
\endaligned$$

  In the ball model, the operator $\widetilde{\Delta_k}$ is
given by
$$\Delta_{{\Bbb B}^{2}, k}=\Delta_{{\Bbb B}^{2}}+(1-|w_1|^2-|w_2|^2)
  \left[k \left(w_1 \frac{\partial}{\partial w_1}-\overline{w_1} 
  \frac{\partial}{\partial \overline{w_1}}+w_2 \frac{\partial}
  {\partial w_2}-\overline{w_2} \frac{\partial}{\partial 
  \overline{w_2}}\right)+k^2\right].\tag 2.32$$
In the coordinates $(r, \varphi, \theta_1, \theta_2)$,
$$\aligned
  w_1 \frac{\partial}{\partial w_1}
=&\frac{1}{2} r \cos^2 \varphi \frac{\partial}{\partial r}-
  \frac{i}{2} \frac{\partial}{\partial \theta_1}-\frac{1}{2}
  \sin \varphi \cos \varphi \frac{\partial}{\partial \varphi},\\ 
  w_2 \frac{\partial}{\partial w_2}
=&\frac{1}{2} r \sin^2 \varphi \frac{\partial}{\partial r}-
  \frac{i}{2} \frac{\partial}{\partial \theta_2}+\frac{1}{2}
  \sin \varphi \cos \varphi \frac{\partial}{\partial \varphi}.  
\endaligned$$
Thus,
$$w_1 \frac{\partial}{\partial w_1}-\overline{w_1} \frac{\partial}
  {\partial \overline{w_1}}+w_2 \frac{\partial}{\partial w_2}-
  \overline{w_2} \frac{\partial}{\partial \overline{w_2}}=-i
  \left(\frac{\partial}{\partial \theta_1}+\frac{\partial}
  {\partial \theta_2}\right).$$

  In the coordinates $(r, \varphi, \theta_1, \theta_2)$, 
let $f_{k}=f_{k}(r)$. By $\Delta_{{\Bbb B}^2, k} f_{k}=s(s-2) f_{k}$, 
we have
$$\frac{1}{4}(1-r^2)^2 \frac{d^2 f_k}{d r^2}+\frac{(1-r^2)(3-r^2)}{4r}
  \frac{df_k}{dr}+[k^2(1-r^2)-s(s-2)] f_k=0,$$
i.e.,
$$r^2 \frac{d^2 f_{k}}{d r^2}+r \frac{3-r^2}{1-r^2} \frac{df_{k}}{dr}
  +\frac{4r^2}{(1-r^2)^2} [k^2(1-r^2)+s(2-s)] f_{k}=0.$$
The indicial equation is $\nu(\nu-1)+3 \nu=0$. Since
$\frac{(1-r^2)(3-r^2)}{4r}$ has a pole at $r=0$, $\nu=-2$.
Let $f_{k}(r)=\sum_{n=0}^{\infty} a_n r^{n-2}$ with $a_0=1$.
Then 
$$\lim_{r \to 0} r^3 f_{k}(r)=0, \quad \lim_{r \to 0} r^3 
  f_{k}^{\prime}(r)=-2.\tag 2.33$$
Set $v=r^2$ and $f_{k}(r)=\widetilde{f}_{k}(v)$, then
$$\frac{d^2 \widetilde{f}_{k}}{d v^2}+\left(\frac{2}{v}+
  \frac{-1}{v-1}\right) \frac{d \widetilde{f}_{k}}{dv}+
  \frac{1}{v(v-1)} \left[\frac{s(2-s)}{v-1}-k^2\right] 
  \widetilde{f}_{k}=0.\tag 2.34$$
Here, 
$$a=0, b=1, c=\infty, \alpha_1=0, \alpha_2=-1, \beta_1=s, 
  \beta_2=2-s, \gamma_1=|k|, \gamma_2=-|k|.$$
$$P \left\{\matrix
       0 & 1   & \infty\\
       0 & s   & |k|\\
	  -1 & 2-s & -|k| 
       \endmatrix; r^2 \right\}=(r^2-1)^{s} F(|k|+s, -|k|+s; 2; r^2).$$
Now, we get the solution
$$f_{k}(r)=(r^2-1)^{s} F(s+|k|, s-|k|; 2s-1; 1-r^2).\tag 2.35$$
By \cite{Y}, Lemma 4.8, we have
$$\left[\frac{d^2}{d \sigma^2}+\left(\frac{1}{\sigma}+\frac{2}
 {\sigma-1}\right) \frac{d}{d \sigma}+\frac{1}{\sigma(\sigma-1)}
 \left(\frac{k^2}{\sigma}+s(2-s)\right)\right] \phi_{s, k}(\sigma)
 =0.\tag 2.36$$
While, $\sigma=2 \delta-1=2 \frac{1+r^2}{1-r^2}-1=\frac{1+3r^2}
{1-r^2}=\frac{1+3v}{1-v}$, and
$\frac{d}{dv}=\frac{1}{4}(\sigma+3)^2 \frac{d}{d \sigma}$.
The equation (2.34) is transformed to the equation
$$\left\{(\sigma-1)(\sigma+3) \frac{d^2}{d \sigma^2}+[3(\sigma-1)+8]
  \frac{d}{d \sigma}+\left[s(2-s)+\frac{4k^2}{\sigma+3}\right]\right\}
  \widetilde{\phi}_{s, k}(\sigma)=0.$$
Set $\sigma=4t-3$, then
$$\left\{\frac{d^2}{d t^2}+\left(\frac{1}{t}+\frac{2}{t-1}\right)
  \frac{d}{dt}+\frac{1}{t(t-1)} \left[\frac{k^2}{t}+s(2-s)\right]
  \right\} \widetilde{\widetilde{\phi}}_{s, k}(t)=0,$$
which is the same as (2.36).

  By the above argument, we see that by the Cayley transform, 
the integral
$$\int_{{\frak S}_{2, \epsilon}} \phi_{s, k}(\sigma((-\omega, 0), Z)) 
  (\widetilde{\Delta_k}-\lambda) g(Z) dm(Z)$$  
is transformed to
$$\int_{{\Bbb B}_{\epsilon}^{2}} f_{k}(r) 
  (\Delta_{{\Bbb B}^2, k}-\lambda) h(W) dm(W).$$
Here,
$$\aligned
 &\int_{{\Bbb B}_{\epsilon}^{2}} f_{k}(r) \cdot k (1-|w_1|^2-|w_2|^2) 
  (w_1 \frac{\partial}{\partial w_1}-\overline{w_1} \frac{\partial}
  {\partial \overline{w_1}}+w_2 \frac{\partial}{\partial w_2}-
  \overline{w_2} \frac{\partial}{\partial \overline{w_2}}) 
  h(W) dm(W)\\
=&\int_{\epsilon}^{1} \int_{0}^{\frac{\pi}{2}} \int_{0}^{2 \pi} 
  \int_{0}^{2 \pi} f_{k}(r) k(1-r^2) (-i) \left(\frac{\partial}
  {\partial \theta_1}+\frac{\partial}{\partial \theta_2}\right) 
  h(r \cos \varphi e^{i \theta_1}, r \sin \varphi e^{i \theta_2})\\
 &\times 4r^3 (1-r^2)^{-3} \sin \varphi \cos \varphi d \theta_1
  d \theta_2 d \varphi dr.
\endaligned$$
Note that
$$\aligned
 &\int_{0}^{2 \pi} \frac{\partial h}{\partial \theta_1}(r \cos 
  \varphi e^{i \theta_1}, r \sin \varphi e^{i \theta_2}) d \theta_1
 =h(r \cos \varphi e^{i \theta_1}, r \sin \varphi e^{i \theta_2})
  |_{\theta_1=0}^{2 \pi}=0,\\
 &\int_{0}^{2 \pi} \frac{\partial h}{\partial \theta_2}(r \cos 
  \varphi e^{i \theta_1}, r \sin \varphi e^{i \theta_2}) d \theta_2 
 =h(r \cos \varphi e^{i \theta_1}, r \sin \varphi e^{i \theta_2})
  |_{\theta_2=0}^{2 \pi}=0.
\endaligned$$
Thus, the above integral vanishes. Hence,
$$\int_{{\Bbb B}_{\epsilon}^{2}} f_{k}(r) 
  (\Delta_{{\Bbb B}^{2}, k}-\lambda) h(W) dm(W)
 =\int_{{\Bbb B}_{\epsilon}^{2}} f_{k}(r) [\Delta_{{\Bbb B}^2}
  +k^2(1-r^2)-\lambda] h(W) dm(W).$$
Now, we have
$$\aligned
 &\int_{{\frak S}_{2}} \varphi_{s, k}((-\omega, 0), Z)
  (-\Delta_k+\lambda) g(Z) dm(Z)\\
=&\lim_{\epsilon \to 0} \int_{{\frak S}_{2, \epsilon}}
  \varphi_{s, k}((-\omega, 0), Z) (-\Delta_k+\lambda) g(Z) dm(Z),
\endaligned$$
and
$$\aligned
 &\int_{{\frak S}_{2, \epsilon}} \varphi_{s, k}((-\omega, 0), Z)
  (-\Delta_k+\lambda) g(Z) dm(Z)\\
=&\int_{{\frak S}_{2, \epsilon}} \left(\frac{\overline{z_1}-\omega}
  {z_1-\overline{\omega}}\right)^{-k} \phi_{s, k}(\sigma((-\omega, 0), 
  Z)) (-\Delta_k+\lambda) g(Z) dm(Z)\\
=&\int_{{\frak S}_{2, \epsilon}} \phi_{s, k}(\sigma) 
  (-\widetilde{\Delta_k}+\lambda) \left[\left(\frac{\overline{z_1}
  -\omega}{z_1-\overline{\omega}}\right)^{-k} g(Z)\right] dm(Z).
\endaligned$$
By the Cayley transformation, the function 
$\left(\frac{\overline{z_1}-\omega}{z_1-\overline{\omega}}\right)^{-k}
g(Z)$ is transformed to the function $h(W)$. Hence, the above integral
is equal to
$$\aligned
 &\int_{{\Bbb B}_{\epsilon}^{2}} f_{k}(r) [-\Delta_{{\Bbb B}^{2}}-
  k^2(1-r^2)+\lambda] h(W) dm(W)\\
=&\int_{{\Bbb B}_{\epsilon}^{2}} \{h(W) [\Delta_{{\Bbb B}^{2}}+
  k^2(1-r^2)] f_{k}(r)-f_{k}(r) [\Delta_{{\Bbb B}^{2}}+k^2(1-r^2)] 
  h(W)\} dm(W)\\
=&-\int_{{\Bbb B}_{\epsilon}^{2}} [f_{k}(r) \Delta_{{\Bbb B}^{2}} h(W)-
  h(W) \Delta_{{\Bbb B}^{2}} f_{k}(r)] dm(W)\\
=&-\int_{{\Bbb B}_{\epsilon}^{2}} d \omega
 =\int_{S_{\epsilon}} \omega.  
\endaligned$$
Here, $\omega=T(f_{k}(r)) h(W)+f_{k}(r) \overline{T(h(W))}$. 
In fact, 
$$\aligned
  \int_{0}^{2 \pi} \frac{\partial}{\partial \theta_j} h(W) d \theta_j
=&\int_{0}^{2 \pi} \frac{\partial}{\partial \theta_j} h(r \cos 
  \varphi e^{i \theta_1}, r \sin \varphi e^{i \theta_2}) d \theta_j\\
=&h(r \cos \varphi e^{i \theta_1}, r \sin \varphi e^{i \theta_2})
  |_{\theta_j=0}^{2 \pi}=0, \quad j=1, 2.
\endaligned$$
By (2.33), we have
$$\aligned
 \lim_{\epsilon \to 0} \int_{S_{\epsilon}} \omega
=&\lim_{\epsilon \to 0} \int_{0}^{\frac{\pi}{2}} \int_{0}^{2 \pi}
  \int_{0}^{2 \pi} [-r^3 (1-r^2)^{-1} \sin \varphi \cos \varphi
  f^{\prime}(r) h(W)\\
 &+f(r) (-r^3) (1-r^2)^{-1} \sin \varphi \cos \varphi 
  \frac{\partial}{\partial r} \overline{h(W)}] 
  d \theta_1 d \theta_2 d \varphi\\  
=&h(0) \int_{0}^{\frac{\pi}{2}} \int_{0}^{2 \pi} \int_{0}^{2 \pi}
  2 \sin \varphi \cos \varphi d \theta_1 d \theta_2 d \varphi\\
=&4 \pi^2 h(0)\\
=&4 \pi^2 \left(\frac{-\overline{\omega}-\omega}{-\omega-
  \overline{\omega}}\right)^{-k} g((-\omega, 0))\\
=&4 \pi^2 g((-\omega, 0)).     
\endaligned$$
Hence,
$$\int_{{\frak S}_{2}} \varphi_{s, k}((-\omega, 0), Z) 
  (-\Delta_k+\lambda) g(Z) dm(Z)=4 \pi^2 g((-\omega, 0)).$$
  
  In general, there exists $\gamma \in SU(2, 1)$, such that
$Z^{\prime}=\gamma(-\omega, 0)$, $Z=\gamma(\widetilde{Z})$.
$$\varphi_{s, k}(Z^{\prime}, Z)=J(\gamma, (-\omega, 0))^{2k}
  \varphi_{s, k}((-\omega, 0), \widetilde{Z}) J(\gamma, 
  \widetilde{Z})^{-2k}.$$
By \cite{Y}, Proposition 4.2, we have
$$(\Delta_k g)(Z)=(\Delta_k g)(\gamma(\widetilde{Z}))
 =J(\gamma, \widetilde{Z})^{2k} \Delta_{k}[g(\gamma(\widetilde{Z})) 
  J(\gamma, \widetilde{Z})^{-2k}].$$
Thus,
$$\aligned
 &\int_{{\frak S}_{2}} \varphi_{s, k}(Z^{\prime}, Z) 
  (-\Delta_k+\lambda) g(Z) dm(Z)\\
=&\int_{{\frak S}_{2}} J(\gamma, (-\omega, 0))^{2k} \varphi_{s, k}
  ((-\omega, 0), \widetilde{Z}) J(\gamma, \widetilde{Z})^{-2k}\\
 &\times J(\gamma, \widetilde{Z})^{2k} [-\Delta_k(\widetilde{Z})+
  \lambda] [g(\gamma(\widetilde{Z})) J(\gamma, \widetilde{Z})^{-2k}] 
  dm(\widetilde{Z})\\
=&J(\gamma, (-\omega, 0))^{2k} \int_{{\frak S}_{2}} 
  \varphi_{s, k}((-\omega, 0), \widetilde{Z}) [-\Delta_k(\widetilde{Z})
  +\lambda] [g(\gamma(\widetilde{Z})) J(\gamma, \widetilde{Z})^{-2k}] 
  dm(\widetilde{Z})\\
=&J(\gamma, (-\omega, 0))^{2k} \cdot 4 \pi^2 g(\gamma(-\omega, 0)) 
  J(\gamma, (-\omega, 0))^{-2k}\\
=&4 \pi^2 g(\gamma(-\omega, 0))=4 \pi^2 g(Z^{\prime}).
\endaligned$$

  Now, we have
$$\aligned
  4 \pi^2 g(Z^{\prime})
&=\int_{{\frak S}_{2}} \varphi_{s, k}(Z^{\prime}, Z) (-\Delta_k+
  \lambda) g(Z) dm(Z)\\
&=\sum_{\gamma \in \Gamma} \int_{\gamma({\Cal F})} \varphi_{s, k}
  (Z^{\prime}, Z) (-\Delta_k+\lambda) g(Z) dm(Z)\\
&=\sum_{\gamma \in \Gamma} \int_{{\Cal F}} \varphi_{s, k}(Z^{\prime},
  \gamma(Z)) [-\Delta_{k}(\gamma(Z))+\lambda] g(\gamma(Z)) 
  dm(\gamma(Z)).
\endaligned$$
Since $g \in C^{2}(\Gamma \backslash {\frak S}_{2}, k)$, 
one has $g(\gamma(Z))=J(\gamma, Z)^{2k} g(Z)$, and
$$\aligned
  \Delta_{k}(\gamma(Z)) g(\gamma(Z))
&=(\Delta_{k} g)(\gamma(Z))=J(\gamma, Z)^{2k} 
  \Delta_{k}[g(\gamma(Z)) J(\gamma, Z)^{-2k}]\\
&=J(\gamma, Z)^{2k} \Delta_{k} g(Z).
\endaligned$$
Hence,
$$\aligned
 &[-\Delta_{k}(\gamma(Z))+\lambda] g(\gamma(Z))
  =-J(\gamma, Z)^{2k} \Delta_k g(Z)+\lambda J(\gamma, Z)^{2k} g(Z)\\
=&J(\gamma, Z)^{2k} (-\Delta_k+\lambda) g(Z).
\endaligned$$
Therefore,
$$\aligned
  4 \pi^2 g(Z^{\prime})
&=\sum_{\gamma \in \Gamma} \int_{{\Cal F}} \varphi_{s, k}(Z^{\prime},
  \gamma(Z)) J(\gamma, Z)^{2k} (-\Delta_k+\lambda) g(Z) dm(Z)\\
&=\int_{{\Cal F}} G_{s, k}(Z^{\prime}, Z) (-\Delta_k+\lambda) g(Z)
  dm(Z), 
\endaligned$$
where
$$G_{s, k}(Z^{\prime}, Z)=\sum_{\gamma \in \Gamma} \varphi_{s, k}
  (Z^{\prime}, \gamma(Z)) J(\gamma, Z)^{2k}.$$
This completes the proof of Theorem 2.11.
\flushpar
$\qquad \qquad \qquad \qquad \qquad \qquad \qquad \qquad \qquad 
 \qquad \qquad \qquad \qquad \qquad \qquad \qquad \qquad \qquad
 \quad \square$

{\smc Proposition 2.12}. {\it Set $L[u]=\Delta_k u-s(s-2) u$. The
eigenfunctions of the associated integral equation
$$\phi(Z)-\mu \int_{{\Cal F}} G_{s, k}(Z, \xi) \phi(\xi) dm(\xi)=0
  \tag 2.37$$
are complete in $L^2(\Gamma \backslash {\frak S}_{2}, k)$.} 
 
{\it Proof}. Apply Theorem 2.11 in the context of the analogue
of the Hilbert-Schmidt theorem and the completeness theorem of
eigenfunction expansions (see \cite{Gar}, p.382-385). 
\flushpar
$\qquad \qquad \qquad \qquad \qquad \qquad \qquad \qquad \qquad
 \qquad \qquad \qquad \qquad \qquad \qquad \qquad \qquad \qquad
 \quad \square$ 
 
{\smc Proposition 2.13}. {\it Let $L[u]=\Delta_k u-s(s-2) u$. The
equations $L[u]+\mu u=0$ and
$$\phi(Z)=\frac{\mu}{4 \pi^2} \int_{{\Cal F}} G_{s, k}(Z, \xi)
  \phi(\xi) dm(\xi)\tag 2.38$$
have identical eigenfunctions for $L^2(\Gamma \backslash 
{\frak S}_{2}, k)$.} 
 
{\it Proof}. Suppose $u \in C^2(\Gamma \backslash {\frak S}_{2}, k)$ 
satisfy that $L[u]+\mu u=0$. As an application of Theorem 2.11, we 
have
$$u(Z)=\frac{\mu}{4 \pi^2} \int_{{\Cal F}} G_{s, k}(Z, \xi)
  u(\xi) dm(\xi).$$
Conversely, assume that 
$\phi \in L^2(\Gamma \backslash {\frak S}_{2}, k)$ 
satisfy $\phi=\frac{\mu}{4 \pi^2} G_{s, k} \circ \phi$.
One can modify the proof of the classical Poisson equation in
potential theory (see \cite{Gar}, p.170-173) to see that
$\phi \in C^2(\Gamma \backslash {\frak S}_{2}, k)$ and
$L[\phi]=-\mu \phi(Z)$.
\flushpar
$\qquad \qquad \qquad \qquad \qquad \qquad \qquad \qquad \qquad
 \qquad \qquad \qquad \qquad \qquad \qquad \qquad \qquad \qquad
 \quad \square$
 
{\smc Theorem 2.14}. {\it The equations $\Delta_k f=\mu f$ and
$$\phi(Z)=\frac{1}{4 \pi^2} [s(s-2)-\mu] \int_{{\Cal F}} 
  G_{s, k}(Z, \xi) \phi(\xi) dm(\xi)\tag 2.39$$
have identical eigenfunctions for $L^2(\Gamma \backslash 
{\frak S}_{2}, k)$. These eigenfunctions are complete in
$L^2(\Gamma \backslash {\frak S}_{2}, k)$.}

{\it Proof}. This is the consequence of Proposition 2.12 and
Proposition 2.13. 
\flushpar
$\qquad \qquad \qquad \qquad \qquad \qquad \qquad \qquad \qquad
 \qquad \qquad \qquad \qquad \qquad \qquad \qquad \qquad \qquad
 \quad \square$

{\it Definition} {\smc 2.15}. Let $\{ \phi_n \}_{n=0}^{\infty}$ be
the complete orthonormal family of $L^2(\Gamma \backslash 
{\frak S}_{2}, k)$ eigenfunctions generated by the Fredholm
theory of Theorem 2.14 for fixed $s$. Thus,
$$\Delta_k \phi_n=\lambda_n \phi_n, \quad
  \phi_n(Z)=\frac{1}{4 \pi^2} [s(s-2)-\lambda_n]
  \int_{{\Cal F}} G_{s, k}(Z, \xi) \phi_n(\xi) dm(\xi), n \geq 0.
  \tag 2.40$$ 
 
  Let
$$\aligned
 &L^2(\Gamma \backslash {\frak S}_{2} \times \Gamma \backslash
  {\frak S}_{2}, k)=\{F(Z, W): F \in L^2({\Cal F} \times {\Cal F}),\\
 &F(\gamma_1(Z), \gamma_2(W))=J(\gamma_1, Z)^{2k} F(Z, W)
  J(\gamma_2, W)^{-2k} \quad \text{for} \quad (\gamma_1, 
  \gamma_2) \in \Gamma \times \Gamma \}.
\endaligned$$   
$L^2(\Gamma \backslash {\frak S}_{2} \times \Gamma \backslash
{\frak S}_{2}, k)$ is a Hilbert space with inner product
$$(F, G)=\int_{{\Cal F}} \int_{{\Cal F}} F(Z, W) \overline{G(Z, W)}
         dm(Z) dm(W).$$ 
Note that
$$L^2(\Gamma \backslash {\frak S}_{2}, k)=\bigoplus_{n=0}^{\infty}
  [\phi_n(Z)],\tag 2.41$$ 
$$L^2(\Gamma \backslash {\frak S}_{2} \times \Gamma \backslash
  {\frak S}_{2}, k)=\bigoplus_{m=0}^{\infty} \bigoplus_{n=0}^{\infty}
  [\phi_{m}(Z) \overline{\phi_{n}(W)}].\tag 2.42$$  
We have
$$G_{k}(Z, W)=\sum_{n=0}^{\infty} \Lambda_{k}(\lambda_n) \phi_{n}(Z) 
              \overline{\phi_{n}(W)},\tag 2.43$$ 
where the convergence is uniform and absolute (see \cite{H}, p.12,
Proposition 3.4, p.14, Proposition 3.8, p.346, Proposition 4.8 and 
p.374), and the eigenvalues $\Lambda_{k}(\lambda_n)$ are real. 
Moreover,
$$\sum_{n=0}^{\infty} |\Lambda_{k}(\lambda_n)|<\infty.\tag 2.44$$ 
$$\int_{{\Cal F}} G_{k}(Z, Z) dm(Z)
 =\sum_{n=0}^{\infty} \Lambda_{k}(\lambda_n)
 =\text{Tr}(L_k).\tag 2.45$$
We have
$$\aligned
 &\text{Tr}(L_k)=\int_{{\Cal F}} G_{k}(Z, Z) dm(Z)\\ 
=&\int_{{\Cal F}} \sum_{\gamma \in \Gamma} K_{k}(Z, \gamma(Z))
  J(\gamma, Z)^{2k} dm(Z)\\
=&\sum_{\{ \gamma \}} \sum_{R \in \{ \gamma \}} \int_{{\Cal F}}
  K_{k}(Z, R(Z)) J(R, Z)^{2k} dm(Z).
\endaligned$$
Here the element $R$ can be expressed as $R=\sigma^{-1} \gamma 
\sigma$ with $\sigma \in \Gamma$. In fact, $\sigma^{-1} \gamma 
\sigma=\tau^{-1} \gamma \tau$ if and only if $\tau \sigma^{-1}
\gamma \sigma \tau^{-1}=\gamma$ if and only if $\sigma \tau^{-1}
\in Z(\gamma)$ if and only if $\sigma \in Z(\gamma) \tau$, where
$Z(\gamma)$ is the centralizer of $\gamma$ in $\Gamma$. Hence,
$$\text{Tr}(L_k)
 =\sum_{\{ \gamma \}} \sum_{\sigma \in Z(\gamma) \backslash \Gamma}
  \int_{{\Cal F}} K_{k}(Z, \sigma^{-1} \gamma \sigma(Z)) J(\sigma^{-1}
  \gamma \sigma, Z)^{2k} dm(Z).$$
Note that 
$$J(\gamma \sigma, Z)=J(\sigma, \sigma^{-1} \gamma 
  \sigma(Z)) J(\sigma^{-1} \gamma \sigma, Z)$$ 
and
$$K(\sigma(Z), \gamma \sigma(Z))=J(\sigma, Z)^{2k} K(Z, \sigma^{-1}
  \gamma \sigma(Z)) J(\sigma, \sigma^{-1} \gamma \sigma(Z))^{-2k},$$
we have
$$\aligned
  \text{Tr}(L_k)
=&\sum_{\{ \gamma \}} \sum_{\sigma \in Z(\gamma) \backslash \Gamma}
  \int_{{\Cal F}} K_{k}(\sigma(Z), \gamma \sigma(Z)) J(\sigma, Z)^{-2k}
  J(\sigma, \sigma^{-1} \gamma \sigma(Z))^{2k}\\
 &\times J(\gamma \sigma, Z)^{2k} J(\sigma, \sigma^{-1} \gamma 
  \sigma(Z))^{-2k} dm(Z)\\
=&\sum_{\{ \gamma \}} \sum_{\sigma \in Z(\gamma) \backslash \Gamma}
  \int_{{\Cal F}} K_{k}(\sigma(Z), \gamma \sigma(Z))
  \frac{J(\gamma \sigma, Z)^{2k}}{J(\sigma, Z)^{2k}} dm(Z)\\
=&\sum_{\{ \gamma \}} \sum_{\sigma \in Z(\gamma) \backslash \Gamma}
  \int_{{\Cal F}} K_{k}(\sigma(Z), \gamma \sigma(Z)) J(\gamma,
  \sigma(Z))^{2k} dm(Z)\\
=&\sum_{\{ \gamma \}} \sum_{\sigma \in Z(\gamma) \backslash \Gamma}
  \int_{\sigma({\Cal F})} K_{k}(\xi, \gamma(\xi)) J(\gamma, \xi)^{2k}
  dm(\xi).
\endaligned$$ 
Since ${\Cal F}$ is compact and $\Phi$ has compact support, the
$(\gamma, \sigma)$ sums above are of finite length. By the same 
method as in \cite{H}, p.23, Proposition 5.1 (Chapter 1), we can
prove that 
$\bigcup_{\sigma \in Z(\gamma) \backslash \Gamma} \sigma({\Cal F})
=\text{FR}[Z(\gamma)]$, the fundamental region for $Z(\gamma)$. 
Moreover, the integral $\int_{\text{FR}[Z(\gamma)]} K_{k}(Z, 
\gamma(Z)) J(\gamma, Z)^{2k} dm(Z)$ is independent of the choice
of fundamental region. In fact, for $\tau \in Z(\gamma)$,
$$\aligned
 &K_{k}(\tau(Z), \gamma \tau(Z)) J(\gamma, \tau(Z))^{2k}\\
=&K_{k}(\tau(Z), \tau \gamma(Z)) J(\gamma, \tau(Z))^{2k}\\
=&J(\tau, Z)^{2k} K_{k}(Z, \gamma(Z)) J(\tau, \gamma(Z))^{-2k}
  J(\gamma, \tau(Z))^{2k}\\
=&K_{k}(Z, \gamma(Z)) J(\tau, \gamma(Z))^{-2k} J(\gamma \tau, 
  Z)^{2k}\\
=&K_{k}(Z, \gamma(Z)) J(\tau, \gamma(Z))^{-2k} J(\tau \gamma,
  Z)^{2k}\\
=&K_{k}(Z, \gamma(Z)) J(\gamma, Z)^{2k}.    
\endaligned$$
Hence,
$$\text{Tr}(L_k)=\sum_{\{ \gamma \}} \int_{\text{FR}[Z(\gamma)]}
  K_{k}(Z, \gamma(Z)) J(\gamma, Z)^{2k} dm(Z).\tag 2.46$$ 
Since $\Gamma \backslash {\frak S}_{2}$ is compact and smooth, 
there are no parabolic terms and elliptic terms. We have
$$\text{Tr}(L_k)=\int_{{\Cal F}} K_{k}(Z, Z) dm(Z)+
  \sum_{\aligned
  &\text{$\{ \gamma \}$ hyperbolic}\\
  &\text{Tr}(\gamma)>2
  \endaligned} \int_{\text{FR}[Z(\gamma)]} K_{k}(Z, \gamma(Z))
  J(\gamma, Z)^{2k} dm(Z).\tag 2.47$$
 
\vskip 0.5 cm
\centerline{\bf 3. Trace formulas and zeta functions of weight $k$}
\vskip 0.5 cm

3.1. {\it Basic notions}
\vskip 0.2 cm

  According to \cite{Hi2}, let $M$ be a bounded domain in ${\Bbb C}^n$
endowed with the Bergmann hermitian metric. This is a K\"{a}hler 
metric which is invariant under complex analytic homeomorphisms
of $M$. Let $I(M)$ be the group of all such homeomorphisms and
$Y=M/\Gamma$ the quotient space defined by the action of a 
subgroup $\Gamma$ of $I(M)$. The identification map $p: M \to Y$
is a complex analytic covering map of a compact complex manifold
$Y$ if (a) $\Gamma$ is properly discontinuous, i.e. any compact 
set in $M$ intersects only a finite number of its images under 
$\Gamma$; (b) $M/\Gamma$ is compact; (c) $\Gamma$ acts freely, 
i.e. only the identity element of $\Gamma$ has fixed points.
Properties (a), (b), (c) imply that the canonical line bundle
$K_Y$ is a positive line bundle over $Y$. Therefore $Y$ is an
algebraic manifold. The following special case of a theorem of 
Borel shows that there always exist algebraic manifolds 
$M/\Gamma$ (see \cite{Hi2, Theorem 22.2.2} or \cite{Bo}).

{\smc Theorem} (Borel). {\it Let $M$ be a bounded homogeneous
symmetric domain, and $I(M)$ the group of complex analytic
homeomorphisms of $M$. Then
\roster
\item $I(M)$ contains a subgroup $\Gamma$ which satisfies (a),
      (b) and (c);
\item if $\Gamma$ is a subgroup of $I(M)$ which satisfies (a)
      and (b), and which does not consist only of the identity
	  element, then $\Gamma$ has a proper normal subgroup of
	  finite index which satisfies (a), (b) and (c).
\endroster} 

  Let $M$ be a complete locally symmetric Riemannian manifold of
negative curvature. Assume that the real dimension of $M$ is even, 
$M$ has finite volume and $M$ is compact, since the complex
hyperbolic space ${\Bbb H}_{\Bbb C}^{n}$ is homogeneous, the
Gauss-Bonnet formula or Hirzebruch proportionality theorem
\cite{Hi2, Theorem 22.2.1} tells us that there is a constant 
$\kappa_{n}$ such that $\text{Vol}(M)=\kappa_n e(M)$, where
$e(M)=\sum_{i=0}^{2n} (-1)^{i} \dim H^i(M, {\Bbb R})$ with
$\dim_{{\Bbb R}} M=2n$ is the Euler number of $M$. When 
the holomorphic sectional curvature is normalized to be $-1$, hence 
the sectional curvature is between $-\frac{1}{4}$ and $-1$, we have 
(see \cite{HP}) $$\kappa_n=\frac{(-\pi)^n 2^{2n}}{(n+1)!}.$$
In particular, $\kappa_2=\frac{8}{3} \pi^2$.

  Yau's Theorem \cite{Yau} tells us that all varieties $X$ on which
$K$ is ample carry a unique K\"{a}hler-Einstein metric. This
result shows that if $X$ is a surface on which $K$ is ample, then
the Chern numbers satisfy $c_1^2 \leq 3 c_2$, the equality holds
if and only if $X$ is isomorphic to $\Gamma \backslash {\Bbb B}^2$,
where $\Gamma \subset PSU(2, 1)$ is a discrete torsion-free cocompact
subgroup. Hirzebruch \cite{Hi1} showed that the surfaces 
$\Gamma \backslash {\Bbb B}^2$ did satisfy $c_1^2=3c_2$.   
Furthermore, Borel showed that (see \cite{BPV}, p.177)
for infinitely many $a \in {\Bbb N}$ there exists a surface $X$ 
which has the unit ball in ${\Bbb C}^{2}$ as its universal 
covering, such that $c_1^2(X)=3a$, $c_2(X)=a$.

  In our case, since $M=\Gamma \backslash {\frak S}_{2}$ is covered 
by the unit ball, $c_1^2=3 c_2$. Noether's formula tells us that 
(see \cite{BPV}, p.20) 
$$\chi(X)=\sum_{i=0}^{2} (-1)^{i} \dim H^{i}(X, {\Cal O}_{X})
 =1-q(X)+p_g(X)=\frac{1}{12}(c_1^2(X)+c_2(X)),$$
where $X$ is a compact complex surface. Hence
$c_1^2+c_2 \equiv 0 (\text{mod $12$})$. 
Consequently, $c_2(M) \equiv 0 (\text{mod $3$})$.  
Therefore, $\text{Vol}(M)=\frac{8}{3} \pi^2 c_2(M)$ is a 
multiple of $\pi^2$.
  
\vskip 0.4 cm
3.2. {\it The computation of the identity} 
\vskip 0.2 cm

  We have
$$\int_{{\Cal F}} K_{k}(Z, Z) dm(Z)
 =\int_{{\Cal F}} H_{k}(Z, Z) \Phi(0) dm(Z)
 =\Phi(0) \text{Vol}({\Cal F}),$$
where ${\Cal F}$ is the fundamental region of $\Gamma$.
  
  In the case of weight $k=0$, by
$\Phi(u)=\frac{1}{\pi} \int_{u}^{\infty} \frac{1}{\sqrt{v-u}}
         d(\sqrt{v+1} P^{\prime}(v))$,
we have
$$\Phi(0)=\frac{1}{\pi} \int_{0}^{\infty} \frac{1}{\sqrt{v}}
         d(\sqrt{v+1} P^{\prime}(v)).$$
Since $v=\sinh^2 \frac{\eta}{2}$, one has
$v+1=\cosh^2 \frac{\eta}{2}$, $dv=\frac{1}{2} \sinh \eta d \eta$.
$P(v)=g(\eta)$ implies that $P^{\prime}(v) dv=g^{\prime}(\eta) d 
\eta$. So $P^{\prime}(v)=g^{\prime}(\eta) \frac{d \eta}{d v}=
\frac{2 g^{\prime}(\eta)}{\sinh \eta}$. Thus, 
$\sqrt{v+1} P^{\prime}(v)=\frac{g^{\prime}(\eta)}
      {\sinh \frac{\eta}{2}}$,
$$\Phi(0)=\frac{1}{\pi} \int_{0}^{\infty} \frac{1}{\sinh \frac{\eta}{2}}
       d \left(\frac{g^{\prime}(\eta)}{\sinh \frac{\eta}{2}}\right).$$
While,
$g(\eta)=\frac{1}{4 \pi^2} \int_{-\infty}^{\infty} h(r)
 e^{-ir \eta} dr$ and $h(r)=h(-r)$,
we have
$$g^{\prime}(\eta)=-\frac{1}{4 \pi^2} \int_{-\infty}^{\infty}
           r h(r) \sin r \eta dr.$$
Hence,
$$\Phi(0)=\frac{1}{4 \pi^3} \int_{-\infty}^{\infty} r h(r) dr
  \int_{0}^{\infty} \frac{\frac{1}{2} \cosh \frac{\eta}{2} 
  \sin r \eta-r \sinh \frac{\eta}{2} \cos r \eta}{\sinh^3
  \frac{\eta}{2}} d \eta.$$
Denote the inner integral as $F(r)$. Set $u=\frac{1}{2} \eta$ and 
$G(2r)=F(r)$, then
$$G(r)=\int_{0}^{\infty} \frac{\cosh u \sin ru-r \sinh u \cos ru}
       {\sinh^3 u} du.$$
We have $G(r)=\frac{\pi}{4} r^2 \frac{1+e^{-\pi r}}{1-e^{-\pi r}}$
(see (A.1) in Appendix). So $F(r)=G(2r)=\frac{\pi r^2}{\tanh \pi r}$.
Hence, $\Phi(0)=\frac{1}{4 \pi^2} \int_{-\infty}^{\infty} h(r)
\frac{r^3}{\tanh \pi r} dr$ and
$$\text{Vol}(M) \Phi(0)=\frac{2}{3} c_2(M) \int_{-\infty}^{\infty} 
  h(r) \frac{r^3}{\tanh \pi r} dr.\tag 3.1$$

  In the case of weight $k=1$, by 
$\Phi(u)=\frac{u+1}{\pi} \int_{u}^{\infty}
 \frac{1}{\sqrt{v-u}} d \left[ \left(\frac{P_1(v)}{\sqrt{v+1}}
 \right)^{\prime}\right]$,  
we have
$$\Phi(0)=\frac{1}{\pi} \int_{0}^{\infty} \frac{1}{\sqrt{v}}
  \left(\frac{P_1(v)}{\sqrt{v+1}}\right)^{\prime \prime} dv.$$
Put $v=\sinh^2 \frac{u}{2}$, then $\sqrt{v}=\sinh \frac{u}{2}$,
$\sqrt{v+1}=\cosh \frac{u}{2}$, $P_1(v)=g_1(u)$ and
$\frac{dv}{du}=\frac{1}{2} \sinh u$.
We have $P_1^{\prime}(v) dv=g_1^{\prime}(u) du$, so
$P_1^{\prime}(v)=\frac{2 g_1^{\prime}(u)}{\sinh u}$. 
Similarly,
$P_1^{\prime \prime}(v) dv=\left(\frac{2 g_1^{\prime}(u)}{\sinh u}
 \right)^{\prime} du$, this implies that
$P_1^{\prime \prime}(v)=\frac{2}{\sinh u} \left(\frac{2 
 g_1^{\prime}(u)}{\sinh u}\right)^{\prime}$.
Here
$$\left(\frac{P_1(v)}{\sqrt{v+1}}\right)^{\prime \prime}
 =\frac{P_1^{\prime \prime}(v)}{\sqrt{v+1}}-\frac{P_1^{\prime}(v)}
  {(v+1)^{\frac{3}{2}}}+\frac{3}{4} 
  \frac{P_1(v)}{(v+1)^{\frac{5}{2}}}.$$
Thus,
$$\Phi(0)=\frac{1}{\pi} \int_{0}^{\infty} \left[\frac{4 g_1^{\prime 
  \prime}(u)}{\sinh^2 u}-\frac{4 \cosh u}{\sinh^3 u} g_1^{\prime}(u)
  -\frac{2 g_1^{\prime}(u)}{\sinh u \cosh^2 \frac{u}{2}}+\frac{3}{4}
  \frac{g_1(u)}{\cosh^4 \frac{u}{2}}\right] du.$$
Note that
$g_1(u)=\frac{1}{2 \pi^2} \int_{0}^{\infty} h_{1}(r) \cos ru dr$,
we have
$$\aligned
  \Phi(0)
&=\frac{1}{2 \pi^3} \int_{0}^{\infty} h_{1}(r) dr\\
&\times \int_{0}^{\infty}
  \left[4r \frac{\cosh u \sin ru-r \cos ru \sinh u}{\sinh^3 u}+2r
  \frac{\sin ru}{\sinh u \cosh^2 \frac{u}{2}}+\frac{3}{4}
  \frac{\cos ru}{\cosh^4 \frac{u}{2}} \right] du.
\endaligned$$
Here, by (A.1) in Appendix, we have
$$\int_{0}^{\infty} 4r \frac{\cosh u \sin ru-r \sinh u \cos ru}
  {\sinh^3 u} du=\frac{\pi r^3}{\tanh \frac{\pi}{2} r}.$$
On the other hand,
$2r \int_{0}^{\infty} \frac{\sin ru}{\sinh u \cosh^2 \frac{u}{2}} du
 =r \int_{0}^{\infty} \frac{\sin ru}{\sinh \frac{u}{2} \cosh^3
  \frac{u}{2}} du$.
By (A.3) and (A.4) in Appendix, we obtain
$$\int_{0}^{\infty} \frac{\sin r u}{\sinh \frac{u}{2} \cosh^3
  \frac{u}{2}} du=\frac{\pi}{\tanh \pi r}-\frac{\pi (2r^2+1)}
  {\sinh \pi r}, \quad
  \int_{0}^{\infty} \frac{\cos ru}{\cosh^4 \frac{u}{2}} du=
  \frac{4 \pi r (r^2+1)}{3 \sinh \pi r}.$$
Therefore,
$$\Phi(0)=\frac{1}{2 \pi^3} \int_{0}^{\infty} h_{1}(r) 
  \left[\frac{\pi r^3}{\tanh \frac{\pi}{2} r}+\frac{\pi r}
  {\tanh \pi r}-\frac{\pi r (2r^2+1)}{\sinh \pi r}+
  \frac{\pi r (r^2+1)}{\sinh \pi r}\right] dr.$$
By $\frac{1}{\tanh \frac{\pi}{2} r}-\frac{1}{\sinh \pi r}=
   \frac{1}{\tanh \pi r}$, we have
$\Phi(0)=\frac{1}{2 \pi^2} \int_{0}^{\infty} h_{1}(r) \frac{r(r^2+1)}
 {\tanh \pi r} dr$. Thus,
$$\text{Vol}(M) \Phi(0)=\frac{2}{3} c_2(M) \int_{-\infty}^{\infty} 
  h_{1}(r) \frac{r(r^2+1)}{\tanh \pi r} dr.\tag 3.2$$  

  In the case of weight $k=\frac{1}{2}$, by
$$\Phi(u)=\frac{1}{\pi} \sqrt{u+1} \int_{u}^{\infty}
  \frac{d Q_0^{\prime}(v)}{\sqrt{v-u}},$$
we have
$$\Phi(0)=\frac{1}{\pi} \int_{0}^{\infty} \frac{1}{\sqrt{v}}
  d Q_0^{\prime}(v).$$
By the same method as above, one has
$$\aligned
 &\Phi(0)=\frac{1}{\pi} \int_{0}^{\infty} \frac{1}{\sinh 
  \frac{\eta}{2}} d \left(\frac{2 g_{\frac{1}{2}}^{\prime}(\eta)}
  {\sinh \eta} \right)\\
=&\frac{1}{2 \pi^3} \int_{-\infty}^{\infty} r h_{\frac{1}{2}}(r) dr
  \int_{0}^{\infty} \frac{\cosh \eta \sin r \eta-r \sinh \eta 
  \cos r \eta}{\sinh \frac{\eta}{2} \sinh^2 \eta} d \eta. 
\endaligned$$
Denote the inner integral as $H(r)$, we have (see (A.2) in Appendix)
$$H(r)=\frac{\pi}{2} \left(r^2+\frac{1}{4}\right) \tanh \pi r.$$
Therefore,
$$\Phi(0)=\frac{1}{4 \pi^2} \int_{-\infty}^{\infty} 
  h_{\frac{1}{2}}(r) r \left(r^2+\frac{1}{4}\right) \tanh \pi r dr$$
and
$$\text{Vol}(M) \Phi(0)=\frac{2}{3} c_2(M) \int_{-\infty}^{\infty}
  h_{\frac{1}{2}}(r) r \left(r^2+\frac{1}{4}\right) \tanh \pi r dr.
  \tag 3.3$$
  
\vskip 0.4 cm
3.3. {\it The computation of hyperbolic elements}
\vskip 0.2 cm

  Suppose that the hyperbolic element takes the form  
$\gamma=\text{diag}(\mu e^{i \theta}, e^{-2i \theta},                           
        \mu^{-1} e^{i \theta})$.  
For $Z=(z_1, z_2) \in {\frak S}_{2}$,   
$\gamma(z_1, z_2)=(\mu^2 z_1,\mu e^{-3i \theta} z_2)$,  
$\rho(\gamma(z_1, z_2))=\mu^2 \rho(z_1, z_2)$.
By Proposition 2.1, $Z(\gamma)=<\gamma_1> \times <\gamma_2>$, 
where $(\gamma_1, \gamma_2)$ is a primitive pair and 
$\gamma_1=\text{diag}(r_0, 1, r_0^{-1})$, 
$\gamma_2=\text{diag}(e^{i \varphi_0}, e^{-2i \varphi_0},
 e^{i \varphi_0})$. Thus, 
$$\text{FR}[Z(\gamma)]=\text{FR}[<\gamma_1> \times <\gamma_2>]
 =\text{FR}(\gamma_1) \cap \text{FR}(\gamma_2).$$
One has $\gamma_2(z_1, z_2)=(z_1, e^{-3 i \varphi_0} z_2)$ and 
$\gamma_1(z_1, z_2)=(r_0^2 z_1, r_0 z_2)$.
Hence,
$$\text{FR}(\gamma_2)=\{ (z_1, z_2) \in {\frak S}_{2}:
  0 \leq \arg(z_2) < \frac{6 \pi}{\text{ord}(\gamma_2)} \},$$
$$\text{FR}(\gamma_1)=\{ (z_1, z_2) \in {\frak S}_{2}:
  1 \leq \rho(z_1, z_2) < r_0^2 \}.$$
Therefore,
$$\text{FR}[Z(\gamma)]=\{ (z_1, z_2) \in {\frak S}_{2}:
  1 \leq \rho(z_1, z_2) < r_0^2, \quad 0 \leq \arg(z_2) <
  \frac{6 \pi}{\text{ord}(\gamma_2)} \}.$$
On the other hand,
$$\sigma(Z, \gamma(Z))=\frac{|\rho(Z, \gamma(Z))|^2}{\rho(Z) 
  \rho(\gamma(Z))}
 =\mu^{-2} \rho^{-2} |\overline{z_1}+\mu^2 z_1-\mu e^{-3 i \theta} 
   z_2 \overline{z_2}|^2.$$  
We have
$$\aligned
 &c(\gamma)=\int_{\text{FR}[Z(\gamma)]} \Phi(u(Z, \gamma(Z))) dm(Z)\\
=&\int_{\text{FR}[Z(\gamma)]} \Phi(\rho^{-2} \mu^{-2} |\overline{z_1}+
  \mu^2 z_1-\mu e^{-3 i \theta} z_2 \overline{z_2}|^2-1) \rho^{-3} 
  dz_1 d \overline{z_1} dz_2 d \overline{z_2}. 
\endaligned$$
Set $z_1=\frac{1}{2}(\rho+|z|^2)+it$, $z_2=z$ where $\rho>0$, 
$t \in {\Bbb R}$ and $z \in {\Bbb C}$. Then
$$\aligned
 &|\overline{z_1}+\mu^2 z_1-\mu e^{-3 i \theta} z_2 \overline{z_2}|^2\\
=&|[\frac{1}{2}(\rho+|z|^2)(\mu^2+1)-\mu |z|^2 \cos 3 \theta]+
  i[(\mu^2-1)t+\mu |z|^2 \sin 3 \theta]|^2\\
=&\left[\frac{1}{2}(\mu^2+1)(\rho+|z|^2)-\mu |z|^2 \cos 3 \theta
  \right]^2+\left[(\mu^2-1)t+\mu |z|^2 \sin 3 \theta \right]^2.
\endaligned$$
The norm of $\gamma$ is $N(\gamma)=\mu^2$. Therefore,
$$\aligned
  c(\gamma)
=&\int_{1}^{N(\gamma_1)} \int_{-\infty}^{\infty} \int_{z \in {\Bbb C},
  0 \leq \arg(z) < \frac{6 \pi}{\text{ord}(\gamma_2)}}
  \Phi(\rho^{-2} ([\frac{1}{2}(\mu+\mu^{-1})(\rho+|z|^2)
  -|z|^2 \cos 3 \theta]^2\\
 &+[(\mu-\mu^{-1})t+|z|^2 \sin 3 \theta]^2)-1)
  \rho^{-3} dz d \overline{z} dt d \rho\\
=&\frac{6 \pi}{\text{ord}(\gamma_2)} \int_{1}^{N(\gamma_1)} 
  \int_{-\infty}^{\infty} \int_{0}^{\infty} 
  \Phi(\rho^{-2}([\frac{1}{2}(\mu+\mu^{-1})
  (\rho+r)-r \cos 3 \theta]^2\\
 &+[(\mu-\mu^{-1})t+r \sin 3\theta]^2)-1) \rho^{-3} dr dt d \rho. 
\endaligned$$
Set $u=\frac{r}{\rho}$ and $v=\frac{2t}{\rho}$, then
$$\aligned
  c(\gamma)
=&\frac{3 \pi}{\text{ord}(\gamma_2)} \ln N(\gamma_1) 
  \int_{-\infty}^{\infty} \int_{0}^{\infty} 
  \Phi([\frac{1}{2}(\mu+\mu^{-1})(u+1)-u \cos 3 \theta]^2\\
 &+[\frac{1}{2}(\mu-\mu^{-1})v+u \sin 3 \theta]^2-1) du dv.
\endaligned$$
Set $f(u)=\frac{1}{2}(\mu+\mu^{-1})(u+1)-u \cos 3 \theta$, then
$f^{\prime}(u)=\frac{1}{2}(\mu+\mu^{-1})-\cos 3 \theta \geq 1-
 \cos 3 \theta \geq 0$. So $f(u) \geq f(0)=\frac{1}{2}(\mu+\mu^{-1})$.
Let
$$\xi=\left[\frac{1}{2}(\mu+\mu^{-1})(u+1)-u \cos 3 \theta \right]^2-1,
  \quad
  \eta=\frac{1}{2}(\mu-\mu^{-1})v+u \sin 3 \theta.$$
Then $\xi \in [\frac{1}{4}(\mu-\mu^{-1})^2, \infty)$,
$\eta \in (-\infty, \infty)$.
$$d \xi=\sqrt{\xi+1} (\mu+\mu^{-1}-2 \cos 3 \theta) du, \quad
  d \eta=\frac{1}{2} (\mu-\mu^{-1}) dv+\sin 3 \theta du.$$
Then
$$du \wedge dv
 =\frac{2}{(\mu+\mu^{-1}-2 \cos 3 \theta)(\mu-\mu^{-1})}
  \frac{d \xi \wedge d \eta}{\sqrt{\xi+1}}.$$
Hence,
$$c(\gamma)=\frac{6 \pi \cdot \text{ord}(\gamma_2)^{-1} 
  \ln N(\gamma_1)}{(\mu-\mu^{-1})(\mu+\mu^{-1}-2 \cos 3 \theta)} 
  \int_{\frac{1}{4}(\mu-\mu^{-1})^2}^{\infty} \int_{-\infty}^{\infty} 
  \Phi(\xi+\eta^2) d \eta \frac{d \xi}{\sqrt{\xi+1}}.$$
Set $v=\eta^2$, then
$$c(\gamma)=\frac{6 \pi \cdot \text{ord}(\gamma_2)^{-1} 
  \ln N(\gamma_1)}{(\mu-\mu^{-1})(\mu+\mu^{-1}-2 \cos 3 \theta)} 
  \int_{\frac{1}{4}(\mu-\mu^{-1})^2}^{\infty} \int_{0}^{\infty} 
  \Phi(\xi+v) \frac{dv}{\sqrt{v}} \frac{d \xi}{\sqrt{\xi+1}}.$$
Here
$$\aligned
 &\int_{\frac{1}{4}(\mu-\frac{1}{\mu})^2}^{\infty} \int_{0}^{\infty}
  \Phi(\xi+\eta) \frac{d \eta}{\sqrt{\eta}} \frac{d \xi}{\sqrt{\xi+1}}
=\int_{\frac{1}{4}(\mu-\frac{1}{\mu})^2}^{\infty} \int_{\xi}^{\infty}
  \Phi(w) \frac{d w}{\sqrt{w- \xi}} \frac{d \xi}{\sqrt{\xi+1}}\\
=&\int_{\frac{1}{4}(\mu-\frac{1}{\mu})^2}^{\infty} \Phi(w) dw
  \int_{\frac{1}{4}(\mu-\frac{1}{\mu})^2}^{\infty}
  \frac{d \xi}{\sqrt{(w- \xi)(\xi+1)}}\\
=&\int_{\frac{1}{4}(\mu-\frac{1}{\mu})^2}^{\infty} \Phi(w)
  \arccos\left(\frac{\frac{1}{2}(\mu-\frac{1}{\mu})^2-w+1}{w+1}
  \right)dw\\    
=&P \left(\frac{1}{4} \left(\mu-\frac{1}{\mu}\right)^2 \right)
=g(\eta), 
\endaligned$$
where $\frac{1}{4}(\mu-\frac{1}{\mu})^2=\frac{1}{4} (e^{\eta}+
e^{-\eta}-2)$, i.e., $\eta=\ln \mu^2=\ln N(\gamma)$.
Note that $\mu+\mu^{-1}-2 \cos 3 \theta=|\mu^{\frac{1}{2}} 
e^{\frac{3}{2} i \theta}-\mu^{-\frac{1}{2}} e^{-\frac{3}{2} i 
\theta}|^2$. Therefore,
$$c(\gamma)=\frac{6 \pi \cdot \text{ord}(\gamma_2)^{-1} 
  \ln N(\gamma_1)}{[N(\gamma)^{\frac{1}{2}}-N(\gamma)^{-\frac{1}{2}}] 
  |N(\gamma)^{\frac{1}{4}} e^{\frac{3}{2} i \theta(\gamma)}-
  N(\gamma)^{-\frac{1}{4}} e^{-\frac{3}{2} i \theta(\gamma)}|^2} 
  g(\ln N(\gamma)).\tag 3.4$$

  In the case of weight $k=1$, set
$$I_{k}(\gamma)=\int_{\text{FR}[Z(\gamma)]} H_{k}(Z, \gamma(Z))
  \Phi(u(Z, \gamma(Z))) J(\gamma, Z)^{2k} dm(Z).$$
Then
$$\aligned
 &I_{k}(\gamma)=\\
 &\int_{1}^{N(\gamma_1)} \int_{-\infty}^{\infty} 
  \int_{z \in \Bbb C, 0 \leq \arg(z) < \frac{6 \pi}{\text{ord}
  (\gamma_2)}} H_{k}(Z, \gamma(Z)) \Phi(u(Z, \gamma(Z))) 
  J(\gamma, Z)^{2k} \rho^{-3} d \rho dt dz d \overline{z}.
\endaligned$$
Note that $J(\gamma, Z)=e^{i \theta}$.
Using the same notation as above, by the transform 
$z=\sqrt{r} e^{i \vartheta}$ with $r>0$, 
$\vartheta \in [0, 2 \pi)$, $u=\frac{r}{\rho}$ and 
$v=\frac{2 t}{\rho}$, we have
$$\aligned
 &I_{k}(\gamma)=\frac{3 \pi}{\text{ord}(\gamma_2)} \ln N(\gamma_1)
  e^{2ki \theta}\\
 &\times \int_{-\infty}^{\infty} \int_{0}^{\infty}
  \frac{([\frac{1}{2}(\mu+\mu^{-1})(u+1)-u \cos 3 \theta]+i[\frac{1}{2}
  (\mu-\mu^{-1}) v+u \sin 3 \theta])^{2k}}{|[\frac{1}{2}(\mu+\mu^{-1})
  (u+1)-u \cos 3 \theta]+i[\frac{1}{2} (\mu-\mu^{-1}) v+u 
  \sin 3 \theta]|^{2k}}\\
 &\times \Phi([\frac{1}{2}(\mu+\mu^{-1})(u+1)-u \cos 3 \theta]^2
  +[\frac{1}{2}(\mu-\mu^{-1})v+u \sin 3 \theta]^2-1) du dv.
\endaligned$$
Let $\xi=[\frac{1}{2}(\mu+\mu^{-1})(u+1)-u \cos 3 \theta]^2-1$, 
$\eta=\frac{1}{2}(\mu-\mu^{-1})v+u \sin 3 \theta$ as above,
$$\aligned
  I_{k}(\gamma)
=&\frac{6 \pi \cdot \text{ord}(\gamma_2)^{-1} \ln N(\gamma_1) 
  e^{2ki \theta}}{(\mu-\mu^{-1})(\mu+\mu^{-1}-2 \cos 3 \theta)} 
  \int_{-\infty}^{\infty} \int_{\frac{1}{4}(\mu-\mu^{-1})^2}^{\infty}
  \frac{(\sqrt{\xi+1}+i \eta)^{2k}}{|\sqrt{\xi+1}+i \eta|^{2k}}\\
 &\times \Phi(\xi+\eta^2) \frac{1}{\sqrt{\xi+1}} d \xi d \eta. 
\endaligned$$
By
$$(\sqrt{\xi+1}+i \eta)^{2k}=\sum_{m=0}^{2k} \binom{2k}{m} 
  \sqrt{\xi+1}^{2k-m} i^m \eta^m, \quad
  |\sqrt{\xi+1}+i \eta|^{2k}=(\xi+\eta^2+1)^k,$$
we have
$$\aligned
  I_{k}(\gamma)
=&\frac{6 \pi \cdot \text{ord}(\gamma_2)^{-1} \ln N(\gamma_1) 
  e^{2ki \theta}}{(\mu-\mu^{-1})(\mu+\mu^{-1}-2 \cos 3 \theta)} 
  \int_{-\infty}^{\infty} \int_{\frac{1}{4}(\mu-\mu^{-1})^2}^{\infty} 
  \sum_{m=0}^{k} \binom{2k}{2m} (-1)^{m}\\
 &\times (\xi+1)^{k-m-\frac{1}{2}} \eta^{2m} (\xi+\eta^2+1)^{-k} 
  \Phi(\xi+\eta^2) d \xi d \eta.
\endaligned$$
By the transform $v=\eta^2$, one has
$$\aligned
  I_{k}(\gamma)
=&\frac{6 \pi \cdot \text{ord}(\gamma_2)^{-1} \ln N(\gamma_1) 
  e^{2ki \theta}}{(\mu-\mu^{-1})(\mu+\mu^{-1}-2 \cos 3 \theta)} 
  \int_{0}^{\infty} \int_{\frac{1}{4}(\mu-\mu^{-1})^2}^{\infty} 
  \sum_{m=0}^{k} \binom{2k}{2m} (-1)^m\\
 &\times (\xi+1)^{k-m-\frac{1}{2}} v^{m-\frac{1}{2}} (\xi+v+1)^{-k} 
  \Phi(\xi+v) d \xi dv. 
\endaligned$$
Here,
$$\aligned
 &\int_{0}^{\infty} \int_{\frac{1}{4}(\mu-\mu^{-1})^2}^{\infty}
  \sum_{m=0}^{k} \binom{2k}{2m} (-1)^{m} \frac{(\xi+1)^{k-m-
  \frac{1}{2}} \eta^{m-\frac{1}{2}}}{(\xi+\eta+1)^k} \Phi(\xi+\eta) 
  d \xi d \eta\\ 
=&\int_{\frac{1}{4}(\mu-\mu^{-1})^2}^{\infty} \int_{\xi}^{\infty}
  \sum_{m=0}^{k} \binom{2k}{2m} (-1)^{m} (\xi+1)^{k-m-\frac{1}{2}} 
  (w-\xi)^{m-\frac{1}{2}} (w+1)^{-k} \Phi(w) dw d \xi\\ 
=&\int_{\frac{1}{4}(\mu-\mu^{-1})^2}^{\infty} \Phi(w) (w+1)^{-k} dw
  \int_{\frac{1}{4}(\mu-\mu^{-1})^2}^{w} \sum_{m=0}^{k} \binom{2k}{2m}
  (-1)^{m}\\
 &\times (\xi+1)^{k-m- \frac{1}{2}} (w- \xi)^{m-\frac{1}{2}} d \xi. 
\endaligned$$
Note that
$$P_{k}(v)=\int_{v}^{\infty} \Phi(w) (w+1)^{-k} dw
           \int_{v}^{w} \sum_{m=0}^{k} (-1)^{m} \binom{2k}{2m}
		   (\xi+1)^{k-m-\frac{1}{2}} (w- \xi)^{m- \frac{1}{2}} d \xi.$$
Then
$$\aligned
 &I_{k}(\gamma)\\
=&\frac{6 \pi \cdot \text{ord}(\gamma_2)^{-1} \ln N(\gamma_1) 
  e^{2ki \theta(\gamma)}}{[N(\gamma)^{\frac{1}{2}}-
  N(\gamma)^{-\frac{1}{2}}] |N(\gamma)^{\frac{1}{4}} 
  e^{\frac{3}{2} i \theta(\gamma)}-N(\gamma)^{-\frac{1}{4}} 
  e^{-\frac{3}{2} i \theta(\gamma)}|^{2}} P_{k} 
  (\frac{1}{4}[N(\gamma)^{\frac{1}{2}}-N(\gamma)^{-\frac{1}{2}}]^{2}).
\endaligned\tag 3.5$$
In fact, $\gamma$ can be expressed as $\gamma_1^m \gamma_2^q$, 
where $(\gamma_1, \gamma_2)$ is a primitive pair (see Proposition
2.1). Thus,
$$c(\gamma)=\frac{6 \pi \cdot \text{ord}(\gamma_2)^{-1} 
  \ln N(\gamma_1)}{[N(\gamma_1)^{\frac{m}{2}}
  -N(\gamma_1)^{-\frac{m}{2}}] |N(\gamma_1)^{\frac{m}{4}}
  e^{\frac{3}{2} i q \theta(\gamma_2)}-N(\gamma_1)^{-\frac{m}{4}}
  e^{-\frac{3}{2} i q \theta(\gamma_2)}|^2} g(m \ln N(\gamma_1)).$$
Similar expression can be done for $I_{k}(\gamma)$.
$$I_{k}(\gamma)=\frac{6 \pi \cdot \text{ord}(\gamma_2)^{-1} 
  \ln N(\gamma_1) e^{2kiq \theta(\gamma_2)}}
  {[N(\gamma_1)^{\frac{m}{2}}-N(\gamma_1)^{-\frac{m}{2}}] 
  |N(\gamma_1)^{\frac{m}{4}} e^{\frac{3}{2} i q \theta(\gamma_2)}
  -N(\gamma_1)^{-\frac{m}{4}} e^{-\frac{3}{2} i q 
  \theta(\gamma_2)}|^2} g_{k}(m \ln N(\gamma_1)).$$

\vskip 0.4 cm
3.4. {\it The trace formulas of weight $k$}
\vskip 0.2 cm

From now on, we need the following assumption. 
\roster
\item $h_{k}(r)$ is an analytic function on 
      $|\text{Im}(r)| \leq 1+\delta$.
\item $h_{k}(-r)=h_{k}(r)$.
\item $|h_{k}(r)| \leq M [1+|\text{Re}(r)|]^{-4-\delta}$.
\endroster
Here $\delta$ and $M$ are positive constants. 
 
  We claim that the Fourier transform $g_{k}(\eta)=\frac{1}{4 \pi^2} 
\int_{-\infty}^{\infty} h_{k}(r) e^{-ir \eta} dr$, $\eta \in 
{\Bbb R}$ satisfies the inequality  
$$|g_{k}(\eta)| \leq \frac{M}{2 \pi^2} e^{-(1+\delta) |\eta|}.$$  
In fact, suppose that $\eta<0$. An application of the Cauchy
integral theorem shows that  
$$\int_{-\infty}^{\infty} h_{k}(r) e^{-ir \eta} dr
 =\int_{\text{Im}(r)=1+\delta} h_{k}(r) e^{-ir \eta} dr.$$  
The right hand side is bounded in absolute value by
$$\aligned
 &\int_{-\infty}^{\infty} |h(x+i(1+\delta))| e^{(1+\delta) \eta} dx
  \leq e^{(1+\delta) \eta} \int_{-\infty}^{\infty} 
  \frac{M}{(1+|x|)^{4+\delta}} dx\\
\leq & 2 e^{(1+\delta) \eta} \int_{0}^{\infty} \frac{M}{(1+x)^2} dx
      =2 M e^{(1+\delta) \eta}.
\endaligned$$
Note that
$$|N(\gamma_1)^{\frac{m}{4}} e^{\frac{3}{2}iq \theta(\gamma_2)}-
  N(\gamma_1)^{-\frac{m}{4}} e^{-\frac{3}{2}iq \theta(\gamma_2)}|^2
  \geq [N(\gamma_1)^{\frac{m}{4}}-N(\gamma_1)^{-\frac{m}{4}}]^2.$$
We have
$$\aligned
 &\sum_{\{(\gamma_1, \gamma_2)\}} \sum_{q=1}^{\text{ord}(\gamma_2)} 
  \sum_{m=1}^{\infty} \frac{\text{ord}(\gamma_2)^{-1} \cdot 
  \ln N(\gamma_1) \cdot g_{k}(m \ln N(\gamma_1))}
  {[N(\gamma_1)^{\frac{m}{2}}-N(\gamma_1)^{-\frac{m}{2}}]
   |N(\gamma_1)^{\frac{m}{4}} e^{\frac{3}{2}iq \theta(\gamma_2)}-
   N(\gamma_1)^{-\frac{m}{4}} e^{-\frac{3}{2}iq \theta(\gamma_2)}|^2}\\
\leq & \sum_{\{(\gamma_1, \gamma_2)\}} \sum_{m=1}^{\infty} 
  \frac{\ln N(\gamma_1) \cdot g_{k}(m \ln N(\gamma_1))}
  {[N(\gamma_1)^{\frac{m}{2}}-N(\gamma_1)^{-\frac{m}{2}}]
   [N(\gamma_1)^{\frac{m}{4}}-N(\gamma_1)^{-\frac{m}{4}}]^2}.
\endaligned$$
As $m \to \infty$, $N(\gamma_1)^{\frac{m}{4}}-N(\gamma_1)^{-\frac{m}
{4}} \to \infty$. Hence, there exists $m_0 \in {\Bbb N}$, such that
when $m \geq m_0$, $N(\gamma_1)^{\frac{m}{4}}-N(\gamma_1)^{-\frac{m}
{4}} \geq 1$. Thus,
$$\aligned
 &\sum_{m=m_0}^{\infty} \frac{\ln N(\gamma_1) \cdot 
  g_{k}(m \ln N(\gamma_1))}{[N(\gamma_1)^{\frac{m}{2}}-
  N(\gamma_1)^{-\frac{m}{2}}][N(\gamma_1)^{\frac{m}{4}}-
  N(\gamma_1)^{-\frac{m}{4}}]^{2}}\\
\leq &\sum_{m=m_0}^{\infty} \frac{\ln N(\gamma_1)}{N(\gamma_1)^{
  \frac{m}{2}}-N(\gamma_1)^{-\frac{m}{2}}} \frac{M}{2 \pi^2}
  e^{-(1+\delta)m \ln N(\gamma_1)},
\endaligned$$
which is absolutely convergent (see \cite{H}, p.31, Proposition 7.4).
Now, we get the Main Theorem 1.

{\smc Theorem 3.1}. {\it Assume that $h_{k}(r)$ is an 
analytic function on $|\text{Im}(r)| \leq 1+\delta$, 
$h_{k}(-r)=h_{k}(r)$ and $|h_{k}(r)| \leq M [1+|\text{Re}
(r)|]^{-4-\delta}$, where $\delta$ and $M$ are positive 
constants. Moreover, suppose that  
$$g_{k}(\eta)=\frac{1}{4 \pi^2} \int_{-\infty}^{\infty}
  h_{k}(r) e^{-ir \eta} dr, \quad \eta \in {\Bbb R}.$$ 
Then the trace formulas of weight $k$ for $SU(2, 1)$ take 
the following form:
$$\aligned
  \text{Tr}(L_{k})=
 &\sum_{n=0}^{\infty} h_{k}(r_n)
 =\frac{2}{3} c_2(M) \int_{-\infty}^{\infty} 
  h_{k}(r) r(r^2+k^2) \frac{\cosh 2 \pi r+\cos 2 k \pi}
  {\sinh 2 \pi r} dr+\\
 &+2 \pi \sum_{\{(\gamma_1, \gamma_2)\}} 
  \frac{3}{\text{ord}(\gamma_2)} \sum_{q=1}^{\text{ord}
  (\gamma_2)} \sum_{m=1}^{\infty} \frac{\ln N(\gamma_1)}
  {[N(\gamma_1)^{\frac{m}{2}}-N(\gamma_1)^{-\frac{m}{2}}]}\\
 &\times \frac{e^{2kiq \theta(\gamma_2)}}{|N(\gamma_1)^{\frac{m}{4}} 
  e^{\frac{3}{2} i q \theta(\gamma_2)}-N(\gamma_1)^{-\frac{m}{4}} 
  e^{-\frac{3}{2} i q \theta(\gamma_2)}|^2} g_{k}(m \ln N(\gamma_1)).
\endaligned\tag 3.6$$
where $k=0, \frac{1}{2}, 1$ and $h_0=h$, $g_0=g$.}		  

  $M=\Gamma \backslash {\frak S}_{2}$ with $\Gamma$ a uniform,
torsion-free discrete subgroup of $G=SU(2, 1)$, is a compact complex
algebraic surface. It is known that (see \cite{S} or \cite{M}), 
on $S(M)$, the unit tangent bundle of $M$, we have geodesic 
flow $\varphi_t$. 
If $(x, \xi) \in S(M)$ is the origin and tangent vector of a closed 
geodesic $\gamma$ starting at $x$ with tangent $\xi$ and of length 
$L$, then $(x, \xi)$ is a fixed point of $\varphi_{L}$. So 
$d \varphi_{L}: T(S(M), (x, \xi)) \to T(S(M), (x, \xi))$. 
It preserves the geodesic flow vector field
and hence induces a transformation, the Poincar\'{e} map, 
$P_{\gamma}$ normal to that direction, i.e., $P_{\gamma}$ is a
linear transformation on a vector space of dimension $2n-2$ if
$\dim_{{\Bbb R}} M=n$. For symmetric spaces $P_{\gamma}$ can be 
computed using the Jacobi equation: For $P_{\gamma}$ there is the 
map: $(V(0), V^{\prime}(0)) \to (V(L), V^{\prime}(L))$ where $V$ 
is any Jacobi field along $\gamma$, normal to $\xi$. In our case,
$\dim_{{\Bbb R}} M=4$, $P_{\gamma}$ is a linear transformation on 
a vector space of dimension $6$.

  Given $A, B \in M_{m}({\Bbb C})$ we write $A \sim B$ if $A$ and
$B$ are similar; that is, there exists $Q \in M_{m}({\Bbb C})$ such
that $QAQ^{-1}=B$. For $\gamma=\text{diag}(e^{u} e^{i \theta}, 
e^{-2i \theta}, e^{-u} e^{i \theta})$, with 
$N(\gamma)=\mu^2=e^{2 u}$, $u>0$.
$\overline{\gamma}(z_1, z_2)=(e^{2u} z_1, e^{u} e^{3i \theta} z_2)$.
The associated Poincar\'{e} map $P_{\gamma}$:
$$P_{\gamma} \sim \left(\matrix
  e^{2u} &         &          &       \\
 &\left(\matrix
   e^{u} e^{3i \theta} & *                     \\	  
   0                   & e^{u} e^{3i \theta}
  \endmatrix\right)  &          &          \\   
 &        & e^{-2u}   &        \\         
 &        &           &  \left(\matrix
        e^{-u} e^{-3i \theta} & *                       \\  
        0                     & e^{-u} e^{-3i \theta} 
       \endmatrix\right)
 \endmatrix\right).$$
Hence,
$$\aligned
  |\det(I-P_{\gamma})|
=&|(1-e^{2u})(1-e^{-2u})| \cdot |(1-e^{u} e^{3i \theta})
  (1-e^{-u} e^{-3i \theta})|^2\\
=&|e^{u}-e^{-u}|^2 \cdot |e^{\frac{1}{2} u} e^{\frac{3}{2}
  i \theta}-e^{-\frac{1}{2} u} e^{-\frac{3}{2} i \theta}|^4.
\endaligned$$
Thus, $|\det(I-P_{\gamma})|^{\frac{1}{2}}=(e^{u}-e^{-u})
|e^{\frac{1}{2} u} e^{\frac{3}{2} i \theta}-e^{-\frac{1}{2} 
u} e^{-\frac{3}{2} i \theta}|^2$.

  Using the Poincar\'{e} map, we can give another formulation 
of the trace formula as follows.
  
{\smc Theorem 3.2}. {\it The other formulation of the trace 
formulas of weight $k$ for $SU(2, 1)$ takes the form:
$$\aligned
 &\text{Tr}(L_{k})=\frac{2}{3} c_2(M) \int_{-\infty}^{\infty} 
  h_{k}(r) r(r^2+k^2) \frac{\cosh 2 \pi r+\cos 2 k \pi}
  {\sinh 2 \pi r} dr+\\
 &+2 \pi \sum_{\{(\gamma_1, \gamma_2)\}} \frac{3}{\text{ord}(\gamma_2)} 
  \sum_{q=1}^{\text{ord}(\gamma_2)} \sum_{m=1}^{\infty} 
  \frac{\ln N(\gamma_1)}{|\det(I-P_{\gamma_1^{m} 
  \gamma_2^{q}})|^{\frac{1}{2}}} e^{2kiq \theta(\gamma_2)} 
  g_{k}(m \ln N(\gamma_1)),
\endaligned\tag 3.7$$
where $k=0, \frac{1}{2}, 1$ and $h_0=h$, $g_0=g$.}
 
\vskip 0.4 cm
3.5. {\it The zeta functions of weight $k$}
\vskip 0.2 cm

{\it Definition} {\smc 3.3}. The zeta function of weight $k$ 
associated with the group $<\gamma_1> \times <\gamma_2>$ 
generated by $\gamma_1$ and $\gamma_2$, where $<\gamma_1>$ 
is an infinite cyclic group, $<\gamma_2>$ is a finite cyclic 
group, is given by
$$\aligned
 &Z_{<\gamma_1> \times <\gamma_2>}(k, s)\\
=&\prod_{j=0}^{\infty} \prod_{l=0}^{\infty} \prod_{n=0}^{\infty} 
  \left[1-N(\gamma_1)^{-s-j-\frac{l+n}{2}}\right]
  ^{\frac{3}{\text{ord}(\gamma_2)} \sum_{q=1}^{\text{ord}(\gamma_2)} 
  \exp([2k+3(l-n)]iq \theta(\gamma_2))},
\endaligned\tag 3.8$$
for $\text{Re}(s)>2$ and $2k \in {\Bbb Z}$.
  
  For $\Gamma$ a uniform, torsion-free discrete subgroup of 
$SU(2, 1)$, we define
$$Z_{k}(s)=Z_{\Gamma}(k, s)=\prod_{\{(\gamma_1, \gamma_2)\}} 
  Z_{<\gamma_1> \times <\gamma_2>}(k, s), \tag 3.9$$
where $\{(\gamma_1, \gamma_2)\}$ runs through a set of 
representatives of primitive conjugacy classes of $\Gamma$.

{\smc Proposition 3.4}. {\it The zeta functions of weight $k$ 
satisfy the following relations:
$$\aligned
 &\frac{Z_{k}^{\prime}(s)}{Z_{k}(s)}=\sum_{\{(\gamma_1, \gamma_2)\}} 
  \frac{3}{\text{ord}(\gamma_2)}
  \sum_{q=1}^{\text{ord}(\gamma_2)} \sum_{m=1}^{\infty} 
  \frac{\ln N(\gamma_1)}{[N(\gamma_1)^{\frac{m}{2}}-N(\gamma_1)
  ^{-\frac{m}{2}}]}\\
 &\times \frac{e^{2kiq \theta(\gamma_2)}}{|N(\gamma_1)^{\frac{m}{4}} 
  e^{\frac{3}{2} i q \theta(\gamma_2)}-N(\gamma_1)^{-\frac{m}{4}} 
  e^{-\frac{3}{2} i q \theta(\gamma_2)}|^2} 
  \frac{1}{N(\gamma_1)^{m(s-1)}}.
\endaligned\tag 3.10$$}

{\it Proof}. We have
$$\aligned
S:=&\sum_{q=1}^{\text{ord}(\gamma_2)} \sum_{m=1}^{\infty} 
  \frac{3 \cdot \text{ord}(\gamma_2)^{-1} e^{2kiq \theta(\gamma_2)}
  \ln N(\gamma_1) \cdot N(\gamma_1)^{-m(s-1)}}{[N(\gamma_1)^{
  \frac{m}{2}}-N(\gamma_1)^{-\frac{m}{2}}] |N(\gamma_1)^{\frac{m}{4}} 
  e^{\frac{3}{2} i q \theta(\gamma_2)}-N(\gamma_1)^{-\frac{m}{4}} 
  e^{-\frac{3}{2} i q \theta(\gamma_2)}|^2}\\
=&\sum_{q=1}^{\text{ord}(\gamma_2)} \sum_{m=1}^{\infty} 
  \frac{3 \cdot \text{ord}(\gamma_2)^{-1} e^{2kiq \theta(\gamma_2)} 
  \ln N(\gamma_1) \cdot N(\gamma_1)^{-ms}}{[1-N(\gamma_1)^{-m}]
  [1-N(\gamma_1)^{-\frac{1}{2} m} e^{3iq \theta(\gamma_2)}] 
  [1-N(\gamma_1)^{-\frac{1}{2} m} e^{-3iq \theta(\gamma_2)}]}\\ 
=&\frac{3}{\text{ord}(\gamma_2)} \sum_{q=1}^{\text{ord}(\gamma_2)} 
  \sum_{m=1}^{\infty} \ln N(\gamma_1) e^{2kiq \theta(\gamma_2)} 
  N(\gamma_1)^{-ms} \sum_{j=0}^{\infty} N(\gamma_1)^{-mj}\\
 &\times \sum_{l=0}^{\infty} N(\gamma_1)^{-\frac{m}{2} l}
  e^{3iql \theta(\gamma_2)} \sum_{n=0}^{\infty} 
  N(\gamma_1)^{-\frac{m}{2} n} e^{-3iqn \theta(\gamma_2)}.
\endaligned$$
A straightforward calculation gives
$$\aligned  
S=&\frac{3}{\text{ord}(\gamma_2)} \ln N(\gamma_1) 
  \sum_{j=0}^{\infty} \sum_{l=0}^{\infty} 
  \sum_{n=0}^{\infty} \sum_{q=1}^{\text{ord}(\gamma_2)} 
  e^{[2k+3(l-n)]iq \theta(\gamma_2)} \sum_{m=1}^{\infty} 
  N(\gamma_1)^{-m(s+j+\frac{l+n}{2})}\\  
=&\frac{3}{\text{ord}(\gamma_2)} \ln N(\gamma_1) 
  \sum_{j=0}^{\infty} \sum_{l=0}^{\infty} 
  \sum_{n=0}^{\infty} \frac{N(\gamma_1)^{-(s+j+\frac{l+n}{2})}}
  {1-N(\gamma_1)^{-(s+j+\frac{l+n}{2})}} 
  \sum_{q=1}^{\text{ord}(\gamma_2)} e^{[2k+3(l-n)]iq 
  \theta(\gamma_2)}\\
=&\frac{3}{\text{ord}(\gamma_2)} \sum_{j=0}^{\infty} 
  \sum_{l=0}^{\infty} \sum_{n=0}^{\infty} 
  \frac{d}{ds}\left\{\log \left[1-N(\gamma_1)^{-s-j-\frac{l+n}{2}}
  \right]\right\} \sum_{q=1}^{\text{ord}(\gamma_2)} 
  e^{[2k+3(l-n)]iq \theta(\gamma_2)}\\
=&\frac{d}{ds} \log \left\{\prod_{j=0}^{\infty} \prod_{l=0}^{\infty}
  \prod_{n=0}^{\infty} \left[1-N(\gamma_1)^{-s-j-\frac{l+n}{2}}
  \right]^{\frac{3}{\text{ord}(\gamma_2)} 
  \sum_{q=1}^{\text{ord}(\gamma_2)} 
  \exp([2k+3(l-n)]iq \theta(\gamma_2))}\right\}. 
\endaligned$$
By $Z_{k}(s)=\prod_{\{(\gamma_1, \gamma_2)\}} Z_{<\gamma_1> 
\times <\gamma_2>}(k, s)$, we have
$\frac{Z_{k}^{\prime}(s)}{Z_{k}(s)}=\sum_{\{(\gamma_1, \gamma_2)\}} 
 \frac{Z_{<\gamma_1> \times <\gamma_2>}^{\prime}(k, s)}
 {Z_{<\gamma_1> \times <\gamma_2>}(k, s)}$.
\flushpar 
$\qquad \qquad \qquad \qquad \qquad \qquad \qquad \qquad \qquad
 \qquad \qquad \qquad \qquad \qquad \qquad \qquad \qquad \qquad
 \quad \square$

  Using the Poincar\'{e} map, we can reformulate the Proposition 
3.4 as follows.

{\smc Theorem 3.5}. {\it The expansions of zeta functions of 
weight $k$ for $SU(2, 1)$ are expressed as:
$$\frac{Z_{k}^{\prime}(s)}{Z_{k}(s)}=\sum_{\{(\gamma_1, \gamma_2)\}} 
  \frac{3}{\text{ord}(\gamma_2)} \sum_{q=1}^{\text{ord}(\gamma_2)} 
  \sum_{m=1}^{\infty} \frac{\ln N(\gamma_1)}{|\det(I-P_{\gamma_1^{m} 
  \gamma_2^{q}})|^{\frac{1}{2}}} \frac{e^{2kiq \theta(\gamma_2)}}
  {N(\gamma_1)^{m(s-1)}}.\tag 3.11$$}

\vskip 0.4 cm
3.6. {\it The functional equations of zeta functions}
\vskip 0.2 cm

  Now, we derive the functional equations of zeta functions of
weight $k$. In fact, by the trace formulas which we develop as 
above, there are three zeta functions, which we denote by $Z_0(s)$, 
$Z_{\frac{1}{2}}(s)$ and $Z_1(s)$ according to the weight $k=0$,
$k=\frac{1}{2}$ or $k=1$. We will prove that all of them 
satisfy the analogue of the Riemann hypothesis. Put 
$$\aligned
 &g_1(\eta)=g_{\frac{1}{2}}(\eta)=g(\eta)\\
=&\frac{1}{4 \pi (\alpha_1^2-\alpha_2^2)}
  \left(\frac{1}{\alpha_2} e^{-\alpha_2 |\eta|}-\frac{1}{\alpha_1}
  e^{-\alpha_1 |\eta|}\right)-\frac{1}{4 \pi (\beta_1^2-\beta_2^2)} 
  \left(\frac{1}{\beta_2} e^{-\beta_2 |\eta|}-\frac{1}{\beta_1}
  e^{-\beta_1 |\eta|}\right)
\endaligned$$
for $\eta \in {\Bbb R}$ and $1 < \text{Re}(\alpha_1) < 
\text{Re}(\alpha_2)<\text{Re}(\beta_1)<\text{Re}(\beta_2)$.
By the following formula 
$$\int_{-\infty}^{\infty} e^{-a |x|} e^{irx} dx=\frac{2a}{r^2+a^2},$$
we obtain
$$h_1(r)=h_{\frac{1}{2}}(r)=h(r)=\frac{1}{(r^2+\alpha_1^2)
  (r^2+\alpha_2^2)}-\frac{1}{(r^2+\beta_1^2)(r^2+\beta_2^2)}.$$
Note that
$$h_1(r)=h_{\frac{1}{2}}(r)=h(r)=\frac{(\beta_1^2+\beta_2^2-
  \alpha_1^2-\alpha_2^2) r^2+(\beta_1^2 \beta_2^2-\alpha_1^2 
  \alpha_2^2)}{(r^2+\alpha_1^2)(r^2+\alpha_2^2)(r^2+\beta_1^2)
  (r^2+\beta_2^2)}=O(r^{-6}), \quad r \to \infty.$$
The trace formulas yield
$$\aligned
 &\sum_{n=0}^{\infty} \left[\frac{1}{(r_n^2+\alpha_1^2)
  (r_n^2+\alpha_2^2)}-\frac{1}{(r_n^2+\beta_1^2)(r_n^2+\beta_2^2)}
  \right]\\
=&\frac{2}{3} c_2(M) \int_{-\infty}^{\infty} 
  \left[\frac{1}{(r^2+\alpha_1^2)(r^2+\alpha_2^2)}-\frac{1}
  {(r^2+\beta_1^2)(r^2+\beta_2^2)}\right] \frac{\cosh 2 \pi r
  +\cos 2k \pi}{\sinh 2 \pi r}\\
 &\times r(r^2+k^2) dr+\frac{1}{2(\alpha_1^2-\alpha_2^2)} 
  \sum_{\{(\gamma_1, \gamma_2)\}} \frac{3}{\text{ord}(\gamma_2)} 
  \sum_{q=1}^{\text{ord}(\gamma_2)} \sum_{m=1}^{\infty} 
  \frac{\ln N(\gamma_1)}{[N(\gamma_1)^{\frac{m}{2}}
  -N(\gamma_1)^{-\frac{m}{2}}]}\\
 &\times \frac{e^{2kiq \theta(\gamma_2)}}{|N(\gamma_1)^{\frac{m}{4}} 
  e^{\frac{3}{2} i q \theta(\gamma_2)}-N(\gamma_1)^{-\frac{m}{4}} 
  e^{-\frac{3}{2} i q \theta(\gamma_2)}|^2}
  \left[\frac{1}{\alpha_2} \frac{1}{N(\gamma_1)^{m \alpha_2}}-
  \frac{1}{\alpha_1} \frac{1}{N(\gamma_1)^{m \alpha_1}}\right]\\  
 &-\frac{1}{2(\beta_1^2-\beta_2^2)} \sum_{\{(\gamma_1, \gamma_2)\}} 
  \frac{3}{\text{ord}(\gamma_2)} \sum_{q=1}^{\text{ord}(\gamma_2)} 
  \sum_{m=1}^{\infty} \frac{\ln N(\gamma_1)}{[N(\gamma_1)^{\frac{m}{2}}
  -N(\gamma_1)^{-\frac{m}{2}}]}\\
 &\times \frac{e^{2kiq \theta(\gamma_2)}}{|N(\gamma_1)^{\frac{m}{4}}
  e^{\frac{3}{2} i q \theta(\gamma_2)}-N(\gamma_1)^{-\frac{m}{4}} 
  e^{-\frac{3}{2} i q \theta(\gamma_2)}|^2}
  \left[\frac{1}{\beta_2} \frac{1}{N(\gamma_1)^{m \beta_2}}-
  \frac{1}{\beta_1} \frac{1}{N(\gamma_1)^{m \beta_1}}\right].   
\endaligned$$
Set $\alpha_1=s-1$ with $\text{Re}(s)>2$, this gives
$$\aligned
 &\sum_{n=0}^{\infty} \left[\frac{1}{\alpha_2^2-(s-1)^2} \left(\frac{1}
  {r_n^2+(s-1)^2}-\frac{1}{r_n^2+\alpha_2^2}\right)-\frac{1}
  {\beta_2^2-\beta_1^2} \left(\frac{1}{r_n^2+\beta_1^2}-
  \frac{1}{r_n^2+\beta_2^2}\right)\right]\\
=&\frac{2}{3} c_2(M) \int_{-\infty}^{\infty} 
  [\frac{1}{\alpha_2^2-(s-1)^2} \left(\frac{1}{r^2+(s-1)^2}-
  \frac{1}{r^2+\alpha_2^2}\right)\\
 &-\frac{1}{\beta_2^2-\beta_1^2} \left(\frac{1}{r^2+\beta_1^2}-
  \frac{1}{r^2+\beta_2^2}\right)] r(r^2+k^2) \frac{\cosh 2 \pi r
  +\cos 2k \pi}{\sinh 2 \pi r} dr\\
 &+\frac{1}{2[(s-1)^2-\alpha_2^2]} \left[\frac{1}{\alpha_2}
  \frac{Z_{k}^{\prime}(\alpha_2+1)}{Z_{k}(\alpha_2+1)}-\frac{1}{s-1}
  \frac{Z_{k}^{\prime}(s)}{Z_{k}(s)}\right]\\
 &-\frac{1}{2(\beta_1^2-\beta_2^2)} \left[\frac{1}{\beta_2}
  \frac{Z_{k}^{\prime}(\beta_2+1)}{Z_{k}(\beta_2+1)}-\frac{1}{\beta_1}
  \frac{Z_{k}^{\prime}(\beta_1+1)}{Z_{k}(\beta_1+1)}\right]. 
\endaligned$$
Therefore, we have the following result.

{\smc Proposition 3.6}. {\it For $\text{Re}(s)>1$, the following
identities hold:
$$\aligned
 &\frac{1}{s-1} \frac{Z_{k}^{\prime}(s)}{Z_{k}(s)}\\
=&\frac{1}{\alpha_2} \frac{Z_k^{\prime}(\alpha_2+1)}{Z_k(\alpha_2+1)}-
  \frac{(s-1)^2-\alpha_2^2}{\beta_1^2-\beta_2^2} \left[\frac{1}
  {\beta_2} \frac{Z_k^{\prime}(\beta_2+1)}{Z_k(\beta_2+1)}-\frac{1}
  {\beta_1} \frac{Z_k^{\prime}(\beta_1+1)}{Z_k(\beta_1+1)}\right]\\
&+\frac{4}{3} c_2(M) \int_{-\infty}^{\infty} 
 [-\left(\frac{1}{r^2+(s-1)^2}-\frac{1}{r^2+\alpha_2^2}\right)\\
&+\frac{(s-1)^2-\alpha_2^2}{\beta_1^2-\beta_2^2} \left(
  \frac{1}{r^2+\beta_1^2}-\frac{1}{r^2+\beta_2^2}\right)]
  r(r^2+k^2) \frac{\cosh 2 \pi r+\cos 2k \pi}{\sinh 2 \pi r} dr\\
&-2 \sum_{n=0}^{\infty} \left[-\left(\frac{1}{r_n^2+(s-1)^2}-
  \frac{1}{r_n^2+\alpha_2^2}\right)+\frac{(s-1)^2-\alpha_2^2}
  {\beta_1^2-\beta_2^2} \left(\frac{1}{r_n^2+\beta_1^2}-
  \frac{1}{r_n^2+\beta_2^2}\right)\right],  
\endaligned\tag 3.12$$
where $k=0, \frac{1}{2}, 1$.}

  The preceding theorem allows us to continue 
$\frac{Z^{\prime}(s)}{Z(s)}$ into $\text{Re}(s)<1$.

  Let us consider the $r_n$ sum first. It is trivial to 
check the absolute convergence away from the points 
$s_n=1+i r_n$, $\widetilde{s_n}=1-i r_n$. The singularity 
which occurs at $s_n$, $\widetilde{s_n}$ is given as follows:
\roster
\item $r_n \neq 0$. $s_n$ contributes $\frac{1}{2s_n-2} 
 \frac{1}{s-s_n}+O(1)$; $\widetilde{s_n}$ contributes
 $\frac{1}{2 \widetilde{s_n}-2} \frac{1}{s- \widetilde{s_n}}+O(1)$.  
\item $r_n=0$. $s_n=1$ contributes $\frac{1}{(s-1)^2}+O(1)$.
\endroster

  Let
$$\aligned
  W_{k}(\xi)
=&\int_{-\infty}^{\infty} \left[\frac{1}{r^2+\alpha_2^2}-
  \frac{1}{r^2+\xi^2}+\frac{\alpha_2^2-\xi^2}{\beta_2^2-
  \beta_1^2} \left(\frac{1}{r^2+\beta_1^2}-\frac{1}{r^2
  +\beta_2^2}\right) \right]\\
 &\times r(r^2+k^2) \frac{\cosh 2 \pi r+\cos 2k \pi}
  {\sinh 2\pi r} dr.
\endaligned$$
We will give the analytic continuation of $W_{k}(\xi)$ into 
$\text{Re}(\xi)<0$. Without loss of generality, we can only
consider the case of $k=0$. 
		   
$$\frac{1}{\tanh \pi r}=\frac{e^{\pi r}+e^{-\pi r}}{e^{\pi r}-
  e^{-\pi r}}=\frac{1+e^{-2 \pi r}}{1-e^{-2 \pi r}}$$
is an odd function with period $i$. Let $r=r_1+i r_2$, $r_1>0$,
$r_2 \in {\Bbb R}$, we have 
$$\frac{1}{\tanh \pi r}=1+O(e^{-2 \pi r_1}) \quad \text{as}
  \quad r_1 \to \infty.$$
The poles of $\frac{1}{\tanh \pi r}$ are simple and are located
at the points $k i$, $k \in {\Bbb Z}$.

  We will move the path of integration in $W(\xi)$ to
$\text{Im}(r)=N- \frac{1}{2}$, where $N$ is a large integer.
We assume here that $0<\text{Re}(\xi)<2$.

  Singularities are encountered for $r=\pm i \xi$, $\pm i \alpha_2$,
$\pm i \beta_1$, $\pm i \beta_2$ and $k i$ $(k \in {\Bbb Z})$. The
Cauchy residue theorem yields
$$\aligned
 &\int_{-\infty}^{\infty} [\cdots] \frac{r^3}{\tanh \pi r} dr
 =\int_{-\infty+i(N- \frac{1}{2})}^{\infty+i(N- \frac{1}{2})}
  [\cdots] \frac{r^3}{\tanh \pi r} dr+2 \pi i \sum_{0 \leq k < N} 
  \text{Res}[r=ki]\\
 &+2 \pi i \text{Res}[r=i \xi]+2 \pi i \text{Res}[r=i \alpha_2]
  +2 \pi i \text{Res}[r=i \beta_1]+2 \pi i \text{Res}[r=i \beta_2].
\endaligned$$
Here
$$\text{Res}[r=k i]=\frac{R_k^3}{\pi} \left[\frac{1}{R_k^2+
  \alpha_2^2}-\frac{1}{R_k^2+\xi^2}+\frac{\alpha_2^2-\xi^2}
  {\beta_2^2-\beta_1^2} \left(\frac{1}{R_k^2+\beta_1^2}-
  \frac{1}{R_k^2+\beta_2^2}\right)\right],$$ 
where $R_k=ki$.
$$\text{Res}[r=i \xi]=-\frac{i}{2} \frac{\xi^2}{\tan \pi \xi},
  \quad
  \text{Res}[r=i \alpha_2]=\frac{i}{2} \frac{\alpha_2^2}
  {\tan \pi \alpha_2}.$$
$$\text{Res}[r=i \beta_1]=\frac{i}{2} \frac{\beta_1^2}
  {\tan \pi \beta_1} \frac{\alpha_2^2-\xi^2}{\beta_2^2-\beta_1^2},
  \quad
  \text{Res}[r=i \beta_2]=-\frac{i}{2} \frac{\beta_2^2}
  {\tan \pi \beta_2} \frac{\alpha_2^2-\xi^2}{\beta_2^2-\beta_1^2}.$$
Thus,
$$\aligned
 &\int_{-\infty}^{\infty} [\cdots] \frac{r^3}{\tanh \pi r} dr
 =\int_{-\infty+i(N-\frac{1}{2})}^{\infty+i(N-\frac{1}{2})}
  [\cdots] \frac{r^3}{\tanh \pi r} dr\\
 &+2i \sum_{0 \leq k<N} R_k^3 \left[\frac{1}{R_k^2+\alpha_2^2}-
  \frac{1}{R_k^2+\xi^2}+\frac{\alpha_2^2-\xi^2}{\beta_2^2-\beta_1^2}
  \left(\frac{1}{R_k^2+\beta_1^2}-\frac{1}{R_k^2+\beta_2^2}\right) 
  \right]\\ 
 &+\pi \left(\frac{\xi^2}{\tan \pi \xi}-\frac{\alpha_2^2}
  {\tan \pi \alpha_2}\right)-\pi \frac{\alpha_2^2-\xi^2}{\beta_2^2
  -\beta_1^2} \left(\frac{\beta_1^2}{\tan \pi \beta_1}-
  \frac{\beta_2^2}{\tan \pi \beta_2}\right).
\endaligned$$
The analytic continuation of $W(\xi)$ is obtained.

  If $|\text{Re}(\xi)|<N-\frac{1}{2}$, we can define
$$\aligned
  W_0(\xi)
=&\int_{-\infty+i(N-\frac{1}{2})}^{\infty+i(N-\frac{1}{2})}
  [\cdots] \frac{r^3}{\tanh \pi r} dr\\
 &+2 i \sum_{0 \leq k<N} R_k^3 \left[\frac{1}{R_k^2+\alpha_2^2}
  -\frac{1}{R_k^2+\xi^2}+\frac{\alpha_2^2-\xi^2}{\beta_2^2-\beta_1^2}
  \left(\frac{1}{R_k^2+\beta_1^2}-\frac{1}{R_k^2+\beta_2^2}\right) 
  \right]\\
 &+\pi (\xi^2 \cot \pi \xi-\alpha_2^2 \cot \pi \alpha_2)    
  -\pi \frac{\alpha_2^2-\xi^2}{\beta_2^2-\beta_1^2} (\beta_1^2
  \cot \pi \beta_1-\beta_2^2 \cot \pi \beta_2).
\endaligned$$
The $i(N-\frac{1}{2})$ integral is holomorphic for 
$|\text{Re}(\xi)|<N-\frac{1}{2}$.

  Let $u=r-i(N- \frac{1}{2})$, note that
$\tanh [\pi(u+i(N- \frac{1}{2}))]=\frac{1}{\tanh \pi u}$, we have
$$\aligned
 &\int_{-\infty+i(N- \frac{1}{2})}^{\infty+i(N- \frac{1}{2})} 
  [\cdots] \frac{r^3}{\tanh \pi r} dr\\ 
=&\int_{-\infty}^{\infty} [\frac{1}{(u+i(N-\frac{1}{2}))^2+\alpha_2^2}
  -\frac{1}{(u+i(N- \frac{1}{2}))^2+\xi^2}+\frac{\alpha_2^2-\xi^2}
  {\beta_2^2-\beta_1^2}\\
 &\times (\frac{1}{(u+i(N- \frac{1}{2}))^2+\beta_1^2}-\frac{1}{(u+i(N- 
  \frac{1}{2}))^2+\beta_2^2})] (u+i(N- \frac{1}{2}))^3 \tanh \pi u du. 
\endaligned$$
Breaking up this integral into contributions for $|u| \leq N$ and
$|u|>N$, then its value is $O(N^{-2})$ for fixed $\xi$, $\alpha_2$,
$\beta_1$ and $\beta_2$. Thus
$$\aligned
  W_0(\xi)
=&2 \sum_{k=0}^{\infty} k^3 \left[\frac{1}{\alpha_2^2-k^2}-
  \frac{1}{\xi^2-k^2}+\frac{\alpha_2^2-\xi^2}{\beta_2^2-\beta_1^2}
  \left(\frac{1}{\beta_1^2-k^2}-\frac{1}{\beta_2^2-k^2}\right)\right]\\  
 &+\pi (\xi^2 \cot \pi \xi-\alpha_2^2 \cot \pi \alpha_2) 
  -\pi \frac{\alpha_2^2-\xi^2}{\beta_2^2-\beta_1^2} \left(
  \beta_1^2 \cot \pi \beta_1-\beta_2^2 \cot \pi \beta_2 \right).
\endaligned\tag 3.13$$
Similarly, 
$$\aligned
  W_1(\xi)
=&2 \sum_{k=0}^{\infty} (k^3-k) \left[\frac{1}{\alpha_2^2-k^2}-
  \frac{1}{\xi^2-k^2}+\frac{\alpha_2^2-\xi^2}{\beta_2^2-\beta_1^2}
  \left(\frac{1}{\beta_1^2-k^2}-\frac{1}{\beta_2^2-k^2}\right)\right]\\  
 &+\pi [(\xi^2-1) \cot \pi \xi-(\alpha_2^2-1) \cot \pi \alpha_2]\\ 
 &-\pi \frac{\alpha_2^2-\xi^2}{\beta_2^2-\beta_1^2} [(\beta_1^2-1) 
  \cot \pi \beta_1-(\beta_2^2-1) \cot \pi \beta_2].
\endaligned\tag 3.14$$
$$\aligned
  W_{\frac{1}{2}}(\xi)
=&\sum_{k=0}^{\infty} k(k+1)(2k+1) [\frac{1}{\alpha_2^2-(k+
  \frac{1}{2})^{2}}-\frac{1}{\xi^2-(k+\frac{1}{2})^{2}}+
  \frac{\alpha_2^2-\xi^2}{\beta_2^2-\beta_1^2}\\
 &\times \left(\frac{1}{\beta_1^2-(k+\frac{1}{2})^2}-
  \frac{1}{\beta_2^2-(k+\frac{1}{2})^2}\right)]+
  \pi [(\frac{1}{4}-\xi^2) \tan \pi \xi-(\frac{1}{4}-\alpha_2^2)\\
 &\times \tan \pi \alpha_2]-\pi \frac{\alpha_2^2-\xi^2}{\beta_2^2-
  \beta_1^2} [(\frac{1}{4}-\beta_1^2) \tan \pi \beta_1-(\frac{1}{4}
  -\beta_2^2) \tan \pi \beta_2].
\endaligned\tag 3.15$$

{\smc Proposition 3.7}. {\it The following formulas hold:
$$\aligned
 &\frac{1}{s-1} \frac{Z_k^{\prime}(s)}{Z_k(s)}=\frac{4}{3} 
  c_2(M) W_{k}(s-1)+\\
 &+\frac{1}{\alpha_2} \frac{Z_k^{\prime}(\alpha_2+1)}{Z_k(\alpha_2+1)}
  -\frac{(s-1)^2-\alpha_2^2}{\beta_1^2-\beta_2^2} \left[\frac{1}
  {\beta_2} \frac{Z_k^{\prime}(\beta_2+1)}{Z_k(\beta_2+1)}-\frac{1}
  {\beta_1} \frac{Z_k^{\prime}(\beta_1+1)}{Z_k(\beta_1+1)}\right]+\\
 &-2 \sum_{n=0}^{\infty} \left[\frac{1}{r_n^2+\alpha_2^2}-
  \frac{1}{r_n^2+(s-1)^2}+\frac{(s-1)^2-\alpha_2^2}{\beta_1^2-
  \beta_2^2} \left(\frac{1}{r_n^2+\beta_1^2}-\frac{1}{r_n^2+
  \beta_2^2}\right)\right],
\endaligned\tag 3.16$$
where $k=0$, $\frac{1}{2}$, $1$ and $\text{Re}(s)>1$.}

{\smc Theorem 3.8}. {\it The functional equations of zeta 
functions of weight $k$ are given by
$$\aligned
 &Z_k(s)=Z_k(2-s) \exp \left[\frac{8}{3} \pi c_2(M) 
  \int_{0}^{s-1} v(v^2-k^2) \cot \pi v dv \right], \quad 
  k=0, 1;\\
 &Z_k(s)=Z_k(2-s) \exp \left[\frac{8}{3} \pi c_2(M)
  \int_{0}^{s-1} v(k^2-v^2) \tan \pi v dv \right], \quad
  k=\frac{1}{2}, 
\endaligned\tag 3.17$$
where $M=\Gamma \backslash {\frak S}_{2}$ is a compact complex
algebraic surface with $c_1^2=3 c_2$.}

{\it Proof}. When $k=0$, by Proposition 3.7, we have
$$\frac{1}{s-1} \left[\frac{Z_0^{\prime}(s)}{Z_0(s)}+
  \frac{Z_0^{\prime}(2-s)}{Z_0(2-s)}\right]=\frac{8}{3} \pi
  c_2(M) (s-1)^2 \cot \pi (s-1).$$
An elementary calculation gives
$$\log \frac{Z_0(s)}{Z_0(2-s)}=\frac{8}{3} \pi c_2(M)
  \int_{0}^{s-1} v^3 \cot \pi v dv.$$
The cases of $k=1$ and $k=\frac{1}{2}$ can be proved in the same way.
\flushpar
$\qquad \qquad \qquad \qquad \qquad \qquad \qquad \qquad \qquad
 \qquad \qquad \qquad \qquad \qquad \qquad \qquad \qquad \qquad
 \quad \square$

{\smc Corollary 3.9}. {\it The functional equations of 
zeta functions of weight $k$ can be expressed as
$$Z_0(s)=Z_0(2-s) [\sin \pi (s-1)]^{\frac{8}{3} c_2(M) 
  (s-1)^3} \exp \left[-8 c_2(M) \int_{0}^{s-1} v^2
  \log \sin \pi v dv \right],\tag 3.18$$ 
$$\aligned
  Z_1(s)
=&Z_1(2-s) [\sin \pi (s-1)]^{\frac{8}{3} c_2(M) s(s-1)(s-2)}\\
 &\times \exp \left[-8 c_2(M) \int_{0}^{s-1} (v^2-\frac{1}{3}) 
  \log \sin \pi v dv \right],
\endaligned\tag 3.19$$
and
$$\aligned
  Z_{\frac{1}{2}}(s)
=&Z_{\frac{1}{2}}(2-s) [\cos \pi (s-1)]^{\frac{8}{3} c_2(M)
  (s- \frac{1}{2})(s-1)(s- \frac{3}{2})}\\
 &\times \exp \left[-8 c_2(M) \int_{0}^{s-1} (v^2-\frac{1}{12})
  \log \cos \pi v dv \right].   
\endaligned\tag 3.20$$}

{\it Proof}. If $k=0$, it is obtained by the following integral	   
$$\int_{0}^{s-1} \pi v^3 \cot \pi v dv=(s-1)^3 \log \sin \pi (s-1)
 -3 \int_{0}^{s-1} v^2 \log \sin \pi v dv.$$
The other two cases can be proved in the same way.
\flushpar
$\qquad \qquad \qquad \qquad \qquad \qquad \qquad \qquad \qquad
 \qquad \qquad \qquad \qquad \qquad \qquad \qquad \qquad \qquad
 \quad \square$

\vskip 0.4 cm
3.7. {\it The analogue of the Riemann hypothesis}
\vskip 0.2 cm

{\smc Theorem 3.10}. {\it Let 
$$Z_{k}(s)=\prod_{\{(\gamma_1, \gamma_2)\}} \prod_{j=0}^{\infty} 
  \prod_{l=0}^{\infty} \prod_{n=0}^{\infty} 
  \left[1-N(\gamma_1)^{-s-j-\frac{l+n}{2}}\right]^{\frac{3}
  {\text{ord}(\gamma_2)} \sum_{q=1}^{\text{ord}(\gamma_2)} 
  \exp([2k+3(l-n)]iq \theta(\gamma_2))}$$
where $\text{Re}(s)>2$ and $k=0, \frac{1}{2}, 1$. Then	   
\roster
\item $Z_k(s)$ are actually entire functions.
\item The identity of Proposition 3.7 holds for all $s$.
\item $Z_k(s)$ $(k=0, 1)$ have trivial zeros $s=-m$, $m \geq 0$, with
      multiplicity $\frac{8}{3} c_2(M) (m+1)^{3}$;
	  $Z_k(s)$ $(k=\frac{1}{2})$ has trivial zeros $s=-m+\frac{1}{2}$,
	  $m \geq 0$, with multiplicity $\frac{4}{3} c_2(M) m(m+1)(2m+1)$. 
\item The nontrivial zeros of $Z_k(s)$ are located at $1 \pm i r_n$,
      i.e., the analogue of the Riemann hypothesis is true.
\endroster}

{\it Proof}. By (3.16), $Z_k^{\prime}(s)/Z_k(s)$ continues 
meromorphically to the whole complex plane ${\Bbb C}$.
By the calculation on the $r_n$ sum in (3.12), we know that
the $r_n$ sum contributes principal parts $(s-s_n)^{-1}$,
$(s- \widetilde{s}_{n})^{-1}$ for each $n \geq 0$. This implies 
that the analogue of the Riemann hypothesis is true.
By $\pi \cot \pi z=\frac{1}{z}+\sum_{n=1}^{\infty} \left(
\frac{1}{z-n}+\frac{1}{z+n}\right)$ and (3.13), when the
weight $k=0$, we have
$$\aligned
 &-2 \sum_{m=0}^{\infty} \frac{m^3}{(s-1)^2-m^2}+\pi (s-1)^2
  \cot \pi (s-1)\\
=&-\sum_{m=1}^{\infty} m^2 \left(\frac{1}{s-1-m}-\frac{1}{s-1+m}
  \right)+(s-1)+\\
 &+(s-1)^2 \sum_{m=1}^{\infty} \left(\frac{1}{s-1-m}
  +\frac{1}{s-1+m}\right)\\
=&(s-1)+\sum_{m=1}^{\infty} (s-1+m)+\sum_{m=1}^{\infty}
  \frac{(s-1)^2+m^2}{s-1+m}.
\endaligned$$
Here $\sum_{m=1}^{\infty} \frac{(s-1)^2+m^2}{s-1+m}=
\sum_{m=0}^{\infty} \frac{(s-1)^2+(m+1)^2}{s+m}$. The
$m$-sum contributes 
$$(m+1) \cdot \frac{4}{3} c_2(M) \cdot 2(m+1)^2 (s+m)^{-1}
 =\frac{8}{3} c_2(M) \cdot (m+1)^{3} (s+m)^{-1}.$$
Hence, $Z_{0}(s)$ has trivial zeros $s=-m$ with multiplicity
$\frac{8}{3} c_2(M) (m+1)^{3}$.
  
  The cases of $k=1$ and $k=\frac{1}{2}$ can be proved in the 
same way.
\flushpar
$\qquad \qquad \qquad \qquad \qquad \qquad \qquad \qquad \qquad
 \qquad \qquad \qquad \qquad \qquad \qquad \qquad \qquad \qquad
 \quad \square$
 		  
\vskip 0.5 cm
\centerline{\bf Appendix}
\vskip 0.5 cm

  In this appendix, we calculate the following four integrals.
$$\int_{0}^{\infty} \frac{\cosh u \sin ru-r \cos ru \sinh u}
  {\sinh^3 u} du, \quad
  \int_{0}^{\infty} \frac{\cosh u \sin ru-r \cos ru \sinh u}
  {\sinh \frac{u}{2} \sinh^2 u} du,$$
$$\int_{0}^{\infty} \frac{\sin ru}{\sinh \frac{u}{2} \cosh^3
  \frac{u}{2}} du, \quad
  \int_{0}^{\infty} \frac{\cos ru}{\cosh^4 \frac{u}{2}} du.$$
  
  For the first integral, we denote 
$$G(r)=\int_{0}^{\infty} \frac{\cosh u \sin ru-r \sinh u \cos ru}
       {\sinh^3 u} du.$$
Let 
$$f(z)=\frac{\cosh z-i r \sinh z}{\sinh^3 z} e^{i r z},$$
$\text{Im} f(u)=\frac{\cosh u \sin ru-r \sinh u \cos ru}
 {\sinh^3 u}$ for $u \in {\Bbb R}$. 
In fact,
$f(z)=\frac{4[(1-ir)e^{4z}+(1+ir)e^{2z}]}{(e^{2z}-1)^{3}} e^{irz}$.
$z=0$ and $z=\pi i$ are two poles of order three of $f$. Thus,
$$\text{Res}(f, 0)=\frac{1}{2} \lim_{z \to 0} \frac{d^2}{dz^2}
  [z^3 f(z)], \quad
  \text{Res}(f, \pi i)=\frac{1}{2} \lim_{z \to \pi i}
  \frac{d^2}{dz^2}[(z-\pi i)^3 f(z)].$$
Set 
$$g_1(z)=z^3 f(z)=h(z) p(z)^3, \quad 
  g_2(z)=(z- \pi i)^3 f(z)=h(z) q(z)^3,$$
where
$$h(z)=4[(1-ir)e^{4z}+(1+ir)e^{2z}] e^{irz}, \quad
  p(z)=\frac{z}{e^{2z}-1}, \quad q(z)=\frac{z-\pi i}{e^{2z}-1}.$$
Then
$$\lim_{z \to 0} p(z)=\frac{1}{2}, \quad
  \lim_{z \to 0} p^{\prime}(z)=-\frac{1}{2}, \quad 
  \lim_{z \to 0} p^{\prime \prime}(z)=\frac{1}{3}.$$
$$\lim_{z \to \pi i} q(z)=\frac{1}{2}, \quad 
  \lim_{z \to \pi i} q^{\prime}(z)=-\frac{1}{2}, \quad
  \lim_{z \to \pi i} q^{\prime \prime}(z)=\frac{1}{3}.$$
$$h(0)=8, \quad h^{\prime}(0)=24, \quad 
  h^{\prime \prime}(0)=80+8 r^2.$$
$$h(\pi i)=8 e^{-\pi r}, \quad h^{\prime}(\pi i)=24 e^{-\pi r}, 
  \quad h^{\prime \prime}(\pi i)=(80+8 r^2) e^{-\pi r}.$$
Moreover,
$$\lim_{z \to 0} g_{1}^{\prime \prime}(z)=
  \frac{1}{8} h^{\prime \prime}(0)-\frac{3}{4} h^{\prime}(0)
  +h(0)=r^2,$$
$$\lim_{z \to \pi i} g_2^{\prime \prime}(z)=
  \frac{1}{8} h^{\prime \prime}(\pi i)-\frac{3}{4}
  h^{\prime}(\pi i)+h(\pi i)=r^2 e^{-\pi r}.$$
Hence, 
$$\text{Res}(f, 0)=\frac{1}{2} r^2, \quad  
  \text{Res}(f, \pi i)=\frac{1}{2} r^2 e^{-\pi r}.$$

  Now, we consider the following contour integral: 
I-II-III-$\gamma_2$-IV-$\gamma_1$-I, where
I is an ordered line segment from $[\varepsilon, 0]$ to $[R, 0]$, 
II is an ordered line segment from $[R, 0]$ to $[R, \pi]$,
III is an ordered line segment from $[R, \pi]$ to $[\varepsilon, \pi]$,
$\gamma_2$ is an ordered arc from $[\varepsilon, \pi]$ 
to $[0, \pi-\varepsilon]$, IV is an ordered line segment from  
$[0, \pi-\varepsilon]$ to $[0, \varepsilon]$, and
$\gamma_1$ is an ordered arc from $[0, \varepsilon]$ to 
$[\varepsilon, 0]$.

  In II, $\int_{II} f(z) dz \to 0$ as $R \to \infty$.
In III, $z=x+\pi i$. Hence,
$$\int_{III} f(z) dz=-e^{-\pi r} \int_{\varepsilon}^{R}
  \frac{4[(1-ir) e^{4x}+(1+ir) e^{2x}]}{(e^{2x}-1)^{3}}
  e^{irx} dx.$$
For $\gamma_1$, $z=0$ is the center of this arc. 
$\int_{\gamma_1} f(z) dz \to -\frac{1}{4} \cdot 2 \pi i 
 \text{Res}(f, 0)=-\frac{\pi i}{4} r^2$.
For $\gamma_2$, $z=\pi i$ is the center of this arc.
$\int_{\gamma_2} f(z) dz \to -\frac{1}{4} \cdot 2 \pi i 
 \text{Res}(f, \pi i)=-\frac{\pi i}{4} r^2 e^{-\pi r}$.
In IV,
$$\aligned
 \int_{IV} f(z) dz
=&-\int_{\varepsilon}^{\pi-\varepsilon} \frac{4[(1-ir) e^{4iy}+
  (1+ir) e^{2iy}]}{(e^{2iy}-1)^{3}} e^{ir \cdot iy} d(i y)\\
=&\int_{\varepsilon}^{\pi-\varepsilon} \frac{\cos y+r \sin y}
  {\sin^3 y} e^{-r y} dy.
\endaligned$$
By Cauchy theorem, we have
$(\int_{I}+\int_{II}+\int_{III}+\int_{\gamma_2}+\int_{IV}+
 \int_{\gamma_1}) f(z) dz=0$, i.e.,
$$\aligned
 &(1-e^{-\pi r}) \int_{\varepsilon}^{R} \frac{4[(1-ir) e^{4x}+
  (1+ir) e^{2x}]}{(e^{2x}-1)^{3}} e^{irx} dx- \frac{\pi i}{4} r^2\\
-&\frac{\pi i}{4} r^2 e^{-\pi r}+\int_{\varepsilon}^{\pi-\varepsilon}
  \frac{\cos y+r \sin y}{\sin^3 y} e^{-r y} dy+o(1)=0.
\endaligned$$
Take imaginary parts and let $\varepsilon \to 0$, $R \to \infty$, 
we get 
$$G(r)=\frac{\pi}{4} r^2 \frac{1+e^{-\pi r}}{1-e^{-\pi r}}.
       \tag A.1$$ 

  For the second integral, we denote
$$H(r)=\int_{0}^{\infty} \frac{\cosh u \sin ru-r \sinh u \cos ru}
  {\sinh \frac{u}{2} \sinh^2 u} du.$$
Let
$$f(z)=\frac{\cosh z-ir \sinh z}{\sinh \frac{z}{2} \sinh^2 z}
       e^{irz}.$$
By the same method as above, we have
$$H(r)=\frac{\pi}{2} \left(r^2+\frac{1}{4}\right) \tanh \pi r.
  \tag A.2$$

  For the third integral
$$\int_{0}^{\infty} \frac{\sin ru}{\sinh \frac{u}{2} \cosh^3
  \frac{u}{2}} du,$$ 
let
$$f(z)=\frac{16 e^{irz} e^{2z}}{(e^z-1)(e^z+1)^3}.$$

  Now, we consider the following contour integral: 
I-II-III-$\gamma_2$-IV-$\gamma_3$-V-$\gamma_1$-I, where
I is an ordered line segment from $[\varepsilon, 0]$ to $[R, 0]$, 
II is an ordered line segment from $[R, 0]$ to $[R, 2 \pi]$,
III is an ordered line segment from $[R, 2 \pi]$ to 
$[\varepsilon, 2 \pi]$,
$\gamma_2$ is an ordered arc from $[\varepsilon, 2 \pi]$ 
to $[0, 2 \pi-\varepsilon]$, IV is an ordered line segment from  
$[0, 2 \pi-\varepsilon]$ to $[0, \pi+\varepsilon]$, $\gamma_3$ is
an ordered arc from $[0, \pi+\varepsilon]$ to $[0, \pi-\varepsilon]$,
V is an ordered line segment from $[0, \pi-\varepsilon]$ to
$[0, \varepsilon]$, and $\gamma_1$ is an ordered arc from 
$[0, \varepsilon]$ to $[\varepsilon, 0]$.
For II: $\int_{II} f(z) dz \to 0$ as $R \to \infty$.
For III: $z=x+2 \pi i$.
$\int_{III} f(z) dz=-e^{-2 \pi r} \int_{\varepsilon}^{R}
 \frac{16 e^{irx} e^{2x}}{(e^x-1)(e^x+1)^3} dx$.
For I: $z=x$.
$\int_{I} f(z) dz=\int_{\varepsilon}^{R} \frac{16 e^{irx} e^{2x}}
 {(e^x-1)(e^x+1)^3} dx$.
$z=0$ is a pole of order one, $\text{Res}(f, 0)=2$.
$\int_{\gamma_{1}} f(z) dz \to -\frac{1}{4} \cdot 2 \pi i
 \text{Res}(f, 0)=-\pi i$.
$z=2 \pi i$ is a pole of order one, 
$\text{Res}(f, 2 \pi i)=2 e^{-2 \pi r}$.
$\int_{\gamma_{2}} f(z) dz \to -\frac{1}{4} \cdot 2 \pi i
 \text{Res}(f, 2 \pi i)=-\pi i e^{-2 \pi r}$.
$z=\pi i$ is a pole of order three,
$$\text{Res}(f, \pi i)=\frac{1}{2} \lim_{z \to \pi i} \frac{d^2}{dz^2}
 \left[(z- \pi i)^3 \frac{16 e^{irz} e^{2z}}{(e^z-1)(e^z+1)^{3}}
 \right]=-2(2 r^2+1) e^{-\pi r}.$$
$$\int_{\gamma_{3}} f(z) dz \to -\frac{1}{2} \cdot 2 \pi i
  \text{Res}(f, \pi i)=2 \pi i (2 r^2+1) e^{-\pi r}.$$
For IV: $z=iy$,
$$\int_{IV} f(z) dz=\int_{2 \pi-\varepsilon}^{\pi+\varepsilon}
  \frac{16 e^{ir (iy)} e^{2 iy}}{(e^{iy}-1)(e^{iy}+1)^{3}} d(iy)
 =-\int_{\pi+\varepsilon}^{2 \pi-\varepsilon} \frac{e^{-ry}}
  {\sin \frac{y}{2} \cos^3 \frac{y}{2}} dy.$$
For V: $z=iy$, 
$\int_{V} f(z) dz=-\int_{\varepsilon}^{\pi-\varepsilon}
 \frac{e^{-ry}}{\sin \frac{y}{2} \cos^3 \frac{y}{2}} dy$.
By the residue theorem, we have
$$\aligned
 &(1-e^{-2 \pi r}) \int_{\varepsilon}^{R} \frac{16 e^{irx} e^{2x}}
  {(e^x-1)(e^x+1)^3} dx- \left(\int_{\pi+\varepsilon}^{2 \pi-
  \varepsilon}+\int_{\varepsilon}^{\pi-\varepsilon}\right)
  \frac{e^{-ry}}{\sin \frac{y}{2} \cos^3 \frac{y}{2}} dy\\
 &-\pi i- \pi i e^{-2 \pi r}+2 \pi i (2r^2+1) e^{-\pi r} 
  +o(1)=0.
\endaligned$$
Take imaginary parts and let $\varepsilon \to 0$ and $R \to \infty$,
we obtain
$$\int_{0}^{\infty} \frac{\sin r x}{\sinh \frac{x}{2} \cosh^3
  \frac{x}{2}} dx=\frac{\pi}{\tanh \pi r}-\frac{\pi (2r^2+1)}
  {\sinh \pi r}.\tag A.3$$

  For the fourth integral
$$\int_{0}^{\infty} \frac{\cos ru}{\cosh^4 \frac{u}{2}} du,$$ 
let
$$f(z)=\frac{16 e^{irz} e^{2z}}{(e^z+1)^{4}}.$$

  Now, we consider the following contour integral: 
I-II-III-IV-$\gamma$-V-I, where
I is an ordered line segment from $[0, 0]$ to $[R, 0]$, 
II is an ordered line segment from $[R, 0]$ to $[R, 2 \pi]$,
III is an ordered line segment from $[R, 2 \pi]$ to 
$[0, 2 \pi]$, IV is an ordered line segment from $[0, 2 \pi]$
to $[0, \pi+\varepsilon]$, $\gamma$ is an ordered arc from 
$[0, \pi+\varepsilon]$ to $[0, \pi-\varepsilon]$, V is an ordered 
line segment from $[0, \pi-\varepsilon]$ to $[0, 0]$.
For I: 
$\int_{I} f(z) dz=\int_{0}^{R} \frac{16 e^{irx} e^{2x}}{(e^x+1)^4} dx$.
For II: $\int_{II} f(z) dz \to 0$ as $R \to \infty$.
For III: $z=x+2 \pi i$. 
$$\int_{III} f(z) dz=-\int_{0}^{R} \frac{16 e^{ir(x+2 \pi i)} 
  e^{2(x+2 \pi i)}}{(e^{x+2 \pi i}+1)^{4}} dx
 =-e^{-2 \pi r} \int_{0}^{R} \frac{16 e^{irx} e^{2x}}{(e^x+1)^4} dx.$$
For IV: $z=iy$. 
$$\int_{IV} f(z) dz=\int_{2 \pi}^{\pi+\varepsilon} \frac{16
  e^{ir(iy)} e^{2iy}}{(e^{iy}+1)^{4}} d(iy)
 =-i \int_{\pi+\varepsilon}^{2 \pi} \frac{e^{-ry}}{\cos^4 
  \frac{y}{2}} dy.$$
For V: $z=iy$. $\int_{V} f(z) dz=-i \int_{0}^{\pi-\varepsilon}
   \frac{e^{-ry}}{\cos^4 \frac{y}{2}} dy$.  
$z=\pi i$ is a pole of order four.
$$\text{Res}(f, \pi i)=\frac{1}{3!} \lim_{z \to \pi i}
  \frac{d^3}{dz^3} \left[(z-\pi i)^4 \frac{16 e^{irz} e^{2z}}
  {(e^z+1)^4}\right]=-\frac{8}{3} ir(r^2+1) e^{-\pi r}.$$
$$\int_{\gamma} f(z) dz \to -\frac{1}{2} \cdot 2 \pi i
  \text{Res}(f, \pi i)=-\frac{8}{3} \pi r(r^2+1) e^{-\pi r}.$$
The residue theorem tell us that
$$(1-e^{-2 \pi r}) \int_{0}^{R} \frac{16 e^{irx} e^{2x}}{(e^x+1)^4}
  dx-i \left(\int_{0}^{\pi-\varepsilon}+\int_{\pi+\varepsilon}^{2 \pi}
  \right) \frac{e^{-ry}}{\cos^4 \frac{y}{2}} dy- \frac{8}{3} \pi r
  (r^2+1) e^{-\pi r}+o(1)=0.$$
Take real parts and let $R \to \infty$, we obtain
$$\int_{0}^{\infty} \frac{\cos rx}{\cosh^4 \frac{x}{2}} dx=
  \frac{4 \pi r (r^2+1)}{3 \sinh \pi r}.\tag A.4$$

\vskip 1.0 cm
{\smc Department of Mathematics, Peking University}

{\smc Beijing 100871, P. R. China}

{\it E-mail address}: yanglei\@math.pku.edu.cn
\vskip 1.0 cm
\Refs

\item{[Ar]} {\smc J. Arthur}, The Selberg trace formula for groups
             of $F$-rank one, Ann. of Math. {\bf 100} (1974), 326-385.

\item{[APS]} {\smc M. F. Atiyah, V. K. Patodi and I. M. Singer},			 
			 Spectral asymmetry and Riemannian geometry I, 
			 Proc. Camb. Phil. Soc. {\bf 77} (1975), 43-69.
			 			 
\item{[BM]} {\smc D. Barbasch and H. Moscovici}, $L^2$-index and the
            Selberg trace formula, J. Funct. Anal. {\bf 53} (1983),
			151-201.

\item{[BPV]} {\smc W. Barth, C. Peters and A. Van de Ven},
             {\it Compact Complex Surfaces}, Ergeb. Math.
			 Grenzgeb. (3) {\bf 4}, Springer-Verlag, 1984.

\item{[Bo]} {\smc A. Borel}, Compact Clifford-Klein forms of
             symmetric spaces, Topology {\bf 2} (1963), 111-122.
	  
\item{[CV]} {\smc P. Cartier and A. Voros}, Une nouvelle
             interpr\'{e}tation de la formule des traces de Selberg,
 	         in: {\it The Grothendieck Festschrift, Vol. II}, 
			 Progress in Math. {\bf 87}, Birkh\"{a}user, 1990,
			 1-67.

\item{[DP]} {\smc E. D'Hoker and D. H. Phong}, On determinants of
            Laplacians on Riemann surfaces, Comm. Math. Phys.
			{\bf 104} (1986), 537-545.

\item{[DL]} {\smc M. Duflo and J. P. Labesse}, Sur la formule des
            traces de Selberg, Ann. Scient. \'{E}c. Norm. Sup.
			(4) {\bf 4} (1971), 193-284.

\item{[El]} {\smc J. Elstrodt}, Die Resolvente zum Eigenwertproblem
            der automorphen Formen in der hyperbolischen Ebene, I,
			Math. Ann. {\bf 203} (1973), 295-330; II, Math. Z. 
			{\bf 132} (1973), 99-134; III, Math. Ann. {\bf 208} 
			(1974), 99-132.

\item{[ElGM]} {\smc J. Elstrodt, F. Grunewald and J. Mennicke}, {\it
            Groups Acting on Hyperbolic Space, Harmonic Analysis and
			Number Theory}, Springer Monographs in Mathematics,
			Springer-Verlag, 1998. 

\item{[Ga]} {\smc R. Gangolli}, Zeta functions of Selberg's type
            for compact space forms of symmetric spaces of rank
			one, Illinois J. Math. {\bf 21} (1977), 1-42.
			
\item{[GW]} {\smc R. Gangolli and G. Warner}, Zeta functions of
            Selberg's type for some noncompact quotients of
			symmetric spaces of rank one, Nagoya Math. J. 
			{\bf 78} (1980), 1-44.			

\item{[Gar]} {\smc P. R. Garabedian}, {\it Partial Differential
             Equations}, AMS Chelsea Publishing, 1998.
			 
\item{[G]} {\smc G. Giraud}, Sur certaines fonctions automorphes
           de deux variables, Ann. \'{E}c. Norm. (3) {\bf 38}
           (1921), 43-164.

\item{[Go]} {\smc W. M. Goldman}, {\it Complex Hyperbolic 
            Geometry}, Oxford Mathematical Monographs,
			Oxford, 1999. 

\item{[H]} {\smc D. A. Hejhal}, {\it The Selberg Trace Formula for
           $PSL(2, {\Bbb R})$}, Vol. {\bf 1}, Lecture Notes in Math.
		   {\bf 548}, Springer-Verlag, 1976.

\item{[HP]} {\smc S. Hersonsky and F. Paulin}, On the volumes of
            complex hyperbolic manifolds, Duke Math. J. {\bf 84}
            (1996), 719-737.

\item{[Hi1]} {\smc F. Hirzebruch}, Automorphe Formen und der Satz
             von Riemann-Roch, Symp. Intern. Top. Alg. 1956, 129-144,
			 Universidad de Mexico 1958.

\item{[Hi2]} {\smc F. Hirzebruch}, {\it Topological Methods in
             Algebraic Geometry}, 3rd ed., Grundlehren Math. Wiss. 
			 {\bf 131}, Springer-Verlag, 1978.

\item{[Ku]} {\smc T. Kubota}, {\it Elementary Theory of Eisenstein 
            Series}. Wiley, New York, 1973.        

\item{[Kr]} {\smc N. Kurokawa}, Special values of Selberg zeta 
            functions, Contemp. Math. {\bf 83} (1989), 133-150.

\item{[Ma]} {\smc H. Maass}, Die Differentialgleichungen in der 
            Theorie der elliptischen Modulfunktionen, Math. Ann.
			{\bf 125} (1953), 235-263. 

\item{[M]} {\smc J. J. Millson}, Closed geodesics and the 
            $\eta$-invariant, Ann. of Math. {\bf 108} (1978), 1-39.

\item{[R]} {\smc M. S. Raghunathan}, {\it Discrete Subgroups of
            Lie Groups}, Ergeb. Math. Grenzgeb. {\bf 68},
      		Springer-Verlag, 1972.

\item{[RS]} {\smc D. B. Ray and I. M. Singer}, Analytic torsion for
            complex manifolds, Ann. of Math. {\bf 98} (1973), 154-177.

\item{[Ro]} {\smc W. Roelcke}, Das Eigenwertproblem der automorphen
             Formen in der hyperbolischen Ebene, I, Math. Ann. 
			 {\bf 167} (1966), 292-337; II, Math. Ann. {\bf 168} 
			 (1967), 261-324.

\item{[Sa]} {\smc P. Sarnak}, Determinants of Laplacians, Comm.
            Math. Phys. {\bf 110} (1987), 113-120.

\item{[Se]} {\smc A. Selberg}, Harmonic analysis and discontinuous
            groups in weakly symmetric Riemannian spaces with 
			applications to Dirichlet series, J. Indian Math. Soc.
			B. {\bf 20} (1956), 47-87.
			
\item{[S]} {\smc I. M. Singer}, Eigenvalues of the Laplacian and
            invariants of manifolds, Proc. Internat. Cong. Math.
			Vancouver, B. C. 1974, 187-200.

\item{[St]} {\smc M. Stern}, $L^2$-index theorems on locally symmetric
            spaces, Invent. Math. {\bf 96} (1989), 231-282.

\item{[V]} {\smc A. Voros}, Spectral functions, special functions
             and the Selberg zeta function, Comm. Math. Phys.
			 {\bf 110} (1987), 439-465.

\item{[Y]} {\smc L. Yang}, Product formulas on a unitary group in
            three variables, math.NT/0104259.

\item{[Yau]} {\smc S. T. Yau}, On Calabi's conjecture and some new 
            results in algebraic geometry, Proc. Nat. Acad. Sci.
			USA {\bf 74} (1977), 1798-1799.

\item{[Yo]} {\smc K. Yosida}, {\it Functional Analysis}, 6th ed.,
            Grundlehren Math. Wiss. {\bf 123}, Springer-Verlag, 
			1980. 
			
\endRefs
\end{document}